\documentclass[12pt]{amsart}

\usepackage{amssymb,amsmath,amsthm,amsfonts,amscd,bbm}
\usepackage{enumitem}
\usepackage{pdflscape}
\usepackage[paper=portrait,pagesize]{typearea}
\usepackage{caption}
\usepackage{bm}

\usepackage{geometry} 

\usepackage{ifpdf}
\ifpdf
\usepackage[pdftex]{graphicx}
\else
\usepackage[dvips]{graphicx}
\fi
\usepackage[all]{xy}
\usepackage{hyperref}
\usepackage{xcolor}
\usepackage{graphicx}
\usepackage[outdir=./]{epstopdf}
\usepackage{enumitem}
\usepackage{tensor}
\usepackage{tikz-cd}
\usepackage{multicol}
\usepackage{mathtools}
\usepackage{tocvsec2}
\usepackage{bbm}
\usepackage{longtable}
\usepackage{fullpage}

\newcommand{\doi}[1]{\href{https://dx.doi.org/#1}{{\small DOI:#1}}}

\newcommand{\googlebooks}[1]{(preview at \href{https://books.google.com/books?id=#1}{google books})}

\newcommand{\numdam}[1]{}
\usepackage{mathrsfs}

\usepackage{mathrsfs}
\DeclareMathAlphabet{\mathpzc}{OT1}{pzc}{m}{it}

\def\semicolon{;}
\def\applytolist#1{
    \expandafter\def\csname multi#1\endcsname##1{
        \def\multiack{##1}\ifx\multiack\semicolon
            \def\next{\relax}
        \else
            \csname #1\endcsname{##1}
            \def\next{\csname multi#1\endcsname}
        \fi
        \next}
    \csname multi#1\endcsname}

\def\calc#1{\expandafter\def\csname c#1\endcsname{{\mathcal #1}}}
\applytolist{calc}QWERTYUIOPLKJHGFDSAZXCVBNM;
\def\bbc#1{\expandafter\def\csname bb#1\endcsname{{\mathbb #1}}}
\applytolist{bbc}QWERTYUIOPLKJHGFDSAZXCVBNM;
\def\bfc#1{\expandafter\def\csname bf#1\endcsname{{\mathbf #1}}}
\applytolist{bfc}QWERTYUIOPLKJHGFDSAZXCVBNM;
\def\sfc#1{\expandafter\def\csname s#1\endcsname{{\sf #1}}}
\applytolist{sfc}QWERTYUIOPLKJHGFDSAZXCVBNM;
\def\fc#1{\expandafter\def\csname f#1\endcsname{{\mathfrak #1}}}
\applytolist{fc}QWERTYUIOPLKJHGFDSAZXCVBNM;

\usepackage{tikz}
\usepackage{tikz-cd}

\makeatletter
\def\fixtikzforbreqn#1#2{%
  \protected\edef#1{\noexpand\ifmmode\mathchar\the\mathcode`#2 \noexpand\else#2\noexpand\fi}%
}
\fixtikzforbreqn\tikz@nonactivesemicolon;
\fixtikzforbreqn\tikz@nonactivecolon:
\fixtikzforbreqn\tikz@nonactivebar|
\fixtikzforbreqn\tikz@nonactiveexlmark!
\makeatother

\usepackage{breqn}   

\usetikzlibrary{arrows,backgrounds,patterns.meta}
\usetikzlibrary{positioning,shadings,cd}
\usetikzlibrary{shapes}
\usetikzlibrary{backgrounds}
\usetikzlibrary{decorations,decorations.pathreplacing,decorations.markings,decorations.pathmorphing}
\usetikzlibrary{fit,calc,through}
\usetikzlibrary{external}
\usetikzlibrary{arrows}
\tikzset{vertex/.style = {shape=circle,draw,fill=black,inner sep=0pt,minimum size=5pt}}
\tikzset{edge/.style = {->,> = latex', bend right}}
\tikzset{
	super thick/.style={line width=3pt}
}
\tikzset{
    quadruple/.style args={[#1] in [#2] in [#3] in [#4]}{
        #1,preaction={preaction={preaction={draw,#4},draw,#3}, draw,#2}
    }
}
\tikzstyle{shaded}=[fill=red!10!blue!20!gray!30!white]
\tikzstyle{unshaded}=[fill=white]
\tikzstyle{empty box}=[circle, draw, thick, fill=white, opaque, inner sep=2mm]
\tikzstyle{annular}=[scale=.7, inner sep=1mm, baseline]
\tikzstyle{rectangular}=[scale=.75, inner sep=1mm, baseline=-.1cm]
\tikzstyle{mid>}=[decoration={markings, mark=at position 0.5 with {\arrow{>}}}, postaction={decorate}]
\tikzstyle{mid<}=[decoration={markings, mark=at position 0.5 with {\arrow{<}}}, postaction={decorate}]
\tikzstyle{over}=[double, draw=white, super thick, double=]
\tikzstyle{snake}=[decorate, decoration={snake, segment length=1mm, amplitude=.3mm}]
\tikzstyle{saw}=[decorate, decoration={saw, segment length=.7mm, amplitude=.25mm}]

\tikzstyle{coupon}=[draw, very thick, rectangle, rounded corners=5pt]
\tikzset{Rightarrow/.style={double equal sign distance,>={Implies},->},
triplecd/.style={-,preaction={draw,Rightarrow}},
quadruplecd/.style={preaction={draw,Rightarrow,
shorten >=0pt
},
shorten >=1pt,
-,double,double
distance=0.2pt}}
\tikzset{
    tripleline/.style args={[#1] in [#2] in [#3]}{
        #1,preaction={preaction={draw,#3},draw,#2}
    }
}
\tikzstyle{triple}=[tripleline={[line width=.15mm,black] in
      [line width=.7mm,white] in
      [line width=1mm,black]}] 
\tikzset{
    quadrupleline/.style args={[#1] in [#2] in [#3] in [#4]}{
        #1,preaction={preaction={preaction={draw,#4},draw,#3}, draw,#2}
    }
}
\tikzstyle{quadruple}=[quadrupleline={[line width=.3mm,white] in
      [line width=.6mm,black] in
      [line width=1.2mm,white] in
      [line width=1.5mm,black]}]
 
\newcommand{\roundRbox}[4]{
    \draw[very thick, #1] ($#2$) circle (#3);   
	\node at ($#2$) {#4};
}

\newcommand{\roundNbox}[6]{
	\draw[rounded corners=5pt, very thick, #1] ($#2+(-#3,-#3)+(-#4,0)$) rectangle ($#2+(#3,#3)+(#5,0)$);
	\coordinate (ZZa) at ($#2+(-#4,0)$);
	\coordinate (ZZb) at ($#2+(#5,0)$);
	\node at ($1/2*(ZZa)+1/2*(ZZb)$) {#6};
}

\newcommand{\tikzmath}[2][]
     {\vcenter{\hbox{\begin{tikzpicture}[#1]#2
                     \end{tikzpicture}}}
     }

\makeatletter
\newcommand{\xMapsto}[2][]{\ext@arrow 0599{\Mapstofill@}{#1}{#2}}
\def\Mapstofill@{\arrowfill@{\Mapstochar\Relbar}\Relbar\Rightarrow}
\makeatother

\theoremstyle{plain}
\newtheorem{thm}{Theorem}[section]
\newtheorem*{thm*}{Theorem}
\newtheorem{thmalpha}{Theorem}

\newtheorem{cor}[thm]{Corollary}
\newtheorem{coralpha}[thmalpha]{Corollary}
\newtheorem*{cor*}{Corollary}

\newtheorem*{conj*}{Conjecture}
\newtheorem{lem}[thm]{Lemma}
\newtheorem*{lem*}{Lemma}

\newtheorem{prop}[thm]{Proposition}

\newtheorem*{quest*}{Question}
\newtheorem*{claim*}{Claim}

\theoremstyle{definition}
\newtheorem{defn}[thm]{Definition}

\newtheorem{ex}[thm]{Example}
\newtheorem{sub-ex}[thm]{Sub-Example}
\newtheorem{counter-ex}[thm]{Counter-Example}
\newtheorem{rem}[thm]{Remark}
\newtheorem*{rem*}{Remark}
\newtheorem{remark}[thm]{Remark}

\usepackage{xcolor}
\definecolor{dark-red}{rgb}{0.7,0.25,0.25}
\definecolor{dark-blue}{rgb}{0.15,0.15,0.55}
\definecolor{medium-blue}{rgb}{0,0,.8}
\definecolor{DarkGreen}{RGB}{0,150,0}
\definecolor{rho}{named}{red}
\hypersetup{
   colorlinks, linkcolor={purple},
   citecolor={medium-blue}, urlcolor={medium-blue}
}

\newcommand{\XColor}{red} 
\newcommand{\YColor}{blue}

\newcommand{\HColor}{red} 
\newcommand{\rColor}{gray!20} 

\newcommand{\id}{\operatorname{id}}

\newcommand{\Irr}{\operatorname{Irr}}

\newcommand{\op}{\operatorname{op}}

\newcommand{\Hom}{\operatorname{Hom}}

\newcommand{\End}{\operatorname{End}}

\newcommand{\Aut}{\operatorname{Aut}}

\newcommand{\Ad}{\operatorname{Ad}}

\newcommand{\FPdim}{\operatorname{FPdim}}

\newcommand{\fd}{{\sf fd}}
\newcommand{\Mod}{{\sf Mod}}

\newcommand{\Stab}{{\sf Stab}}

\newcommand{\fssMod}{{\sf Mod_{f}}}
\newcommand{\Bim}{{\sf Bim}}

\newcommand{\Fun}{{\mathsf{Fun}}}

\newcommand{\Hilb}{\mathsf{Hilb}}

\newcommand{\BigHilb}{\mathsf{BigHilb}}

\newcommand{\rCorr}{{\mathsf{C^*Alg}}}

\newcommand{\IL}{{\sf IndLimC^*Alg}}

\newcommand{\AbM}{{\sf AbMonoid}}

\newcommand{\PAbG}{{\sf PAbGrp}}
\newcommand{\QSys}{\mathsf{QSys}}
\newcommand{\Fib}{\mathsf{Fib}}
\newcommand{\triv}{\mathsf{triv}}
\newcommand{\sign}{\mathsf{sign}}
\newcommand{\MPO}{\mathsf{MPO}}
\newcommand{\InDec}{\mathsf{InDec}}
\newcommand{\ExStd}{\mathsf{ExStd}}

\def\altdb{\vadjust{\vbox to 0pt{\vss\hbox{\kern \hsize
\quad{\dbend}}\kern\baselineskip\kern-10pt}}}

\newcommand{\noshow}[1]{}
\renewcommand{\MR}[1]{}

\title{K-theoretic classification of inductive limit actions of fusion categories on AF-algebras}
\author{Quan Chen, Roberto Hern\'{a}ndez Palomares, and Corey Jones}

\newcommand{\Contact}{{
  \bigskip
  \footnotesize

  Quan Chen, \textsc{Department of Mathematics, Vanderbilt University}\par\nopagebreak
  \textit{E-mail address}: \texttt{quan.chen@vanderbilt.edu}

  \medskip

  Roberto Hern\'{a}ndez Palomares, \textsc{University of Waterloo}\par\nopagebreak
  \textit{E-mail address}: \texttt{robertohp.math@gmail.com}

  \medskip

  Corey Jones, \textsc{Department of Mathematics, North Carolina State University}\par\nopagebreak
  \textit{E-mail address}: \texttt{cmjones6@ncsu.edu}
}}

\begin{document}

\maketitle

\begin{abstract}
We introduce a K-theoretic invariant for actions of unitary fusion categories on unital $\rm C^*$-algebras. We show that for inductive limits of finite dimensional actions of fusion categories on AF-algebras, this is a complete invariant. In particular, this gives a complete invariant for inductive limit actions of finite groups on unital AF-algebras. We apply our results to obtain a classification of finite depth, strongly AF-inclusions of unital AF-algebras.
\end{abstract}

\tableofcontents

\section{Introduction}
Fusion categories are algebraic objects generalizing both finite groups and their representation categories.
They have come to play an important role as generalized symmetries in many areas of mathematics and physics, from subfactors to quantum field theory \cite{MR3166042,MR3242743, MR3308880}. 
Underlying many of these applications is the notion of an \textit{action} of a fusion category $\cC$ on an algebra $A$, which is encoded by a unitary tensor functor $F:\cC\rightarrow \Bim(A)$\ \footnote{$\Bim(A)$ denotes some flavor of monoidal category of bimodules, depending on context.} 
\cite{MR3663592, MR3948170,MR4139893,2010.01072, MR4328058, 2105.05587, MR4236062, MR4419534, AranoRokhlin}. 
Recently, there has been an interest in fusion categorical symmetries of spin chains motivated by the search for exotic conformal field theories \cite{MR3719546, PhysRevLett.128.231602, PhysRevLett.128.231603, PhysRevB.96.125104, 2205.15243}.
Considering these concepts from an operator algebraic point of view leads directly to the study of inductive limits of fusion category actions on finite dimensional $\rm C^*$-algebras, 
which we call AF-actions (see Sections \ref{sec:AFActions}, \ref{sec:MPO} and Definition \ref{defn:AFAction}).
In this paper, we provide a K-theoretic classification of these actions up to an equivalence relation generalizing cocycle conjugacy for group actions.
In particular, our K-theoretic invariant is a complete invariant for inductive limit actions of finite groups on (unital) AF-algebras up to cocycle conjugacy.

AF-algebras are inductive limits of finite dimensional $\rm C^*$-algebras. Elliott showed that the triple $(K_{0}(A), K^{+}_{0}(A), [A])$\
\footnote{In this paper, we use the finitely generated projective Hilbert module picture of $K_{0}$ rather than the projection picture, and thus we denote the distinguished element by the Hilbert $A$-module $[A]$ instead of the customary $[1_{A}]$.} is a complete invariant for isomorphism classes of AF-algebras \cite{MR397420}. This seminal result initiated the classification program for simple nuclear $\rm C^*$-algebras and their symmetries, which has seen tremendous progress in recent years (see \cite{2015arXiv150100135G, 2015arXiv150703437E, MR3583354, MR2663778,MR2827854, MR4207309, 2022arXiv220504933G} and the references therein). 
There are several results in the literature which provide K-theoretic classifications of AF-actions of finite groups in certain cases \cite{MR1193924, MR769758, MR1078515, MR1414065, https://doi.org/10.48550/arxiv.1304.0813, MR4206896}. These results suggest the possibility of a full K-theoretic classification of AF-actions of finite groups and more generally, of fusion categories. The main goal of this paper is to provide such a classification.

Given an action $\cC\overset{F}{\curvearrowright} A$ of a unitary fusion category $\cC$ on a (unital) $\rm C^*$-algebra $A$, we introduce a K-theoretic invariant $(\widehat{F}, [A])$ (see Definition \ref{pointed invariant}). 
As a piece of data, $\widehat{F}$ consists of a finite list of (pre)-ordered abelian groups, 
and a finite family of positive homomorphisms between them satisfying compositional relations determined by $\cC$, while $[A]\in K_{0}(A)$ is the usual distinguished element ($K_{0}(A)$ is one of the groups in our list). In Section \ref{subsec:groupactionex}, we provide an elementary description of this invariant in the special case of ordinary group actions. Before giving a more detailed explanation of the invariant, we record our main result:

\begin{thmalpha}\label{thmalp:main}
Suppose $\cC\overset{F}{\curvearrowright} A$ and $\cC\overset{G}{\curvearrowright} B$ are AF-actions. Then $\left(\widehat{F}, [A]\right)\cong \left(\widehat{G}, [B]\right)$ if and only if $\cC\overset{F}{\curvearrowright} A$ is equivalent to $\cC\overset{G}{\curvearrowright} B$. 
\end{thmalpha}

To obtain the invariant $(\widehat{F}, [A])$, the idea is straightforward: we first find the simplest invariant that completely classifies finite dimensional actions of fusion categories and their equivariant morphisms, which are well understood algebraically via a unitary version of Ostrik's Theorem (see Remark \ref{rem:unitaryostrik}).
We then show that the property of being a complete invariant passes to inductive limits, following a slightly modified version of Elliott's intertwining argument. This approach is closely related to the ideas of \cite{MR4207309}.

Q-systems are generalizations of finite dimensional $\rm C^*$-algebras, which are \textit{internal} to a fusion category $\cC$ (see Section \ref{sec: Q-systems}). Q-systems assemble into a category $[\QSys(\cC)]$, whose objects are Q-systems and morphisms are equivalence classes of bimodule objects in $\cC$. Given an action of a fusion category $\cC\overset{F}{\curvearrowright} A$, a Q-system can be viewed as precisely the data required to form a generalized ``crossed product". For a Q-system $P$, this crossed product is denoted $|P|_{F}$ and is called the \textit{realization} of $P$ (see Section \ref{subsec:realization}). In previous work, we showed that realization extends to a functor from $[\QSys(\cC)]$ to the category of (unital) $\rm C^*$-algebras and isomorphism classes of (right finite) correspondences \cite{MR4419534}. Composing this with the usual $K_{0}$ functor from $\rm C^*$-algebras to the category of (pre)-ordered abelian groups ($\PAbG$) yields a functor $\widehat{F}:[\QSys(\cC)]\rightarrow \PAbG$. The natural isomorphism class of $\widehat{F}$ is the main ingredient in our invariant. The other piece of data, $[A]$, simply remembers the class of the trivial Hilbert module in $K_{0}(A)=\widehat{F}(1_{\cC})$.

While the definitions are phrased in the language of categories and functors, we emphasize that $[\QSys(\cC)]$ is fundamentally combinatorial and $\widehat{F}$ is K-theoretic in nature. In the next section, we give an explicit, category-free description of the invariant $(\widehat{F}, [A])$ in the case that $\cC=\Hilb(G)$ for $G$ a finite group. Nevertheless, computing $[\QSys(\cC)]$ from the data of $\cC$ itself can be fairly involved in general.
Even for a finite group $G$, fully expressing the invariant can be complicated (see Section \ref{subsec:groupactionex} and Appendix \ref{App:Examples}).  
In other cases the invariant is very straightforward. A fusion category $\cC$ is called \textit{torsion-free} \cite{MR3941472} if there is only one Morita class of indecomposable Q-system. Examples include the Fibonacci category $\Fib$ and $\Hilb(G, \omega)$, where $[\omega]\in H^{3}(G, \text{U}(1))$ does not trivialize on any subgroup. In this case, the invariant is just the usual invariant of the $\rm C^*$-algebra $(K_{0}(A), K^{+}_{0}(A), [A])$ viewed as right module over the fusion ring of $\cC$. Then an immediate consequence of our classification result is the following: 

\begin{coralpha} If $\cC$ is a torsion-free fusion category and $A$ is an AF $\rm C^*$-algebra with \newline $(K_{0}(A), K^{+}_{0}(A))\subset (\bbR, \bbR^{+})$, then there is at most one equivalence class of AF-action of $\cC$ on $A$.
\end{coralpha}
\noindent As a consequence, we are able to construct and classify all equivalence classes of AF-actions of $\Fib$ on the AF-algebra $A_\varphi,$ with $K_0(A_\varphi)=\bbZ+\varphi\bbZ\subset \bbR,$ where $\varphi$ is the golden ratio (Corollary \ref{cor:FibFib}).

We can also consider our results from the perspective of subfactor theory. The study of finite index inclusion of $\rm{II}_{1}$ factors began with the groundbreaking work of V. Jones \cite{MR0696688}. Hyperfinite von Neumann algebras are a natural $\rm W^*$-analogue of AF-algebras. Popa's deep results show that (strongly) amenable inclusions of hyperfinite $\rm{II}_{1}$ factors are uniquely determined by their standard invariant \cite{MR1055708, MR1278111}. This reduces the classification of finite depth inclusions of hyperfinite $\rm{II}_{1}$ factors to the algebraic problem of classifying their standard invariants (\cite{MR1334479, MR996454, math.QA/9909027}), which in turn can be rephrased as a problem about unitary fusion categories and their Q-systems \cite{MR1966524,MR3948170,MR1424954}. The algebraic perspective (particularly the framework of V. Jones' planar algebras) has led to remarkable progress in the small index subfactor classification program (see \cite{MR3166042, 1509.00038} and references therein). The successful reduction of the a-priori analytic problem of subfactor classification to an algebraic one in the hyperfinite, finite depth case motivates the search for similar types of results in the setting of inclusions of AF-algebras.

The philosophy espoused in \cite{MR4419534} is that classifying ``finite index" inclusions of algebras is closely related to classifying tensor category actions, and thus our results should have something to say about this problem. An analogue of a hyperfinite subfactor in the $\rm C^*$-setting is an inclusion $(A\subseteq B, E_{A})$ of unital AF-algebras with trivial centers and faithful conditional expectation $E_{A}:B\rightarrow A$. 
Such an inclusion is called \textit{finite depth} if $B$, viewed as an $A$-$A$ correspondence, tensor generates a unitary fusion category. 
The standard invariant of the inclusion is the pair $(\cC, Q)$ where $\cC:=\langle {}_{A} B_{A}\rangle\subseteq \Bim(A)$ is the fusion category generated by $B$, and $Q={}_{A} B_{A}$ is the Q-system in $\cC$ corresponding to the algebra structure on $B$ as in \cite{MR4419534}. 
An inclusion $(A\subseteq B, E_{A})$ is called \textit{strongly AF} if the inclusion $F:\cC\rightarrow \Bim(A)$ is an AF-action. 
The we define the \textit{Extended Standard Invariant} $\ExStd(A\subseteq B, E_{A}):=(\cC, Q, \widehat{F}, [A])$, where $(\cC, Q)$ is the standard invariant as above, and $(\widehat{F}, [A])$ is the invariant described earlier.

The following classification result, in the spirit of Popa's theorem \cite{MR1055708} for finite depth hyperfinite subfactors, is a corollary of our main Theorem.

\begin{thmalpha}\label{thmalp:ExStd&AFaction}
Let $(A\subseteq B, E_A)$ and $(C\subseteq D, E_C)$ be finite depth strongly AF-inclusions of AF $\rm C^*$-algebras. Then $(A\subseteq B, E_A)$ is equivalent to $(C\subseteq D, E_C)$ if and only if $\ExStd(A\subseteq B, E_A)\cong \ExStd(C\subseteq D, E_C)$.
\end{thmalpha}

There is one final perspective on our results that we would like to highlight. An action of a fusion category $\cC$ on a $\rm C^*$-algebra can be viewed as a special case of a $\rm C^*$-algebra object internal to $\cC$ in the sense of \cite{MR3687214}. There are natural notions of nuclearity for a $\rm C^*$-algebra object, which can be translated into properties of the action. This suggests one may be able to develop a classification program for simple nuclear $\rm C^*$-algebras internal to $\cC$, extending the standard version. We believe our work serves as a starting point for this endeavor.

\subsection{The invariant for group actions}\label{subsec:groupactionex}

Our paper makes heavy use of the language of category and 2-category theory. This is largely unavoidable, since this is the natural setting for fusion category actions on $\rm C^*$-algebras. However, for the reader that is not familiar with these concepts, we want to give a description of our invariant in the case of group actions which can be stated in elementary terms.

In this section let $G$ be a finite group.

\begin{defn}
A twisted subgroup is a pair $(H,\mu)$, where $H\le G$ and $\mu\in Z^{2}(H, \text{U}(1))$. Two twisted subgroups $(H,\mu)$ and $(K,\nu)$ are equivalent if there exists a $g\in G$ such that $gHg^{-1}=K$ and $[\mu^{g}]=[\nu]\in H^{2}(K, \text{U}(1))$.
\end{defn}
\noindent In the above and throughout this paper, $Z^{2}(G,\text{U}(1))$ refers to the group of \textit{normalized} cocycles, i.e. $\mu(1,h)=\mu(h,1)=1$ for all $h\in H$. Also $\mu^{g}\in Z^{2}(K, \text{U}(1))$ is given by $\mu^{g}(k_{1},k_{2}):=\mu(g^{-1}k_1g, g^{-1}k_{2}g)$.

Now we recall the \textit{twisted crossed product} construction. Given an action of $G$ on a $\rm C^*$-algebra $A$, then for any twisted subgroup $(H,\mu)$ we have the crossed product algebra
$$A\rtimes_{\mu} H=\sum_{h\in H} u_{h} a_{h},$$ 
where $a_{h}\in A$, and for $h\in H$, the $u_{h}$ are unitaries satisfying $au_{h}=u_{h} h^{-1}(a)$ and $u_{h_1}u_{h_{2}}=\mu(h_{1},h_{2})u_{h_1h_2}$. This is a $\rm C^*$-algebra (note $u^{*}_{h}=\mu(h,h^{-1})u_{h^{-1}}$).
There is a standard conditional expectation $E_A: A\rtimes_\mu H\to A$ given by $u_h a_h\mapsto \delta_{h=1}a_h$.

For the trivial action of $G$ on $\bbC$, the crossed product construction applied to a twisted subgroup $(H,\mu)$ yields the twisted group algebras, denoted $\bbC_{\mu}[H]$. These are finite dimensional $\rm C^*$-algebras.

\begin{defn}\label{def:twisted bimodule} For two twisted subgroups $(H,\mu)$ and $(K,\nu)$ then a $(H, \mu)-(K, \nu)$ bimodule consists of a finite dimensional Hilbert space $\bbC_{\mu}[H]-\bbC_{\nu}[K]$ bimodule $X$, equipped with a $G$-grading $X=\bigoplus_{g\in G} X_{g}$ such that $u_{h}\rhd X_{g}\lhd u_{k}\subseteq X_{hgk}$.

\end{defn}

We have the following proposition:

\begin{prop} Let $G\rightarrow \Aut(A)$ be an action. Then an $(H,\mu)-(K,\nu)$ bimodule defines a bifinite $A\rtimes_{\mu} H-A\rtimes_{\nu} K$ correspondence (\cite[Definition 2.12]{MR4419534}), and thus induces a positive map $(K_{0}(A\rtimes_{\mu} H),K^{+}_{0}(A\rtimes_{\mu} H) ) \rightarrow (K_{0}(A\rtimes_{\nu} K), K^{+}_{0}(A\rtimes_{\nu} K))$.
\end{prop}
\begin{proof}

This follows from our general theory of realization \cite{MR4419534}. The connections between realization and crossed product is elaborated on in Section \ref{subsec:realization}). We will sketch the construction explicitly here for the convenience of the reader.

Let $X$ be an $(H,\mu)-(K,\nu)$ bimodule. Then define the vector space 
$$A\rtimes X := \bigoplus_{g\in G} X_{g}\otimes_{\bbC} A $$
This has a left $A\rtimes_{\mu} H$ action via
$$(u_h a) \rhd (x_{g}\otimes b):= (u_{h}\rhd x_{g})\otimes g^{-1}(a)b,$$
and right $A\rtimes_{\nu} K$ action via
$$(x_{g}\otimes b) \lhd (u_{k} a):= (x_{g}\lhd u_{k})\otimes k^{-1}(b)a$$

It's not hard to show that this bimodule can be made into an $A\rtimes_{\mu} H-A\rtimes_{\nu} K$ bifinite correspondence.

Now, given any finitely generated projective Hilbert $A\rtimes_{\mu} H$-module $Y$, then $Y\boxtimes_{A\rtimes_{\mu} H} A\rtimes X$ is a finitely generated projective right Hilbert $A\rtimes_{\nu} K$ module. Thus
$- \boxtimes_{A\rtimes_{\mu} H} A\rtimes X$ induces a positive map on $K_{0}$ as desired.
\end{proof}

\begin{defn}(The invariant).
Given an action $G\rightarrow \Aut(A)$, we define a cocycle conjugacy invariant which consists of 
\begin{enumerate}
    \item 
    For each conjugacy class of twisted subgroup $(H, \mu)$, the pre-ordered abelian group $(K_{0}(A\rtimes_{\mu} H), K^{+}_{0}(A\rtimes_{\mu} H))$
    \item
    For each isomorphism class $(X,\eta)$ of irreducible $(H,\mu)-(K,\nu)$ bimodule, the induced map $- \boxtimes_{A\rtimes_{\mu} H} A\rtimes X: K_{0}(A\rtimes_{\mu} H)\rightarrow K_{0}(A\rtimes_{\nu} K)$
    \item
    For the trivial subgroup class $(1_{G}, 1)$, the order unit $[A]\in K_{0}(A)$.
\end{enumerate}
\end{defn}

By isomorphism of our invariant, we mean isomorphism pre-order isomorphisms of each of the $K_{0}$ groups that intertwines the actions of the equivariant correspondences. It turns out our invariant is an invariant up to cocycle conjugacy. A corollary of our main result is the following:

\begin{cor}
Let $A$ and $B$ be unital AF-algebras, $G\overset{\alpha}{\curvearrowright} A$ and $G\overset{\beta}{\curvearrowright} B$ be actions which are both inductive limits of finite dimensional actions.  
Then $\alpha$ is cocycle conjugate to $\beta$ if and only if the invariants described above are isomorphic.
\end{cor}

We prove all these results in the much more general context of AF-actions of fusion categories.

\subsection*{Outline}

In Section \ref{sec:PreLim}, we give an overview of some of the required preliminaries concerning the relevant categorical algebra. In Section \ref{sec:Action}, we formally define actions of fusion categories, define our invariant, and prove some basic properties concerning it. In Section \ref{sec:ILActions}, we introduce AF-actions, and describe their connection to symmetries of spin chains. We then prove the main Theorem \ref{thmalp:main} and Theorem \ref{thmalp:ExStd&AFaction}. Finally, in Section \ref{sec:Examples} we consider several examples which demonstrate how to compute the invariant in practice.

\subsection*{Acknowledgements}
The authors would like to thank Sam Evington, Zhengwei Liu, Sergio Gir\'on Pacheco, and David Penneys for helpful comments and enlightening conversations. Corey Jones was supported by NSF DMS-2100531, Roberto Hern\'andez Palomares was partially supported by NSF grant DMS-2000331, and Quan Chen was partially supported by NSF DMS-1654159.

\section{Preliminaries}
\label{sec:PreLim}

In this section we give a brief overview of some basic notions and establish some working notation. For more in-depth descriptions of these topics we refer the reader to \cite[\S 2]{MR4419534}

\begin{defn}
A \emph{$\rm C^*$-category} $\cC$ is a $\bbC$-linear category satisfying the following conditions:
\begin{enumerate}
\item for each pair of objects $a,b\in\cC$, there is a conjugate-linear involution
$*:\cC(a\to b)\to \cC(b\to a)$ such that $(f\circ g)^*=g^*\circ f^*$,
\item there is a Banach norm on morphisms such that $\|f\|^2 = \|f\circ f^*\|= \|f^*\circ f\|,$ for all $f\in \cC(a\to b)$,
\item for each $f\in \cC(a\to b)$, there exists $g\in \cC(a\to a)$ such that $f^*\circ f = g^*\circ g$.
\item
$\cC$ is Cauchy complete: there exists unitary direct sums, and every projection splits (see \cite[Assumption 2.7]{MR4419534}).
\end{enumerate}
\end{defn}

\begin{defn}
A \emph{$\rm C^*$ 2-category} is a 2-category such that every every $\hom$ 1-category is equipped with the structure of a $\rm C^*$-category, and the horizontal composition of 2-morphisms is compatible with the $*$-operations.
\end{defn}

For an unpacked version of the above definition, we refer the reader to \cite[Definition 2.2]{MR4419534}. 
We refer the reader to \cite{MR4261588} for general references on the theory of 2-categories.

\medskip

\noindent \textbf{The $\rm C^*$ 2-category of unital $\rm C^*$-algebras}. The most important $\rm C^*$ 2-category in this paper, $\rCorr,$ consists of all unital $\rm C^*$-algebras $A,B, \hdots$ as objects, 1-morphisms ${}_AX_B\in\rCorr(A\to B)$ are right $\rm C^*$-correspondences which are \textit{finitely generated projective as right Hilbert $\rm C^*$-modules}, and 2-morphisms $f\in \rCorr({}_AX_B\Rightarrow {}_AY_B)$ are adjointable $A$-$B$ bimodular maps. The horizontal composition of 1 and 2-morphisms is the relative tensor product of correspondences ${}_{A}X\boxtimes_{B} Z_{C}$, while vertical composition of 2-morphisms is just ordinary composition. 

\begin{defn}
If $A$ is a unital $\rm C^*$-algebra, we define $\Bim(A)$ to be the full tensor subcategory of dualizable objects in $\rCorr(A\rightarrow A)$, and $\Mod(A)$ the $\rm C^*$-category $\rCorr(\bbC\rightarrow A)$.
\end{defn}

Note that any dualizable right correspondence admits the structure of a left finitely generated projective Hilbert $A$-module (see \cite{MR2085108}).

\medskip

\noindent \textbf{Unitary tensor categories}. A \emph{unitary tensor category} $\cC$ is a semisimple rigid $\rm C^*$-tensor category with simple unit object $1_{\cC}$. Rigidity means that every object has a unitary dual (see \cite[Definition 2.9]{MR4419534} and references therein), and simplicity of the tensor unit is the requirement that $\End_\cC(1_\cC)\cong\bbC$. If the isomorphism classes of simple objects in a unitary tensor category form a finite set, we say the category is a \emph{unitary fusion category}. 

\begin{rem}\label{opremark}
There is a somewhat confusing issue here that will appear frequently in this paper. Given a $\rm C^*$ 2-category with one object $\ast$, we can build a $\rm C^*$-tensor category by taking $\cC=\Hom(\ast\to  \ast)$, and defining $X\otimes Y$ to be the composition of the $1$-morphisms $X$ and $Y$. The problem is that there is an ambiguity since we can compose these 1-morphisms in either order. We will take the convention that our 1-morphisms in a 2-category are always composed from left to right. These seems natural, but there are some unfortunate consequences to this convention. In particular if we take the 2-category of categories, functors and natural transformations, and we consider $\End(\cM)$ as a tensor category, then by our conventions $F\otimes G:=G\circ F$. We urge the reader to keep this in mind.
\end{rem}

\medskip

\noindent \textbf{Unitary tensor functors}. For unitary tensor categories $\cC$ and $\cD,$ a \emph{unitary tensor functor} is a triple $(F, F^1, F^2),$ where $F:\cC\to \cD$ is a $*$-functor, 
$F^1:F(1_{\cC})\to 1_{\cD}$ is a chosen unitary isomorphism, and 
$F^2 =\left\{F^2_{a,b}:F(a)\otimes F(b)\to F(a\otimes b)\right\}_{a,b\in\cC}$ is a natural unitary isomorphism satisfying the standard coherence axioms c.f. \cite[Chapter 2.4]{MR3242743}. Notice we don't require our unitary tensor functors be fully-faithful. In this paper, we will assume our monoidal functors are \textit{strictly unital}, i.e. $F(1_{\cC})=1_{\cD}$ and $F^1=\id_{1_{\cD}}.$

\medskip

\noindent \textbf{Module categories}. Given a unitary tensor category $\cC,$ a (right) unitary $\cC$-\emph{module category} is a $\rm C^*$-category $\cM$ together with a unitary tensor functor $\cC\to\End(\cM).$ 
The target category $\End(\cM)$ is the strict tensor category of all $*$-endofunctors of $\cM$, with the convention that $F\otimes G:=G\circ F$ (see Remark \ref{opremark}). Given a $\cC$-module category $\rho:\cC\rightarrow \End(\cM)$, we will often use the notation $m\lhd c:=\rho(c)(m)$.

Given module-categories $\cM_\cC$ and $\cN_\cC,$ a \emph{$\cC$-module functor} $F:\cM\to\cN$ is a $*$-functor together with a family unitary natural transformations $\{u_{c}: F(\cdot \lhd c)\cong F(\cdot)\lhd c\ |\ c\in \cC\}$ satisfying the obvious coherences (see \cite[Section 7.2]{MR3242743})

The collection of $\cC$-module categories, $\cC$-module functors, and module functor natural transformations assemble into a $\rm C^*$ 2-category.

\begin{defn}If $\cC$ is a unitary fusion category, we denote by $\Mod(\cC)$ the $\rm C^*$ 2-category of right $\cC$-module categories, $\cC$-module functors and $\cC$-module natural transformations. $\fssMod(\cC)$ is the full 2-subcategory whose objects are finitely semisimple as $\rm C^*$-categories.
\end{defn}

\subsection{Q-systems}\label{sec: Q-systems}
\ \\
Recall that if $\cC$ is a unitary tensor category, a \textit{Q-system} is a unitary version of a separable Frobenius algebra object. In particular a (unital) Q-system consists of a pair $(Q, m, i),$ where $Q\in \cC$ and $m\in \cC(Q\otimes Q\to Q)$ and $i\in \cC(1_{\cC}\to Q)$ satisfying the following equations (with associators/unitors suppressed):

\begin{enumerate}[label=(Q\arabic*)]
    \item 
    (Associativity) $m \circ (\id_{Q}\otimes m)=m\circ (m\otimes \id_{Q})$.
    \item
    (Unital) $m\circ (i\otimes \id_{Q})=\id_{Q}=m\circ (\id_{Q}\otimes i)$.
    \item
    (Frobenius) $(m\otimes \id_{Q})\circ (\id_{Q}\otimes m^{*})=m^{*}\circ m= (\id_{Q}\otimes m)\circ (m^{*}\otimes \id_{Q})$.
    \item
    (Separability) $m\circ m^{*}=\id_{Q}$.
\end{enumerate}

If $P, Q$ are Q-systems, a unitary \textit{$P-Q$ bimodule} is an object $X\in \cC$ together with maps $\lambda\in \cC(P\otimes X\to X)$ and $\rho\in \cC(X\otimes Q\to X)$ satisfying the following equations
\begin{enumerate}[label=(B\arabic*)]
    \item 
    (Left $P$-module) $\lambda\circ (\id_P\otimes\lambda) = \lambda\circ (m_P\otimes \id_X)$.
    \item
    (Right $Q$-module) $\rho\circ (\rho\otimes \id_Q) = \rho\circ (\id_X\otimes m_Q)$.
    \item
    (Bimodule) $\rho \circ (\lambda \otimes \id_{Q})= \lambda\circ (\id_{P}\otimes \rho)$.
    \item
    (Unital) $\lambda\circ (i_{P}\otimes \id_{X})=\id_{X}=\rho\circ (\id_{X}\otimes i_{Q})$.
    \item 
    (Frobenius) $(\id_Q\otimes\lambda)\circ (m^*\otimes \id_X) = \lambda^*\circ \lambda$, $(\rho\otimes \id_Q)\circ (\id_X\otimes m^*) =\rho^*\circ \rho$.
    \item 
    (Separability) $\lambda\circ \lambda^* =\id_X =\rho\circ \rho^*$
\end{enumerate}

A \textit{bimodule intertwiner} between two $P-Q$ bimodules $X, Y$ is a morphism $f\in\cC(X\Rightarrow Y)$ intertwining the $P$ and $Q$ actions; i.e. 
$f\circ \lambda = \lambda\circ (\id_P\otimes f),$ and $ f\circ\rho= \rho\circ(f\otimes \id_Q)$.

\begin{defn}
Given a unitary fusion category $\cC$, we denote by $\QSys(\cC)$ the collection of all Q-systems in $\cC$, called \emph{the Q-system completion of $\cC$}, which assembles naturally into a finite rank semisimple $\rm C^*$ 2-category \cite{1812.11933}. 
Namely, $\QSys(\cC)$ consists of the following data:     
\begin{itemize}
    \item Objects in $\QSys(\cC)$ are all Q-systems $P, Q\in\cC,$
    \item 1-morphisms ${}_PX_Q\in\QSys(\cC)(P\to Q)$ are $P$-$Q$ bimodules. Note that \textbf{the source is the left index} and \textbf{the target is the right index}.
    \item 2-morphisms $f\in\QSys(\cC)(X\Rightarrow Y)$ are $P$-$Q$ bimodule intertwiners.
\end{itemize}
For a detailed description of all the higher compositions and coherences, see \cite[\S 3]{MR4419534}.  
\end{defn}
Conceptually, the Q-system completion of a unitary tensor category can be thought of a ``higher unitary idempotent completion" of $\cC,$ initially thought of as a $\rm C^*$ 2-category with one object \cite{1905.09566,  2021arXiv210809872V}. There are finitely many isomorphism classes (i.e. Morita classes) of indecomposable objects (Q-systems).
Recall a Q-system is \textit{indecomposable} if it is not isomorphic to a direct sum of two sub Q-systems.
A Q-system $Q$ is \textit{connected} if $\cC(1_\cC\to Q)$ is one dimensional.  
Then every indecomposable Q-system is Morita equivalent to a connected one (this follows from the unitary version of Ostrik's theorem below), and if $\cC$ is fusion there are finitely many of these up to algebra isomorphism \cite{MR1932664}.

\begin{rem}[Unitary Ostrik's Theorem]\label{rem:unitaryostrik} $\QSys(\cC)$ is (1-contravariantly) equivalent to the $\rm C^*$ 2-category $\fssMod(\cC)$ of finitely semisimple unitary right $\cC$-module categories, via a unitary version of Ostrik's theorem \cite{MR3933035, 2021arXiv210809872V}. The 2-functor witnessing this equivalence takes a Q-system $Q$ and assigns it the right $\cC$-module category ${}_{Q}\cC$ of left $Q$-modules internal to $\cC$.
Namely, the equivalence is given by:
\begin{align*}\begin{split}
    \Psi:\QSys(\cC)&\to \fssMod(\cC)\\
        Q&\mapsto {}_Q\cC\\
        {}_P X_Q &\mapsto \left({}_P X\otimes_Q-:{}_Q\cC \to{}_P\cC  \right)\\
        \left(f:{}_P X_Q\Rightarrow {}_P Y_Q\right)&\mapsto f\otimes_Q \id_{{}_Q\cC}.
\end{split}\end{align*}
\end{rem}

\bigskip
In this paper we will be interested in $[\QSys(\cC)],$ the \textit{decategorification} of the $\rm C^*$ 2-category $\QSys(\cC)$, 
and much more so when $\cC$ is a unitary fusion category. 

\begin{defn}\label{defn:[QSys]}
We define $[\QSys(\cC)]$ as the (1-)category whose objects 
$P,Q,R,\hdots$ are Q-systems in $\cC$, 
and whose morphisms $[X],...,[Z]\in[\QSys(\cC)](P\to Q)$ are isomorphism classes of $P$-$Q$ bimodules in $\QSys(\cC)$. 

For $X\in\QSys(\cC)(P\to Q),$ and $Z\in\QSys(\cC)(Q\to R),$ define
 $[X]\circ [Z]:= [X\otimes_Q Z].$ 
\end{defn}
 
We remark that the additive structure; i.e. $[X]+[Y]:= [X\oplus Y]$, in $[\QSys(\cC)]$ is inherited from the $\rm C^*$ 2-category structure of $\cC.$ Furthermore, composition in $[\QSys(\cC)]$ is bilinear, namely $$([X]+ [Y])\circ [Z]=[X]\circ [Z]+[Y]\circ [Z]\ \text{and}\ [Z]\circ ([X]+ [Y])=[Z]\circ [X]+[Z]\circ [Y].$$
Thus $[\QSys(\cC)]$ is enriched in the symmetric monoidal category of abelian monoids. We also note that every $Q$ system is isomorphic to the direct sum of indecomposables, so it suffices to describe the indecomposables and the morphsims between them when considering $[\QSys(\cC)]$ (see Section \ref{sec:Examples} for a more precise statement).

\begin{ex}{(Pointed fusion categories)}
\label{ex:PointedCatas}
We will describe here the objects and morphisms of the category $[\QSys(\Hilb(G,\omega))]$, where $G$ is a finite group and $\omega\in Z^{3}(G, \text{U}(1))$. We have the following facts about classifying Q-systems (cf \cite{MR3933035}):

\begin{itemize}
\item 
Connected Q-systems in $\Hilb(G,\omega)$ correspond to pairs $(H, \mu)$, where $H\le G$ and \newline$\mu:H\times H\rightarrow \text{U}(1)$ is a trivialization of $\omega|_{H}$. 
\item
Two $Q$-systems $(H, \mu)$ and $(K, \nu)$ are isomorphic if $H=K$ and $\mu$ differs from $\nu$ by a coboundary. 
\item
$(H, \mu)$ and $(K, \nu)$ are Morita equivalent (i.e. isomorphic in $[\QSys]$) if there exists a $g\in G$ such that $(H^{g},\mu^{g})$ is isomorphic to $(K, \nu),$ where $(\cdot)^{g}$ denotes the conjugation automorphism.
\end{itemize}

Since every indecomposable $Q$-system is Morita equivalent to a connected one, and every $Q$-system is isomorphic to a direct sum of indecomposable $Q$-systems, to describe all of $[\QSys(\Hilb(G,\omega))]$ is suffices to describe bimodules between connected Q-systems. In the case when $\omega$ is trivial, there is a very nice description of these in terms of equivariantization which we shall now describe. This is the case that is described in Section \ref{subsec:groupactionex} in more elementary terms.

Let $\cC=\Hilb(G)$, and $(H,\mu)$ and $(K,\nu)$ two connected Q-systems. Then from the definitions, a $(H,\mu)$-$(K,\nu)$ bimodule is an object $X\in \Hilb(G)$ together with a family of isomorphisms
$$\alpha_{h,k}: h\otimes (X\otimes k)\cong X$$ such that the following diagram commutes: 
$$
\begin{tikzcd}
h^{\prime} \otimes ((h\otimes (X\otimes k))\otimes k^{\prime}) \arrow[swap]{d}{\cong}\arrow{rr}{1_{h^{\prime}}\otimes \alpha_{h,k}\otimes 1_{k^{\prime}}} & & h^{\prime}\otimes (X\otimes k^{\prime}) \arrow{d}{\alpha_{h^{\prime}, k^{\prime}}} \\
h^{\prime}h\otimes (X\otimes kk^{\prime}) \arrow{rr}{\mu(h^{\prime},h)\nu(k,k^{\prime})\alpha_{h^{\prime}h,kk^{\prime}}} & & X
\end{tikzcd}.
$$

We can rephrase this structure as follows: We have a natural $H\times K^{\op}$ action on $\Hilb(G)$ (just as a linear category) $G$ by $(h,k)\cdot X:=g\otimes (X\otimes k)$. We define the natural isomorphisms $$m_{(h^{\prime},k^{\prime}),(h,k)}:(h^{\prime},k^{\prime})\cdot((h,k)\cdot X)=h^{\prime} \otimes ((h\otimes (X\otimes k))\otimes k^{\prime})\rightarrow (h^{\prime}h,kk^{\prime})\cdot X. $$

The function $\mu\times \nu^{\op}$ defines a 2-cocycle on $H\times K^{\op}$, and this assembles into a categorical action of the group $H\times K^{\op}$ on $\Hilb(G)$. Then from the above characterization, we see that $(H,\mu)$-$(K,\nu)$ bimodules are precisely objects in the equivariantization of this $H\times K^{\op}$ action. This allows us to classify isomorphism classes of irreducible $(H,\mu)$-$(K,\nu)$ bimodules, by applying \cite[Corollary 2.13]{MR3059899}:

\begin{prop}\label{prop:idcpHKbimod}
Let $\Delta$ be the set of $H$-$K$ double cosets, and for each $\Gamma_{i}\in \Delta$, pick a representative $g_{i}\in \Gamma_{i}.$ 
Then isomorphism classes of irreducible $(H,\mu)$-$(K,\nu)$ bimodules are parameterized by pairs $(\Gamma_{i},\ \pi),$ where 
\begin{enumerate}
    \item    $\Gamma_{i}\in \Delta,$
    \item    $\pi$ is an irreducible $\delta$-projective representation of $\Stab_{H\times K^\op}(g_{i})$, where $\delta=\mu\times \nu^{\op}|_{\Stab(g_{i})}$.
\end{enumerate}
\end{prop}

It is easy to see from our description here that the data of a bimodule between $(H, \mu)$ and $(K,\nu)$ is precisely a $G$-graded $\bbC_{\mu}[H]-\bbC_{\nu}[K]$ bimodule as in Definition \ref{def:twisted bimodule}.

\end{ex}

\section{Actions of fusion categories on $\rm C^*$-algebras}\label{sec:Action}

In this section, we define the notion of an action of a fusion category on a $\rm C^*$-algebra, and various notions of equivalence. These ideas arise from subfactor theory, and have only begun to crystallize in the $\rm C^*$-context fairly recently. Thus we include full definitions and basic results which may be well known to experts, but to our knowledge have not appeared precisely in this form in the literature.

\begin{defn}
An action of a rigid $\rm C^*$-tensor category $\cC$ on a unital $\rm C^*$-algebra $A$ is a unitary tensor functor $F: \cC\rightarrow \Bim(A)\subseteq \rCorr(A\to A)$. 
\end{defn}

We also use the notation $\cC\overset{F}{\curvearrowright} A$ for an action. We can assemble all the actions of $\cC$ on $\rm C^*$-algebra in a $2$-category. Observe that since our tensor category $\cC$ above consists of dualizable objects, then the range of any such action $F$ will necessarily lie inside $\Bim(A),$ the rigid C*-tensor subcategory of dualizable right $A$-$A$ correspondences. 

\begin{defn}\label{defn:rCorrcC}
If $\cC$ is a unitary tensor category, $\rCorr_{\cC}$ is the $\rm C^*$ 2-category whose objects are $\rm C^*$-algebras equipped with a $\cC$ action, with $1$- and $2$-morphisms described below:

\begin{itemize}
\item
If $\cC\overset{F}{\curvearrowright} A$ and $\cC\overset{G}{\curvearrowright} B$ 
are actions of $\cC$ on $\rm C^*$-algebras $A$ and $B$ respectively, 
then a 1-morphism $F\rightarrow G$ 
consists of an element $X\in \rCorr(A\rightarrow B)$ 
and a family of unitary natural $A$-$B$ bimodular isomorphisms $\{u^{c}: F(c)\boxtimes_{A} X\rightarrow X\boxtimes_{B} G(c)\}_{c\in \cC}$ 
such that the following diagram (suppressing associators) commutes:
$$
\begin{tikzcd}
& F(a)\boxtimes_{A} F(b) \boxtimes_{A} X  
\arrow{dr}{F^2_{a,b}\boxtimes 1_X}
\arrow[swap]{dl}{1_{F(a)}\boxtimes u^{b}} &    
\\
F(a)\boxtimes_{A} X\boxtimes_{B} G(b) \arrow[swap]{dd}{u^{a}\boxtimes 1_{G(b)}} 
& &F(a\otimes b)\boxtimes_{A} X \arrow{dd}{u^{a\otimes b}}
\\
& &
\\
X\boxtimes_{B} G(a)\boxtimes_{B}G(b)
\arrow{rr}{1_X\boxtimes G^2_{a,b}} 
&    
& X\boxtimes_{B}G(a\otimes b)
\end{tikzcd}
$$
We call such a family a \textit{$\cC$-equivariant structure} on $X$.
When referring to a 1-morphism $(X,u)$ in $\rCorr_\cC$ we sometimes omit the $\cC$-equivariant structure $u$ if it raises no ambiguity.

\medskip

\item
If $(X,u), (Y,v)\in \rCorr_\cC(F\to G)$, 
and if $f:X\Rightarrow Y$ is an $A$-$B$ correspondence morphism, 
we say that $f\in\rCorr_\cC((X,u)\Rightarrow (Y,v))$ is a 2-morphism 
if moreover the following diagram commutes for all $a\in \cC:$

$$
\begin{tikzcd}
F(a)\boxtimes_{A} X \arrow[swap]{d}{1_{F(a)}\boxtimes f}\arrow{r}{u^{a}}& X\boxtimes_{B} G(a)\arrow{d}{f\boxtimes 1_{G(a)}} \\
F(a)\boxtimes_{A} Y \arrow{r}{v^{a}}& Y\boxtimes_{B} G(a)
\end{tikzcd}.
$$
\medskip

\item
Let $\cC\overset{F}{\curvearrowright} A$, $\cC\overset{G}{\curvearrowright} B$, and $\cC\overset{H}{\curvearrowright} C$ be objects in $\rCorr_\cC$. Let $(X,u)\in \rCorr_\cC(F\to G) $ and $(Y,v)\in \rCorr_\cC(G\to H)$ be $1$-morphisms. Then the composite $1$-morphism is given by 
$$(X\boxtimes_{B} Y, (u\boxtimes 1_{Y})\circ (1_{X}\boxtimes v) ) $$

\medskip

\item
Given correspondences $X, X'\in\rCorr(A\to B),\  Y, Y'\in\rCorr(B\to C)$ and intertwiners $f\in\rCorr_\cC(X\Rightarrow X'),$ and $g\in\rCorr_\cC(Y\Rightarrow Y')$ the horizontal composition of the $2$-morphisms $f$ and $g$ is defined as the ordinary horizontal product $f\boxtimes_B g$ in $\rCorr$.
\end{itemize}
\end{defn}

\begin{rem}There is another way to thing about this 2-category. Given any two $\rm C^*$ 2-categories $\mathcal{E}$ and $\mathcal{F}$, the collection of unitary $2$-functors, $2-\Fun(\mathcal{E}\to\mathcal{F}),$ itself forms a 2-category. If we use $\text{B}\cC$ to denote the 2-category with one object whose endomorphisms are the $\rm C^*$-tensor category $\cC$ then we see
$$2-\Fun(\text{B}\cC\to\rCorr)=\rCorr_{\cC}.$$

Indeed, the definitions presented above for $\rCorr$ are simply an unpacking of the the 2-category structure on $2-\Fun(\text{B}\cC\to\rCorr)$ (e.g. \cite[Section 4]{MR4261588}).

\end{rem}

There is another way to describe $\rCorr_{\cC}$ in terms of $\cC$-module categories.

\begin{prop}\label{prop:rCorr_CEquivCmod}
The 2-category $\rCorr_{\cC}$ is equivalent to the 2-category of $\cC$-module categories which are singly generated as unitary Cauchy complete categories.
\end{prop}

\begin{proof}

First note that unitary Cauchy complete categories $\cM$ which are singly generated are are equivalent to $\Mod(A)$. 
Indeed, let $x\in\cM$ be a generating object and $A:=\End(x)$ be the associated unital $\rm C^*$-algebra of endomorphisms. 
Then any object $c\in\cM$ is isomorphic to a summand of $x^{\oplus n}$, and hence corresponds to a projection in $\End(x^{\oplus n})\cong M_{n}(\End(x)) = M_n(A)$.
Note that every finitely generated projective $A$-module is isomorphic to $pA^n$ for some projection $p\in M_n(A)$.
Therefore, $\cM\cong \Mod(A)$.
Note that $\Fun(\Mod(A)\to\Mod(B))\cong \rCorr(A\to B)$ (see Remark \ref{rem:functor-bimod}), 
so from a tensor functor $F:\cC\to \End(\cM)$, 
we can construct a tensor functor $\cC\to \Bim(A)\subseteq \rCorr(A\rightarrow A)$ by conjugating the equivalence, which is a $\cC$-action on a $\rm C^*$-algebra as desired. Furthermore, composition of functors is isomorphic to the tensor product of bimodules.

Let $I:\cM\to \cN$ be a $\cC$-module functor. Note that $\cM\cong \Mod(A)$ and $\cN\cong \Mod(B)$, so the functor $I:\cM\to \cN$ corresponds to a functor $\Mod(A)\to \Mod(B)$, which gives a $A-B$ correspondence ${}_AX_B$. 
The $\cC$-module structure gives a unitary isomorphism $u^c: X\boxtimes_A F(c)\cong I(- )\lhd c \xrightarrow{\sim} I(- \lhd c)\cong G(c)\boxtimes_B X$.

Conversely, we can realize $\rCorr_\cC$ (contravariantly at the level of 1-morphisms) as the 2-category of right $\rm C^*$ $\cC$-module categories that are equivalent to finitely generated projective modules of a unital $\rm C^*$-algebra $A$.
The later are precisely the unitary Cauchy complete $\rm C^*$-categories with an additive generator.
An action of a fusion category $\cC$ on a $\rm C^*$-algebra $A$ is the same data as a $\cC$-module category structure on a unitarily Cauchy complete category together with a choice of generator,
which, indeed, is $\rCorr(\bbC\to A)$ with right $\cC$ action given by
$H\lhd c:=H\boxtimes_A F(c)$.
We call this $\cC$-module category as $\cM_{A,F}$.
\end{proof}

\begin{rem}\label{rem:functor-bimod}
Let $\cM$ and $\cN$ be $\rm C^*$-categories which are equivalent to $\Mod(A)$ and $\Mod(B)$ respectively. Then we used in the above proof the well-known fact that if $I:\cM\rightarrow \cN$ is a $*$-functor, then there is an $A$-$B$ correspondence $X$ with $I(H)\cong H\boxtimes_{A} X$. To construct $X$, simply take the right $B$-module $I(A)$, and equip this with a left $A$-action by sending $a\mapsto I(a)$ (where the morphism $a\in \End_{\Mod(A)}(A)\cong A$). Now notice that if $I(A)\cong B$, then the assignment $a\mapsto I(a)$ gives a unital $*$-homomorphism $\phi: A\rightarrow \End_{\Mod(B)}(B)\cong B$. Thus $I$ is represented by the correspondence ${}_{\phi}B\in \rCorr(A\rightarrow B)$. (See Definition \ref{defn:hombimod}.)
\end{rem}

Recall that $\fssMod(\cC)$ denotes the $\rm{C}^*$ 2-category of finitely semisimple right $\cC$-module categories. 
Let $\rCorr^{\fd}_{\cC}$ denote the 2-subcategory of $\rCorr_\cC$ with \textit{finite dimensional} $\rm C^*$-algebras.

\begin{cor}\label{cor:fdrCorr_CEquivfssCmod}
$\rCorr^{\fd}_{\cC}\cong \fssMod(\cC)$.
\end{cor}
\begin{proof}
Given $\cC\overset{F}{\curvearrowright} A$, since $A$ is finite dimensional,
$\Mod(A)$ is finitely semisimple and so is $\cM_{A,F}$.
Conversely, let $\{m_1,\cdots, m_N\}$ be a list of all non-isomorphic simple objects in $\cM$ 
and $A:= \End_\cM\left(\bigoplus_{i=1}^N m_i\right),$ which is a finite-dimensional $\rm C^*$-algebra.
\end{proof}

\bigskip

There are a distinguished collection of correspondences which correspond to homomorphisms of $\rm C^*$-algebras. 
\begin{defn}\label{defn:hombimod}
Given a unital $*$-homomorphism $\phi:A\rightarrow B,$ we define $${}_{\phi} B\in \rCorr(A\to B)$$ as the right Hilbert $B$-module $B$, with left action of $A$ given by $a\rhd b:=\phi(a)b$. 
\end{defn}

Up to isomorphism, the correspondences in the above definition are characterized as those which are isomorphic to $B$ as a right $B$-module. We now discuss how to equip such correspondences with the structure of a $1$-morphism in $\rCorr_{\cC}$.

\begin{lem}\label{lem:StrucEqStrucs}
Let $\cC\overset{F}{\curvearrowright} A$ and $\cC\overset{G}{\curvearrowright} B$ be actions on $\rm C^*$-algebras,
and let $\phi: A\rightarrow B$ be an injective, unital $\rm C^*$-algebra homomorphism.
Then there is a bijection between $\cC$-equivariant structures on ${}_{\phi} B$ and a family of linear maps: 
$$\{h^{c}: F(c)\rightarrow G(c)\}_{c\in\cC}$$ such that 
\begin{enumerate}[label=(\arabic*)]
\item $h^{c}(a\rhd x \lhd a')=\phi(a)\rhd h^{c}(x)\lhd \phi(a')$ for $a,a'\in A$ and $x\in F(c)$;
\item
For any morphism $f\in \cC(c\to d)$, $G(f)\circ h^{c}=h^{d}\circ F(f)$;
\item
$\phi(\langle x | y\rangle_A)=\langle h^{c}(x) | h^{c}(y)\rangle_B$ for $x,y\in F(c)$;
\item
$h^{c}(F(c))B=G(c)$;
\item The following diagram commutes:
\[
\begin{tikzcd}
F(c)\otimes F(d)
\arrow[swap]{d}{h^{c}\otimes h^{d}}
\arrow{r}{}
& F(c)\boxtimes_{A} F(d)
\arrow{r}{F^{2}_{c,d}}
& F(c\otimes d)\arrow{d}{h^{c\otimes d}}
 \\
G(c)\otimes G(d)
\arrow{r}{}
& G(c)\boxtimes_{B} G(d)
\arrow{r}{G^{2}_{c,d}}
& G(c\otimes d)
\end{tikzcd}
\]
\end{enumerate}

Furthermore, given a $\rm C^*$-algebra homomorphism $\omega:A\to B$, then two equivariant correspondences of the form $({}_{\phi} B, \{h^{c}\})$ and $({}_{\omega} B, \{k^{c}\})$ are isomorphic if and only if there exists a unitary $v\in B$ such that such that $\Ad(v)\circ \phi=\omega$ and $\Ad(v)\circ h^{c}=k^{c}$ for all $c\in \cC$.
\end{lem}

\begin{proof}
Given the above $\{h^c\}$ data, we define a $\cC$-module structure on $- \boxtimes_{A} {}_{\phi}B$ by 
$u^{c}:F(c)\boxtimes_{A} {}_{\phi} B\rightarrow {}_{\phi} B\boxtimes_{B} G(c)$ by $u^{c}(x\boxtimes b):=1_{B}\boxtimes_{B} h^{c}(x)b$.
Property (1) above ensures $u^{c}$ is $A$-balanced and $A-B$ bimodular.
Naturality follows from (2), 
unitarity from (3) and (4), 
and (5) gives the extra compatibility condition for the $\cC$-module structure.

In the other direction, given a family of unitaries $\{u^{c}: F(c)\boxtimes_{A} {}_{\phi} B\rightarrow {}_{\phi} B\boxtimes_{B} G(c)\}$, we note that we have a canonical $A-B$ bimodule isomorphism $f: {}_{\phi} B\boxtimes_{B} Y \cong {}_{\phi} Y$, which sends $b\boxtimes y\mapsto b\rhd y$. Define 
$$h^{c}(x):=f(u^{c}(x\boxtimes 1_{B})).$$
Clearly these two assignments are inverse to each other by construction. It is straightforward to see that $h^{c}$ satisfies the required criteria.

Note that by directly applying the definition of 2-morphism in $\rCorr_{\cC}$, then $({}_{\omega} B, u_{k})$ is unitarily isomorphic to $({}_{\phi} B, u_{h})$
if and only if there is an $A$-$B$ bimodular unitary $s: {}_{\phi} B\to {}_{\omega} B$ satisfying the following equations for all $c\in \cC$:
$$(s\boxtimes_B \id_{G(c)})\circ u_{h}^{c} = u^{c}_{k}\circ (\id_{F(c)}\boxtimes_A\ s)$$
or equivalently 
$$u^{c}_{k} = (s\boxtimes_B \id_{G(c)})\circ u_{h}^{c}\circ (\id_{F(c)}\boxtimes_A\ s^{*}).$$

Note that since $s: B_{B}\rightarrow B_{B}$, then there is a unique unitary $v\in B$ such that $s(b):=vb$. Furthermore, since $s$ is also left $A$-modular we obtain 
$$v\phi(a)b=s(a\rhd_{\phi} b) = a\rhd_{\omega} s(b) =\omega(a)vb \qquad\Rightarrow \qquad \Ad(v)\circ \phi(-) = \omega(-).$$ 

Now applying both sides of the above equation to the vector $x\boxtimes 1_{B}\in F(c)\boxtimes_{B} {}_{\phi} B$ and then composing with $f_{c}$ yields
\begin{align*}
k^c(x) & = f(u^{c}_{k}(x\boxtimes 1_B)) = f(  (s\boxtimes_B \id_{G(c)})\circ u_{h}^{c}\circ (\id_{F(c)}\boxtimes_A\  s^{*})(x\boxtimes 1_B) ) \\
& = f((s\boxtimes_B \id_{G(c)})\circ u_{h}^{c})(x\boxtimes v^*) ) = f((s\boxtimes_B \id_{G(c)})(1_B\boxtimes h^c(x)v^*) ) \\
 & = f(v\boxtimes h^c(x)v^*) = vh^c(x)v^* = \Ad(v)\circ h^c(x). 
\end{align*}
\end{proof}

We now introduce two notions of equivalence of $\cC$-actions on $\rm C^*$-algebras. 
Let $\cC\overset{F}{\curvearrowright}A$ and $\cC\overset{G}{\curvearrowright}B$ be objects in $\rCorr_{\cC}$.
An isomorphism from $F$ to $G$ is a $1$-morphism $(X, u) :F\rightarrow G$ such that there exists a 1-morphism $ (Y,v):G\rightarrow F$ with $X\boxtimes_{B} Y\cong {}_{A}A_{A}$ and $Y\boxtimes_{A} X\cong {}_{B} B_{B}$. 
An isomorphism $(X, u)$ is called a \textit{strong isomorphism} if there exists an $\rm C^*$-algebra isomorphism $\phi: A \rightarrow B$ such that $X={}_{\phi} B$. 

\begin{defn}\label{defn:EqActions}
Let $\cC$ be a $\rm C^*$-tensor category, and $\cC\overset{F}{\curvearrowright}A$ and $\cC\overset{G}{\curvearrowright}B$ be actions on unital $\rm C^*$-algebras $A$ and $B,$ respectively. 
Then we say $F$ and $G$ are
\begin{itemize}
\item
    \textit{Morita equivalent} if there exists an isomorphism in $\rCorr_\cC$ from $F$ to $G$.
\item
    \textit{Equivalent} if there exists a strong isomorphism in $\rCorr_\cC$ from $F$ to $G$.
\end{itemize}
\end{defn}

\noindent Morita equivalence of two actions implies that $A$ and $B$ are Morita equivalent as algebras, while equivalence implies $A$ and $B$ are isomorphic as algebras. We note that viewing an action of a finite group on a $\rm C^*$-algebra as an action of the fusion category $\Hilb(G)$ as in Section \ref{subsec:realization}, then equivalence of two such actions reduces to cocycle conjugacy. 

\begin{rem}
${}_{\phi} B$ is invertible as a $1$-morphism as in Definition \ref{defn:hombimod} if and only if $\phi:A\rightarrow B$ is an isomorphism of $\rm C^*$-algebras. In light of Lemma \ref{lem:StrucEqStrucs}, to describe an equivalence between actions, it is necessary and sufficient to find an isomorphism between the $\rm C^*$-algebras and a $\cC$-equivariant structure. 
\end{rem}

\begin{rem}
If $\cC\overset{F}{\curvearrowright}A$ and $\cC\overset{G}{\curvearrowright}A$ are two actions on the same algebra, there is an even stronger notion of equivalence, namely we could ask for the functors $F, G:\cC\rightarrow \Bim(A)$ to be unitarily naturally isomorphic. This is the categorical analogue of cocycle equivalence for group actions. This will appear in passing later in the paper, but in general natural isomorphism is too fine of an equivalence relation to have a tractable classification scheme. 
\end{rem}


\subsection{Realization}\label{subsec:realization}

\ \\
Given an action $\cC\overset{F}{\curvearrowright}A$, 
there is a canonical \textit{realization 2-functor} 
$$| - |_{F}: \QSys(\cC)\rightarrow \rCorr.$$
For any Q-system $P$, the realization $|P|_{F}$ should be thought of as a ``generalized crossed-product" \cite{MR4419534}. 
It can be defined (as a vector space) as
$$|P|_{F}:=F(P)\cong \bigoplus_{c\in \Irr(\cC)} \cC(c\to P)\otimes_{\bbC} F(c).$$
The product $F(P)\otimes_{\bbC} F(P)\rightarrow F(P)$ is given by the composition
$$F(P)\otimes_{\bbC} F(P)\rightarrow F(P)\boxtimes_{A} F(P)\xrightarrow{F^2_{P,P}} F(P\otimes P)\xrightarrow{F(m_P)} F(P).$$

In terms of our direct sum decomposition, for homogeneous tensors $x=x_{1}\otimes x_{2}\in \cC(c\to P)\otimes F(c)$ and $y=y_{1}\otimes y_{2}\in \cC(d\to P)\otimes F(d)$ the product becomes
$$x\cdot y:=\sum_{\substack{e\in\Irr(\cC) \\ \alpha\in B(e,cd)}} \left(m_P\circ (x_{1}\otimes y_{1})\circ \alpha \right)\otimes_{\bbC} (F(\alpha^*)\circ  F^{2}_{c,d}(x_{2}\boxtimes y_{2})),$$
where $B(e,cd)$ is an orthonormal basis for the Hilbert space $\cC(e\to c\otimes d)$ with composition inner product. For the norm and $*$-structure, see \cite[Section 4.1]{MR4419534}.

Now, if $P,Q$ are two $Q$-systems and $X\in \QSys(\cC)(P\to Q)$, then we define the realization as a vector space
$$|X|_{F}:= F(X)\cong \bigoplus_{c\in \Irr(\cC)} \cC(c\to X)\otimes_{\bbC} F(c)$$
with left $|P|_{F}$ and right $|R|_{F}$ actions defined similarly as above (for a detailed description together with Hilbert module structure see \cite[Section 4.2]{MR4419534}).


We claim that Q-system realization generalizes the (twisted) crossed-product by a group action. Indeed, let $G$ be a finite group, and suppose we have an action $ G\rightarrow \Aut(A)$.
This defines a unitary tensor functor
\begin{align*}
F:\Hilb(G)&\to \Bim(A) \qquad \text{ with tensorator } &&F^{2}_{g,h}:F(g)\boxtimes_{A} F(h)\cong F(gh)\\
g&\mapsto{}_{g} A  &&\hspace{2.5cm} a\boxtimes b\mapsto h^{-1}(a) b.
\end{align*}
Here ${}_{g} A$ is $A$ as a right Hilbert $A$-module and $a\rhd b=g^{-1}(a)b$.

Now, consider a connected Q-system $P=(H,\mu)\in \Hilb(G)$, where $\mu\in Z^{2}(H, \text{U}(1))$ is a normalized 2-cocycle. The multiplication $m:P\otimes P\rightarrow P$ can be described as
$$m=\sum_{g,h\in H} \mu(g,h) v_{gh}\circ \alpha_{g,h}\circ (v^{*}_{g}\otimes v^{*}_{h}),$$
\noindent where $\{v_{g}\in \Hilb(G)(\bbC_g\to P)\}_{g\in H}$ are isometries implementing the direct sum decomposition $P=\bigoplus_{g\in H} \bbC_{g}$, and $\alpha_{g,h}\in \Hilb(G)(\bbC_{g}\otimes \bbC_{h}\to \bbC_{gh})$ are the canonical unitary isomorphisms.

Now we see that 
$$|P|_{F}\cong \bigoplus_{g\in H} \Hilb(G)(\bbC_{g}\to P)\otimes_{\bbC} {}_{g} A.$$
The right Hilbert $A$ module $\Hilb(G)(\bbC_{g}\to P)\otimes_{\bbC} {}_{g} A$ has a standard singleton Pimsner-Popa basis, namely $v_{g}\otimes 1_{A}$. Thus since $v_{g}$ is an isometry, every element in $|P|_{F}$ can be written uniquely as 
$$\sum_{g\in H} v_{g}\otimes a_{g},$$
where $a_{g}\in A$.
Then from the definition of the product of the realization, we obtain an isomorphism (ultimately of $\rm C^*$-algebras)
$$|P|_{F}\cong A\rtimes_{\mu} H,$$
where the twisted crossed product $A\rtimes_{\mu}H$ is as defined in Section \ref{subsec:groupactionex}. Furthermore, in that section the realization of bimodules between Q-systems is described in elementary terms.

\subsection{The invariant}\label{sec:invariant}

\begin{defn}
$\PAbG$ is the category of pre-ordered abelian groups:
\begin{itemize}
\item 
Objects are pairs $(G,G^{+})$ where $G$ is an abelian group and $G^{+}\subseteq G$ contains $0$, and satisfies $G^{+}-G^{+}=G$ and $G^{+}+G^{+}=G^{+}$. 
\item
Morphisms $(G,G^{+})$ to $(H,H^{+})$ consist of group homomorphisms $f:G\rightarrow H$ such that $f(G^{+})\subseteq H^{+}$.
\end{itemize}
\end{defn}

We call $(G,G^{+})$ a pre-ordered abelian group, $G^{+}$ the positive cone, and the morphism in our category positive homomorphisms. The reason for our terminology is that the relation defined by $a\le b $ if $b-a\in G^{+}$ gives a translation invariant pre-order on $G$ (and every translation invariant pre-order is precisely determined by its positive cone). Furthermore positive homomorphisms are precisely the homomorphisms which respect $\le$.

$\PAbG$ is naturally enriched over the symmetric monoidal category of abelian monoids. Indeed, $\PAbG( (G, G^{+})\to (H, H^{+}) )$ naturally has the structure of an abelian monoid: the pointwise sum of two positive homomorphisms is again a positive homomorphism, and the $0$ homomorphism is the identity. Furthermore composition is clearly bilinear.

The following definition will be useful often in the remainder of this manuscript: 
\begin{defn}
For any two categories $\cC$ and $\cD$ enriched over abelian monoids, we denote the category of additive functors by $\Fun_{+}(\cC\to \cD)$.
\end{defn}

Let $[\rCorr]$ denote the 1-category obtained from the 2-category $\rCorr$ by decategorification (i.e. taking morphisms in the new 1-category to be \textit{isomorphism classes} of 1-morphisms in the original 2-category). The category $[\rCorr]$ is also enriched in abelian monoids, with $[X]+[Y]:=[X\oplus Y]$ and unit $[0]$.

We have the following well-known proposition concerning the familiar $K_{0}$ group from $\rm C^*$-algebra theory:

\begin{prop}\label{prop:Grothendieck}
$K_{0}$ defines an additive functor from $[\rCorr]\rightarrow \PAbG$.
\end{prop}
\begin{proof}

Note that the Grothendieck construction $\cG$ (for example, see \cite{MR1783408}) is in fact a functor from abelian monoids to pre-ordered abelian groups, where the positive cone is given as the image of the monoid in the group completion.
For a given $\rm C^*$-algebra $A$, the pre-ordered abelian group $(K_{0}(A), K^{+}_{0}(A))$ is defined as the Grothendieck completion of the abelian monoid $[\rCorr](\bbC\to A)=[\Mod(A)]$, which by definition is the abelian monoid of isomorphism classes of finitely generated projective $A$-modules. Then for a morphism ${}_A[X]_B\in [\rCorr](A\to B)$, $-\boxtimes_A [X]_B$ induces an abelian monoid homomorphism $[\rCorr](\bbC\to A)\rightarrow [\rCorr](\bbC\to B)$. By functoriality of the Grothendieck construction, this yields a positive homomorphism on the $K_{0}$ groups. 
\end{proof}

Recall that $[\QSys(\cC)]$ is the once decategorified 1-category from Definition \ref{defn:[QSys]}, whose objects are Q-systems in $\cC$ and whose morphisms are isomorphism classes of bimodules in $\cC$. We then have the following definition: 

\begin{defn}
\label{Def:InvariantFhat}
Let $\cC\overset{F}{\curvearrowright}A$ be an action of a fusion category on the $\rm C^*$-algebra $A$. We define the functor 

\begin{align*}
    \widehat{F}:=K_{0}\left(\ \Big[\ |\ -\ |_{F}\ \Big]\ \right):[\QSys(\cC)] & \rightarrow \PAbG\\
    Q &\mapsto \left(K_{0}(|Q|_{F}),K^{+}_{0}(|Q|_{F})\right)\\
    _{P}[X]_{Q}& \mapsto -\boxtimes_{|P|_{F}}\Big[ |X|_{F} \Big]_{|Q|_F}.
\end{align*}
\end{defn}

The functor $\widehat{F}$ will be the main ingredient of our invariant.

\begin{defn}[Pointed Invariant]\label{pointed invariant}
Let $\cC\overset{F}{\curvearrowright}A$ be an action. Then we define the \textit{pointed invariant} to be the pair $(\widehat{F}, [A])$, where $[A]\in \widehat{F}(1_{\cC})\cong K_{0}(A)$. An equivalence of pointed invariants $(\widehat{F}, [A])$ and $(\widehat{G}, [B])$ is a natural isomorphism $\alpha: \widehat{F}\cong \widehat{G}$ such that $\alpha_{1_{\cC}}([A])=[B]$.
\end{defn}

To understand $\widehat{F}$ better we make the following remarks, leading to Theorem \ref{thm:EquivariantCrossedprod}: 

\bigskip

Consider the functor
\begin{align*}\begin{split}
    [\Mod(\cC)(-\to \cM_{A,F})]:[\Mod(\cC)]^\op &\to \AbM\\
    \cN&\mapsto \left(\left[\Mod(\cC)(\cN \to \cM_{A,F}) \right], \oplus\right),
\end{split}\end{align*}
which is a contravariant functor from decategorified $\cC$-module categories $[\Mod(\cC)]$ to the category of abelian monoids $\AbM.$ 
Here, $\cM_{A,F}$ is a right $\cC$-module category whose objects are right fgp $A$-modules with $\cC$-action given by
$X\lhd c:=X\boxtimes_A F(c)$.
Moreover, since $\cM_{A,F}\subset \rCorr(\bbC\to A)$ admits direct sums and has a zero object,
then the $\left[\Mod(\cC)(\cN\to \cM_{A,F})\right]$ becomes an abelian monoid, 
by defining the sum $\oplus$ of $\cC$-module functors object-wisely on $\cN$.

\bigskip

\begin{thm}\label{thm:EquivariantCrossedprod}
Given an action $\cC\overset{F}{\curvearrowright}A$, we have a natural isomorphism $$\eta:\widehat{F}\cong \cG\circ [\Mod(\cC)(\Psi(-)\to \cM_{A,F})],$$
where $\Psi:\QSys(\cC)\to \fssMod(\cC)\hookrightarrow\Mod(\cC)$ is the 1-contravariant 2-functor in Remark \ref{rem:unitaryostrik}. Here, $\cG$ denotes the Grothendieck completion described in Proposition  \ref{prop:Grothendieck}.

\end{thm}
\begin{proof}
It suffices to build a natural isomorphism of functors
$$\eta:D(|\ \cdot \ |_F)\Rightarrow [\Mod(\cC)(\Psi(-)\to \cM_{A,F})],$$
where $D(\cdot)$ denotes the functor from $[\rCorr]\rightarrow \AbM$ that sends a $\rm C^*$-algebra $A$ to its monoid of equivalence classes of finitely generated projective modules $[\rCorr](\bbC\to A)$, and a $\rm C^*$-correspondence ${}_{A}[X]_{B}\mapsto - \boxtimes_{A} [X]_{B}$ as in the proof of Proposition \ref{prop:Grothendieck}. 
Then $\widehat{F}=\cG( D(\ |\ \cdot\ |_{F}\ ))$, so that the desired natural isomorphism is obtained as $\cG(\eta)$.

Let $Q\in \QSys(\cC)$. We need to define $\eta_{Q}:D(|Q|_{F})\rightarrow [\Mod(\cC)(\Psi(Q)\to \cM_{A,F})]$. Let $[Y]\in D(|Q|_{F})=[\rCorr(\bbC\to |Q|_{F})]$ be an equivalence class of finitely generated projective $|Q|_F$-modules. 
Recall that $\Psi(Q) = {}_Q \cC$ is a right $\cC$-module category consisting of left $Q$-module objects $(x,\lambda)$, where $\lambda\in\cC(Q\otimes x\to x)$ satisfies the left $Q$-module compatibility condition (B1) from Section \ref{sec: Q-systems}.
The natural (right) $\cC$-module category structure is given by $(x,\lambda)\lhd c: = (x\otimes c,\lambda\otimes 1_c)$. Given an $(x,\lambda)\in {}_{Q} \cC$ define the functor $L_{Y}: \Psi(Q)\rightarrow \Mod(A)$ by
$$L_{Y}((x,\lambda)):=Y\boxtimes_{|Q|_{F}} |x|_{F}$$
\noindent Here realization $|x|_F$ is a $|Q|_F-A$ bimodule, so $Y\boxtimes_{|Q|_{F}} |x|_{F}\in \Mod(A)$ which is the category underlying $\cM_{A,F}$. To equip $L_{Y}$ with the structure of a $\cC$-module functor, we use the composite natural isomorphism 
\begin{align*}
L_{Y} ((x,\lambda)\lhd c) & = Y\boxtimes_{|Q|_F} (|x\otimes c|_F) \xrightarrow{1_Y\boxtimes |F^{2}_{x,c}|^{-1}} Y\boxtimes_{|Q|_F} (|x|_F\boxtimes_A F(c)) 
\\ & \cong_{assoc} (Y\boxtimes_{|Q|_F} |x|_F)\boxtimes_{A} F(c)= L_{Y}((x,\lambda))\lhd c.     
\end{align*}

It is easy to see these isomorphisms equip $L_{Y}$ with the structure of a $\cC$-module functor whose isomorphism class only depends on the isomorphism class of $Y$.

Now we define 
$$\eta_{Q}([Y]):=[L_{Y}]\in [\Mod(\cC)(\Psi(Q)\to \cM_{A,F})]$$
\noindent Clearly $\eta_{Q}$ is additive, since $\oplus$ commutes with $\boxtimes$ up to isomorphism.

To see naturality, for any $Q,R\in [\QSys(\cC)]$ and $[X]\in [\QSys(\cC)](Q\to R)$, then $|X|_{F}$ is a $|Q|_{F}-|R|_{F}$ correspondence, and we directly see
\begin{align*}
\eta_{R}\circ D(|X|_{F})([Y])&=\eta_{R}([Y\boxtimes_{|Q|_{F}} |X|_{F}]) 
= L_{Y\boxtimes_{|Q|_{F}} {|X|_{F}}}\\
&= (-\circ \Psi(X))\circ \eta_{Q}([Y]).    
\end{align*}
i.e. the following diagram commutes 
\[
\begin{tikzcd}[column sep=4em]
D(|Q|_F)
\arrow[swap]{d}{D(|X|_F)}
\arrow{r}{\eta_Q}
&
{[}\Mod(\cC)(\Psi(Q)\to \cM_{A,F}){]}
\arrow{d}{-\circ \Psi(X)}
\\
D(|R|_F)
\arrow{r}{\eta_R}
& {[}\Mod(\cC)(\Psi(R)\to \cM_{A,F}){]}
\end{tikzcd}
\]
for any $[Y]\in D(|Q|_F)$.
Thus $\eta$ is natural.

We shall now show $\eta$ is an isomorphism by building an inverse $\eta^{-1}_{Q}$ for each $Q\in [\QSys(\cC)]$. 
For a given $\cC$-module functor $I\in\Mod(\cC)({}_Q\cC\to \cM_{A,F})$, 
then for the canonical object $(Q,m)\in {}_Q\cC$ ($Q$ as a module over itself by multiplication)
$I((Q,m))\in \Mod(A)$, but has the additional right $|Q|_F$-module structure defined by 
$$I((Q,m))\boxtimes_A |Q|_F = I((Q,m))\lhd Q\cong I((Q,m)\lhd Q)=I((Q\otimes Q, m)) \to I((Q,m)),$$
where the last morphism is $I(m): I(Q\otimes Q)\to I(Q)$. 
Hence we can view $I((Q,m))\in \Mod(|Q|_{F})$ (cf. \cite[Lemma 5.11]{2021arXiv211106378C}).
Thus $[I((Q,m))]\in D(|Q|_{F})$.
Now we define 
$$\eta_Q^{-1}([I]) := [I((Q,m))].$$

Now we check these are inverses (as set functions):  
starting with $[Y]\in D(|Q|_F)$, $\eta_Q^{-1}\circ\eta_Q([Y]) = [Y\boxtimes_{|Q|_F}|Q|_F]\cong [Y]$.
In the other direction, let $I\in\Mod(\cC)({}_Q\cC\to \cM_{A,F})$, and $(x,\lambda)\in {}_Q\cC.$ Then 
$\eta_Q\circ\eta^{-1}_Q(I)((x,\lambda)) = I(Q,m)\boxtimes_{|Q|_F} |x|_{F} = I(Q\otimes_Q (x,\lambda)) \cong I(x,\lambda)$.
\end{proof}

We note that Theorem \ref{thm:EquivariantCrossedprod} can be viewed as a functorial generalization of Julg's Theorem 
\cite[Theorem 2.6.0]{MR911880}
relating equivariant K-theory and the ordinary K-theory of the crossed product.
If $\cC=\Hilb(G)$ and $\cC\overset{F}{\curvearrowright} A$ arises from a homomorphism $f:G\rightarrow \Aut(A)$, 
we can consider the $Q$-system $Q:=\bbC[G]\in \Hilb(G) $.
The corresponding module category is $\Hilb$ with trivial $G$-action. 
Then $|Q|_{F}=A\rtimes G$, and $\cG([\Mod(\cC)(\Hilb\to \cM_{A,F})])$ is isomorphic to the equivariant K-theory. 

\begin{cor}\label{cor:StrEqPtInv}
Let $\cC\overset{F}{\curvearrowright} A,$ $\cC\overset{G}{\curvearrowright} B$ and $\cC\overset{H}{\curvearrowright} C$ be actions. We then have that: 
\begin{enumerate}
\item 
For any $[X]\in [\rCorr_{\cC}(F\to G)]$, the positive homomorphism $$[- \boxtimes_{A} X]: (K_{0}(A), K^{+}_{0}(A))\rightarrow (K_{0}(B), K^{+}_{0}(B))$$ canonically extends to a natural transformation $\theta_{X}: \widehat{F}\Rightarrow \widehat{G}$,
\item
For the unit correspondence $[A]\in [\rCorr_{\cC}(F\to F)]$, $\theta_{A}=\id_{\widehat{F}}$,
\item
If $[Y]\in [\rCorr_{\cC}(G\to H)]$, then $\theta_{X\boxtimes_{B} Y}=\theta_{Y}\circ \theta_{X}$,
\item
If $\cC\overset{F}{\curvearrowright} A$ is Morita equivalent to $\cC\overset{G}{\curvearrowright} B$, then $\widehat{F}\cong \widehat{G}$,
\item
If $\cC\overset{F}{\curvearrowright} A$ is equivalent to $\cC\overset{G}{\curvearrowright} B$, then there is an isomorphism of pointed invariants $\left(\widehat{F}, [{A}]\right)\cong \left(\widehat{G}, [{B}]\right)$.
\end{enumerate}
\end{cor}

\begin{proof}
The first three statements say that the assignment $F\mapsto \widehat{F}$ and $[\rCorr_{\cC}(F\to G)]\ni [X]\mapsto \theta_{X}$ is a functor from $[\rCorr_{\cC}]$ to $\Fun_{+}([\rCorr_{\cC}]\to \PAbG)$. The fourth and fifth statements follow immediately from the first three.

Given $X\in \rCorr_{\cC}(F\to G)$, then by Proposition \ref{prop:rCorr_CEquivCmod}, $R_{X}:=[-\boxtimes_{A} X]\in [\Mod(\cC)(\cM_{A,F}\to \cM_{B,G})]$. 
Let $\eta_F:\widehat{F}\cong \cG\circ [\Mod(\cC)(\Psi(-)\to \cM_{A,F})]$ and $\eta_G:\widehat{G}\cong \cG\circ [\Mod(\cC)(\Psi(-)\to \cM_{B,G})]$ be the natural isomorphisms from Theorem \ref{thm:EquivariantCrossedprod} for $F$ and $G$, respectively. 
Define $\theta_{X}$ to be the composition of natural transformations
$$\widehat{F}\xrightarrow{\eta_F} \cG\circ[\Mod(\cC)(\Psi(-)\to\cM_{A,F})]\xrightarrow{\cG(R_{X}\circ -)}\cG\circ[\Mod(\cC)(\Psi(-)\to\cM_{B,F})]\xrightarrow{\eta_G^{-1}} \widehat{G},$$
where $\cG:\AbM\to \PAbG$ denotes the Grothendieck functor as in the proof of the previous theorem. Then the results follow from functoriality of $\cG$ and cancelling $\eta^{-1}\circ \eta$. 
\end{proof}

\section{Inductive limit actions on AF $\rm C^*$-algebras}\label{sec:ILActions}

AF-actions are a special class of actions of a fusion category on AF $\rm C^*$-algebras. These are actions which are in a particular sense ``inductive limits" of actions on finite dimensional $\rm C^*$-algebras. As the latter are well understood, it is reasonable to expect that we can leverage this knowledge to obtain classifications for inductive limits of finite dimensional actions, in a spirit similar to Elliott's K-theoretic classification of AF $\rm C^*$-algebras, with the role of $K_{0}$ played by our $\widehat{F}$.

Indeed, as discussed above, one can view the program of classifying AF-actions as a special case of studying AF $\rm C^*$-algebras \textit{internal to the fusion category $\cC$}, in the sense of \cite{MR3687214}. Thus our arguments will be essentially ``$\cC$-equivariant" versions of Elliott's original arguments, with some higher categorical modifications necessary at various points.

The outline of this section is as follows: First, we provide a general framework for inductive limit actions of fusion categories on inductive limit $\rm C^*$-algebras, and show $\widehat{F}$ is continuous with respect to inductive limits. Then we describe a framework for enriched Bratteli diagrams, which gives us a way to easily describe an AF-action $\cC\overset{F}{\curvearrowright} A$ and compute its invariant $\widehat{F}$. Finally, we give a proof that $\widehat{F}$ is a complete invariant amongst the class of AF-actions.

\subsection{Inductive limit $\rm C^*$-algebras}

\begin{defn}\label{defn:IL2Cat}
We define a 2-category $\IL$ of inductive limit $\rm C^*$-algebras as follows:
\begin{enumerate}[label=(IL\arabic*)]\addtocounter{enumi}{-1}
    \item\label{IL:Alg}
        Objects are sequences of unital $\rm C^*$-algebras $\{A_{n}\}_{n\ge 0}$, together with unital inclusions $i^{A}_{n}: A_{n}\rightarrow A_{n+1}$. 
    \item\label{IL:Corr}
        A 1-morphism $\{X_n\}: \{A_{n}\}\rightarrow \{B_{n}\}$ is a sequence $\{X_{n}\}_{n\ge 0}$ where each $X_{n}\in \rCorr(A_{n}\to B_{n})$, together with a sequence of $\bbC$-linear maps $j^{X}_{n}: X_{n}\rightarrow X_{n+1}$ such that $j^X_{n}(X_{n})\lhd B_{n+1}=X_{n+1},$ and for $a\in A_{n},\ b\in B_{n},\ x\in X_{n}$ we have
        $$j^{X}_{n}(a\rhd x\lhd b)=i^{A}_{n}(a)\rhd j^{X}_{n}(x) \lhd i^{B}_{n}(b)\quad \text{ and } \quad i^B_{n}\left(\langle x | y\rangle_{B_{n}}^{X_n}\right)=\left\langle j^X_{n}(x) \big| j^X_{n}(y) \right\rangle_{B_{n+1}}^{X_{n+1}}.$$
    \item\label{IL:Intertwiner}
        If $\{X_{n}\},\{Y_{n}\}$ are 1-morphisms from $\{A_{n}\}$ to $\{B_{n}\}$, a 2-morphism $f=\{f_n\}: \{X_n\}\Rightarrow \{Y_n\}$ consists of a norm-bounded sequence of 2-morphisms $\{f_{n}\in \rCorr(X_{n}\Rightarrow Y_{n})\}$  satisfying for all $x_n\in X_{n}$
        $$j^{Y}_{n}\left(f_{n}(x_n)\right)=f_{n+1}\left(j^{X}_{n}(x_n)\right).$$
\end{enumerate}
 
    \begin{itemize}
    \item
    Composition of $1$-morphisms $\{X_n\}\in\IL(\{A_n\}\to \{B_n\})$ and\newline  $\{Y_n\}\in\IL(\{B_n\}\to \{C_n\})$ is given by 
    $$\{X_{n}\}\boxtimes_{\{B_n\}} \{Y_{n}\}:=\{X_{n}\boxtimes_{B_{n}} Y_{n}\},\quad\text{ where }\quad j^{X\boxtimes_B Y}_n:=j^{X}_{n}\boxtimes_{B_n} j^{Y}_{n}.$$
    
    \item Horizontal product of $2$-morphisms $\{f_n\},\{g_n\}$ is given by $\{f_n\}\boxtimes \{g_n\} := \{f_{n}\boxtimes g_{n}\}$. Composition of 2-morphisms is as usual.
    \end{itemize}
\end{defn}

Note that in the above definition, $i_{n}$ above depends on the object $A$ and $j_{n}$ depends on the 1-morphism $X$, but we will drop this dependence from the notation as it is usually clear from context. Moreover, we write $j_{m,n}:=0$ if $m<n,$ $j_{n,n}:=\id_{X_n},$ and $j_{m,n}= j_{m-1}\circ\hdots\circ j_n$ for $m>n,$ and similarly for the $i_n$ maps.

\begin{remark}\label{rmk:InjConnMaps}

We have chosen our connecting maps $i_n,j_n$ to be injective at each level. Therefore, the inclusion maps $i_n$ from \ref{IL:Alg} in Definition \ref{defn:IL2Cat} are necessarily isometric. Moreover, this implies that the maps $j_n$ from \ref{IL:Corr} are isometric as well, since the $B_n$-valued inner products are assumed non-degenerate. 
\end{remark}

We have the following observations concerning Pimsner-Popa bases of inductive limit correspondences.

\begin{lem}\label{lem:SamePP}
Let $\{b_{i}\}^{m}_{i=1}$ be a Pimsner-Popa basis for $X_{n}$. Then $\{j^{X}_{n}(b_{i})\}$ is a Pimsner-Popa basis for $X_{n+1}$.
\end{lem}

\begin{proof}
We must show the equality $\sum^{m}_{i=1}j^{X}_{n}(b_{i})\langle j^{X}_{n}(b_{i}) | - \rangle_{B_{n+1}}^{X_{n+1}}=\id_{X_{n+1}}$. Since the left hand side is $B_{n+1}$-linear, and $j^{X}_{n}(X_{n})\lhd B_{n+1}=X_{n+1}$ by assumption, it suffices to prove this equality evaluated on the image of $j^{X}_{n}$. We compute
\begin{align*}
    \sum^{m}_{i=1}j^{X}_{n}(b_{i})\langle j^{X}_{n}(b_{i}) | j^{X}_{n}(x)\rangle_{B_{n+1}}^{X_{n+1}} & = \sum^{m}_{i=1}j^{X}_{n}(b_{i})i^{B}_{n}\left(\langle b_{i} | x\rangle_{B_{n}}^{X_n} \right) \\
    & = \sum^{m}_{i=1}j^{X}_{n}(b_{i}(\langle b_{i} | x\rangle_{B_{n}}^{X_n})=j^{X}_{n}(x). 
\end{align*}
\end{proof}

\begin{prop}\label{prop:CompletionFunctor}
There is a natural completion 2-functor $$\varinjlim:\IL\rightarrow \rCorr.$$
\end{prop}
\begin{proof}

The inductive limit $A:=\varinjlim (A_n, i_n, \lambda^A_n)$ is a well-defined familiar $\rm C^*$-algebra, so we will only take care of the higher datum. 
Here, $\{\lambda^A_n: A_n\to \varinjlim A_n\}_n$ are the structure maps. (c.f. \cite{MR1783408})

We shall now describe inductive limit correspondences using an ``increasing union picture''. 
Given $\{X_n\}\in \IL(\{A_n\}\to \{B_n\}),$ consider the product $\prod_{m=1}^\infty X_m$ and the subset $\sum_{m=1}^\infty X_m$ of finitely supported sequences as sets.
Let $\pi: (\prod_{m=1}^\infty X_m)\to (\prod_{m=1}^\infty X_m/\sum_{m=1}^\infty X_m)$ be the quotient map, 
and for each $n\in\bbN$ let $\nu_n: X_n\to \prod_{m=1}^\infty X_m$ defined by $x_n\mapsto (j_{n+k,n}(x_n))_{k\geq -n+1}.$ 
Now, define 
$$\mu_n:= (\pi\circ \nu_n): X_n\to \left(\prod_{m=1}^\infty X_m\Bigg/\sum_{m=1}^\infty X_m\right).$$ 
By direct inspection, we obtain $\mu_n = \mu_{n+1}\circ j_n.$ 
We now define the increasing union $X_0:=\bigcup_{n\geq 1}\mu_n[X_n]$. 
Notice that 
$$\langle \mu_n(x_n)|\ \mu_m(x'_m)\rangle_B:= \lambda^B_{n+m}\left(\left\langle j_{n+m,n}(x_n)|\ j_{n+m,m}(x'_m)\right\rangle^{X_{n+m}}_{B_{n+m}}\right)$$ 
is a well-defined $B$-valued inner product on $X_0$. 
Indeed, this follows since the maps $j_n$ are injective by Remark \ref{rmk:InjConnMaps}. 
We let $X$ be the completion of $X_0$ with respect to the induced $B$-norm, and refer to $\mu_n$ equivalently as $\mu_n^X$ whenever necessary. 

We now equip $X$ with a right $B$-action by bounded operators, and a left $A$-action by right-adjointable operators. 

Let $b_m\in B_m,$ and let $\mu_n(x_n)\in \mu_n[X_n].$ 
We define 
$$x_n\lhd b_m:=\mu_n(x_n)\lhd b_m := \pi\left((j_{n+m+k,n}(x_n)\lhd i_{n+m+k,m}(b_m))_{k\geq (1-m-n)}\right).$$ 
This action is well-defined, and also $||x_n\lhd b_m||_B\leq ||x_n||_B \cdot||b_m||.$
Continuously extending, we see each $B_m$ acts on the right of $X$ by bounded operators. 
To extend this to a right $B$-action on $X$: for $x\in X,$ and $b\in B,$ define $x\lhd b:= \lim_{m\to \infty}x\lhd b_m,$ where $b_m$ is any sequence approximating $b.$ 
This is easily seen to be well-defined, and bounded. 

Using similar arguments and the positivity of the right $B$-valued inner product will provide us the aforementioned commuting left action of $A$ on $X$ by adjointable operators.
Thus, $X:= \varinjlim{ {X_{m}} }$ is an $A$-$B$ correspondence. By Lemma \ref{lem:SamePP}, $\{\varinjlim j^{X}_{n}(b_{i})\}$ is a Pimsner-Popa basis for $\varinjlim X_{n}$, and thus $X\in \rCorr(A\to B)$.

Now let $\{X_n\},\{Y_n\}\in\IL(\{A_n\}\to \{B_n\})$ and $\{f_n\}\in\IL(\{X_n\}\Rightarrow \{Y_n\})$ and we will obtain a map $\varinjlim f_n\in\rCorr(X\Rightarrow Y).$ For $\mu_m(x_m)\in X_0,$ we define 
$$\left(\varinjlim f_n\right)\left(\mu^X_m(x_m)\right):= \mu^Y_{m}(f_{m}(x_m)) = \pi\left(\left(f_{m+k}\left(j^X_{m+k,m}(x_m)\right)\right)_{k\geq -m+1}\right),$$ 
where the last equality follows by Property \ref{IL:Intertwiner}. 
It is clear that for each $k,n\in\bbN$, $\varinjlim f_n$ is $A_k$-$B_m$ bimodular, and by the isometric properties of connecting maps together with the assumption that $\sup ||f_n||<\infty,$ it follows that we can extend $\varinjlim f_m$ to a bounded map on $X.$ 

The resulting map $\varinjlim f_n$ is easily seen to be $A$-$B$ bimodular. 
To verify adjointability, on $Y_0\subset Y:=\varinjlim Y_m$ one can directly define 
$$(\varinjlim f_n)^*(\mu_m^Y(y_m)):=\mu^X_m(f_m^*(y_m)),$$ 
which similarly defines an $A_k$-$B_n$ bimodular bounded map $Y\to X$ for any $k,n\in\bbN.$ 
That these maps are mutual adjoints is a direct computation. 

We now describe the compatibility between $\varinjlim(X_m\boxtimes_{B_m}Y_m)$ and $X\boxtimes_B Y.$ 
For each $m\in\bbN,$ we let $\mu^{X\boxtimes_B Y}_m: X_m\boxtimes_{B_m}Y_m\to \varinjlim(X_n\boxtimes_{B_n}Y_n)$ be the structure map, analogous to those above, and detailed below. 
Now we consider the map
\begin{align*}
    (\varinjlim)^2_{X,Y}:\ \bigcup_m\mu_m^{X\boxtimes_B Y}(X_m\boxtimes_{B_m}Y_m)&\to X\boxtimes_B Y\\
    \mu^{X\boxtimes_B Y}_m(x_m\underset{B_m}{\boxtimes} y_m)=\pi^{X\boxtimes_B Y}\left(\left(j^X_{m+k,m}(x_m)\underset{{B_{m+k}}}{\boxtimes}j^Y_{m+k,m}(y_m)\right)_{k\geq 1-m}\right)&\mapsto \mu^X_m(x_m)\underset{B}{\boxtimes} \mu^Y_m(y_m).
\end{align*}
(Here, $\pi^{X\boxtimes_B Y}$ is the quotient map identifying sequences of tensor products that differ up to a finitely supported sequence, and is analogous to the $\pi$ quotient map from earlier in this proof.) 
The map is well-defined because each $B_k\subset B,$ and it is seen to be isometric and $A_m$-$C_m$ bimodular by inspection.
We can therefore extend this map to an $A$-$C$ bimodular bounded map on all of $\varinjlim(X_m\boxtimes_{B_m}Y_m).$ 
One can directly check unitarity of $(\varinjlim)^2_{X,Y}$ by writing the obvious (densely and well-defined, and bounded) inverse, and the naturality in both variables $X$ and $Y$. 
The remaining coherences for $(\varinjlim,\ \varinjlim^2)$ are also straightforward.
\end{proof}

We will be interested in actions of fusion categories $F:\cC\rightarrow \Bim(A)$ that are built out of functors form $\cC$ to $\IL$. To describe these we introduce some notation.

\begin{rem}\label{rem:ILFunctorData}
For a given unitary tensor functor $F':\cC\to \IL(\{A_n\}\to \{A_n\})$ and $c\in\cC$, 
$F'(c):\{A_n\}\to \{A_n\}$ is a 1-morphism,
so $F'(c)=\{F_n(c)\}$ and $F_n(c)$ is an $A_n-A_n$ bimodule with $j^{c}_n:F_n(c)\to F_{n+1}(c)$ satisfying the conditions \ref{IL:Corr}. 
According to Lemma \ref{lem:StrucEqStrucs}, 
this gives a unitary isomorphism $u_n^c:F_n(c)\boxtimes_{A_{n}}A_{n+1}\xrightarrow{\sim} F_{n+1}(c)$.
Clearly $F_n:\cC\to \Bim(A_n)$ is a unitary tensor functor.

Conversely, for a collection of unitary tensor functors $\{F_n:\cC\to \Bim(A_n)\}$ with the unitary $u_n^c:F_n(c)\boxtimes_{A_{n}}A_{n+1}\xrightarrow{\sim} F_{n+1}(c)$ satisfying obvious coherence conditions,
we can build a unitary tensor functor $\cC\to \IL(\{A_n\}\to \{A_n\})$. Note by Lemma \ref{lem:StrucEqStrucs}, this is equivalent to having a family of linear $j^{c}_{n}:F_{n}(c)\rightarrow F_{n+1}(c)$ satisfying the conditions of Lemma \ref{lem:StrucEqStrucs} and the obvious monoidal coherence.

We will use the notation $\{F_n\}$ to indicate a tensor functor $\cC\to \IL(\{A_n\}\to \{A_n\})$. We denote the resulting composition of $\varinjlim\circ\{F_{n}\}:\cC\rightarrow \Bim(A)$ by $\varinjlim F_{n}$ (where $A=\varinjlim A_{n}$).
\end{rem}

\hspace{.5cm}

One of the main useful features of inductive limit actions is that in principle we can compute invariants --namely $\widehat{F}$-- though their inductive building blocks.

\begin{thm}\label{thm:RealizationContinuous}
Let $\{F_{n}\}: \cC\to\IL(\{A_{n}\}\to \{A_n\})$ be a unitary tensor functor 
and let $Q\in \cC$ be a $Q$-system. 
Then $|Q|_{\varinjlim F_{n}} = \varinjlim |Q|_{F_{n}}$ as inductive limit of $\rm C^*$-algebras. 
\end{thm}
\begin{proof}

Note that by definition for each $m\in\bbN$, $|Q|_{F_m}=F_{m}(Q)$ which unitally embeds in $(\varinjlim F_{m})(Q)=|Q|_{\varinjlim F_{n}}$. Thus we have a unital inclusion of the finite dimensional $\rm C^*$-algebra $\nu_{m}:|Q|_{F_m}\hookrightarrow |Q|_{\varinjlim F_{n}}$, given by $\nu_{m}$ Furthermore, the $j^{Q}_{m}:|Q|_{F_m}\rightarrow |Q|_{F_{m+1}}$ satisfy $\nu_{m+1}\circ j^{Q}_{m}=\nu_{m}$. Thus these assemble into a unique injective $\rm C^*$-algebra homomorphism
$$ \varinjlim |Q|_{F_{n}}\rightarrow |Q|_{(\varinjlim F_{n})}.$$
But this is clearly an equivalence at the level of bimodules hence it must be a bijection, hence a $\rm C^*$-algebra isomorphism.
\end{proof}

We will now show that our invariant is continuous with respect to inductive limits. First, given a functor $\{F_{n}\}: \cC\rightarrow \IL(\{A_{n}\}\to \{A_{n}\})$, we consider the inductive system $\eta_{n}:\widehat{F}_{n}\Rightarrow \widehat{F}_{n+1}$ defined as follows: Recall that in the structure of $\{F_{n}\}$, the family $\{u^{c}_n\}$ (see Remark \ref{rem:ILFunctorData}) equips the $A_n-A_{n+1}$ correspondence $_{i_n} A_{n+1}$ with the structure of a 1-morphism in $\rCorr_{\cC}(F_{n}\to F_{n+1})$. By Corollary \ref{cor:StrEqPtInv}, $- \boxtimes_{A_n} {}_{i_n} A_{n+1}$ uniquely extends to a natural transformation $$\eta_{n}: \widehat{F}_{n}\Rightarrow \widehat{F}_{n+1}.$$

\begin{cor}
\label{Cor:FhatLimitFnhat}
Let $\{F_{n}\}:\cC\rightarrow \IL(\{A_n\}\to \{A_{n}\})$ be as in the statement of the previous theorem. Then $\widehat{\varinjlim F_{n}}\cong \varinjlim \widehat{F}_{n}$, where the latter is the categorical colimit in $\Fun_{+}([\QSys(\cC)]\to\PAbG)$.
\end{cor}
\begin{proof}

We set $F:=\varinjlim F_{n}$. Notice that the inclusion $\lambda_{n}: A_{n}\rightarrow A$ equips the $A_{n}-A$ correspondence $_{\lambda_{n}}A$ with the structure a 1-morphism in $\rCorr(A_{n}\to A)$. Indeed, the isomorphisms $v_{n}^{c}:F_{n}(c)\boxtimes_{A_{n}} {}_{\lambda_{n}}A\cong {}_{\lambda_{n}}A\boxtimes_{A} F(c)$ arise from applying Lemma \ref{lem:StrucEqStrucs} to the inclusions $F_{n}(c)\rightarrow F(c)$ defined in the inductive limit functor. Thus by Corollary \ref{cor:StrEqPtInv} (1), these define natural transformations $\theta_{n}: \widehat{F}_{n}\Rightarrow \widehat{F}.$ Furthermore, since $\lambda_{n}=\lambda_{n+1}\circ i_{n}$, we get 
$$\theta_{n}=\theta_{n+1}\circ \eta_{n}$$
by Corollary \ref{cor:StrEqPtInv} (3). By the universal property of inductive limits, this uniquely defines a natural transformation 
$$\theta: \varinjlim \widehat{F}_{n}\Rightarrow \widehat{F}.$$
But at the level of objects, this is an isomorphism by Theorem \ref{thm:RealizationContinuous}. Thus $\theta$ is a natural isomorphism.
\end{proof}

\subsection{AF-actions}\label{sec:AFActions}
\ \\
We now introduce the main object of study for our paper. 

\begin{defn}\label{defn:AFAction}
An action $F$ of a unitary fusion category $\cC$ on a $\rm C^*$-algebra $A$ is called an \textit{AF-action} if there exists an object $\{A_n\}\in \IL$ with each $A_{n}$ finite dimensional, and a unitary tensor functor $\{F_{n}\}: \mathcal{C}\rightarrow \IL(\{A_n\}\rightarrow \{A_n\})$ such that $F$ is naturally isomorphic to $\varinjlim F_{n}$. 
\end{defn}

We have defined AF-action as a property of an action rather than an additional structure. Our next lemma shows this property is invariant under equivalence.

\begin{lem}\label{lem:AfunderStrong} If $\cC\overset{F}{\curvearrowright} A$ and $\cC\overset{G}{\curvearrowright} B$ are equivalent and $F$ is an AF-action, then $G$ is an AF-action.
\end{lem}
\begin{proof}

If $F$ is AF, then there exists $\{A_{n}\}$ finite dimensional with $\varinjlim A_{n}=A$, and $\{F_{n}\}:\cC\rightarrow \IL(\{A_n\}\to \{A_{n}\})$ with $F$ naturally isomorphic to $H=\varinjlim F_{n}$. Note that $G$ is equivalent to $F$ if and only if $G$ is equivalent to $H$.

Thus there exists an isomorphism $\phi:B\rightarrow A$ such that the tensor functor $H^{\phi}(-):={}_{\phi} A\boxtimes_{A} H(-)\boxtimes_A {}_{\phi^{-1}}B$ is unitarily monoidally naturally isomorphic to $G$. Note that $H^{\phi}(c)\cong H(c)$ as a vector space with $B$ actions $b\rhd x \lhd b^{\prime}=\phi(b)\rhd x \lhd \phi(b^{\prime})$. The right $B$-valued inner product is given by $\langle x|y\rangle_{B}=\phi^{-1}(\langle x|y\rangle_{A})$. 

If we define $B_{n}:=\phi^{-1}(A_{n})$ then that $B=\varinjlim B_{n}$, and $\{F^{\phi}_n\}: \{B_{n}\}\rightarrow \{B_n\}$, then clearly $H^{\phi}=\varinjlim F^{\phi}_{n}$. Thus $G$ is an AF-action.
\end{proof}

AF-actions of the fusion category $\Hilb(G)$ have been studied previously under the guise of inductive limit actions of finite groups \cite{MR1414065, https://doi.org/10.48550/arxiv.1304.0813, MR4206896}. 
There is a way to represent AF-actions which is analogous to Bratteli diagrams for AF $\rm C^*$-algebras \cite{MR312282}.

\begin{defn}\label{defn:EBD}
Let $\cC$ be a unitary fusion category. An \textit{enriched Bratteli diagram} $(\cM_n,I_n,x_n)_{n\ge 0}$ consists of the following data: 
$$(\cM_0,x_0)\xrightarrow{I_0} (\cM_1,x_1)\xrightarrow{I_1} (\cM_2,x_2)\rightarrow\cdots$$
where $\cM_n\in\fssMod(\cC)$, $I_n:\cM_n\to \cM_{n+1}$ is a $\cC$-module functor and $x_n\in\cM_n$ is a generating object such that $I_n(x_n)\cong x_{n+1}$. 
\end{defn}

If $\cC=\Hilb$, then any module category is equivalent to $\Hilb^{\oplus n}$ (with trivial action), 
and any generating object is specified by a positive integer for each $\Hilb$ summand. 
A module functor $\Hilb^{\oplus n}\rightarrow \Hilb^{\oplus m}$ is represented by a bi-partite graph with $n$ source vertices and $m$ target vertices. 
Thus we can recover the usual notion of a Bratteli diagram. 

\hspace{.5cm}

Our claim is that these data characterize a tensor functor $\{F_n\}:\cC\rightarrow \IL(\{A_n\}\to \{A_n\})$ with each $A_{n}$ finite dimensional. 

For a given tensor functor $\{F_n\}:\cC\to \IL(\{A_n\}\to \{A_n\})$, 
let $\cM_n:= \cM_{A_n,F_n}=\Mod(A_n)$ with $\cC$-module structure $X\lhd c:=X\boxtimes_{A_n} F_n(c)$ (see also Proposition \ref{prop:rCorr_CEquivCmod}), and the generating object be ${A_n}_{A_n}$.
Since $A_n$ is finite dimensional, $\cM_n$ is finitely semisimple by Corollary \ref{cor:fdrCorr_CEquivfssCmod}.
Let the connecting $\cC$-module functor $I_n:= -\boxtimes_{A_n} A_{n+1}: \cM_{A_n,F_n}\to \cM_{A_{n+1},F_{n+1}}$,
where the $\cC$-module structure is given by $u_n^c$ such that
$$I_n(X)\lhd c = X\boxtimes_{A_n} A_{n+1}\boxtimes_{A_{n+1}} F_{n+1}(c)\xrightarrow[\sim]{1_X\boxtimes (u_n^c)^*} X\boxtimes_{A_n} F_n(c)\boxtimes_{A_n} A_{n+1} = I_n(X\lhd c),$$
for $X\in \Mod(A_n)$.
It is also clear that for the generating object, 
$$I_n(A_n) = A_n\boxtimes_{A_n} A_{n+1}\cong A_{n+1}.$$

Conversely, given an enriched Bratteli diagram $(\cM_0,x_0)\xrightarrow{I_0} (\cM_1,x_1)\xrightarrow{I_1} (\cM_2,x_2)\rightarrow\cdots$, 
according to Proposition \ref{prop:rCorr_CEquivCmod}, 
we define the finite dimensional $\rm C^*$-algebra $A_n:=\End_{\cM_n}(x_n)$ such that $\cM_n\cong \cM_{A_n,F_n}$ for some unitary tensor functor $F_n:\cC\to \Bim(A)$.
By Remark \ref{rem:functor-bimod},
the $\cC$-module functor $I_n:(\cM_{n},x_n)\to (\cM_{n+1},x_{n+1})$ corresponds to the bimodule ${}_{A_n}{A_{n+1}}_{A_{n+1}}$ with a unital $*$-homomorphism $i_n:A_n\to A_{n+1}$.
The $\cC$-module structure on $I_n$ gives a family of $A_{n}-A_{n+1}$ bimodular unitary isomorphisms 
$$\left\{u^{c}_{n}:F_{n}(c)\boxtimes_{A_{n}} A_{n+1}\rightarrow F_{n+1}(c)\right\}_{c\in\cC}$$
satisfying straightforward coherences (See Definition \ref{defn:rCorrcC}.)
According to Lemma \ref{lem:StrucEqStrucs}, this now gives us $h^{c}_{n}: F_{n}(c)\rightarrow F_{n+1}(c)$ defined by $u^{c}_{n}(x\boxtimes 1_{A_{n+1}})$, which clearly satisfy the required Conditions \ref{IL:Corr} to give $\{F_{n}(c)\}$ the structure of a unitary tensor functor $\{F_n\}:\cC\rightarrow \IL(\{A_n\}\to \{A_n\})$. 

\hspace{.5cm}

By composing the $\varinjlim$ functor on $\IL$, we can say an enriched Bratteli diagram gives an AF-action.
Here, we don't discuss the ``equivalence" between two enriched Bratteli diagrams (or two functors $\cC\to \IL(\{A_n\}\to \{A_n\})$).
In fact, not only will two AF-actions given by ``the same" (under the obvious notion) Bratteli diagrams be equivalent, but we will be able to choose the isomorphisms on the AF $\rm C^*$-algebras in such a way as to preserve the choice of finite dimensional filtration. 
Thus, as in the case for ordinary AF $\rm C^*$-algebras, our enriched Bratteli diagrams are too fine an invariant for classification, but are extremely useful for constructing and representing a given AF-action (see Section \ref{Ex:HilbZ4onZ1/4}).

\subsection{MPO symmetries of spin chains}\label{sec:MPO}
\mbox{}

Recently there has been a great deal of interest from physicists in fusion category symmetries of spin chains. If a critical state on a spin chain has a fusion category $\cC$ of symmetries, then the emergent CFT is expected to realize this fusion category as symmetry defects. This has provided a starting point for the search for exotic quantum field theories \cite{MR3719546, PhysRevLett.128.231602, PhysRevLett.128.231603, PhysRevB.96.125104, 2205.15243}. In addition, fusion categorical symmetry protected topological orders in 1-D have been studied \cite{https://doi.org/10.48550/arxiv.2203.12563}.

Typically fusion category symmetries in the context of physics are described in the language of \textit{tensor networks} (see \cite{Bridgeman_2017}). In particular, fusion categories act on spin chains via \textit{matrix product operators} (MPOs). In this section, we will show that a fusion category action via MPOs yields an AF-action on an AF $\rm C^*$-algebra. Thus our classification result can be interpreted as a classification of general fusion category symmetries on spin chains at the kinematical level (i.e. without states or Hamiltonians). This may be utilized as a first step toward a classification of fusion categorical symmetric states, critical or otherwise.

Below we present a version of the MPO formalism that is distilled from the literature to form a mathematically coherent (and convenient) picture. Our formalism can be found either explicitly or implicitly in nearly all applications of MPO symmetry. The connections with subfactor theory and fusion categories have been pointed out in \cite{MR4109480, 2021arXiv210204562K, 2205.15243}.

Let $\cM$ be an arbitrary monoidal category, and $H\in \cM$ be an object. We define a new monoidal category $\cM_{H}$ as follows:
\begin{itemize}
    \item 
    Objects are pairs $(K, \rho)$, where $K\in \cM$ and $\rho: K\otimes H\rightarrow H\otimes K$ is a (unitary) isomorphism.
    \item
    $\cM_{H}(\ (K,\rho)\to (L,\theta)\ )=\{f\in \cM(K\to L)\ :\ (1_{H}\otimes f)\rho=\theta(f\otimes 1_{H})\}$
    \item
    $(K, \rho)\otimes (L,\theta):=(K\otimes L, (\rho\otimes 1_{L})\circ (1_{K}\otimes \theta))$
\end{itemize}

It is straightforward to check that $\cM_{H}$ assembles into a monoidal category.

Now, let $V$ be a finite set. Then $\BigHilb^{V}$ is the monoidal category of finite dimensional Hilbert spaces bi-graded by a fixed finite set $V$. 
More specifically:  

\begin{itemize}
\item 
Objects are finite dimensional Hilbert spaces with direct sum decomposition $H=\bigoplus_{a,b\in V} H_{ab}$.
\item
Morphisms are linear operators between Hilbert spaces that preserve the bigrading.
\item
The monoidal product $H\otimes_{V} K$ is the bigraded Hilbert space with components $(H\otimes_{V} K)_{ab}:=\bigoplus_{c\in V}  H_{ac}\otimes K_{cb} $,
where the later tensor product is the tensor product of Hilbert spaces. 

We say a $V$ bi-graded Hilbert space is non-degenerate if for each $a$, there is a $b$ with $H_{ab}\ne 0$ and for each $b$, there is a $a$ with $H_{ab}\ne 0$.
\end{itemize}

\begin{defn} Let $H$ be finite dimensional Hilbert spaces bi-graded by a finite set $V$. We define the monoidal category of \textit{matrix product operators} by $\MPO(H, V):=\BigHilb^{V}_{H}$.
\end{defn}

If $(K,\rho)\in \MPO(H, V)$, we call $H$ the \textit{physical Hilbert space}, $K$ the virtual Hilbert space, and $\rho$ the matrix product operator. Our definition is somewhat non-standard, so we will unpack it for the convenience of the reader. Our unpacked version is easily seen to be equivalent to the typical versions that exist in the literature.

If $(K,\rho)\in \MPO(H, V)$, then
\begin{itemize}
\item
$K=\bigoplus_{a,b\in V} K_{ab}$ has a bigrading
\item 
$\rho: K\otimes_{V} H\rightarrow H\otimes_{V} K$ is a bigrading preserving morphism such that $\rho^{*}\rho=1_{K\otimes_{V} H}$ and $\rho\rho^{*}=1_{H\otimes_{V} K}$.
\end{itemize}
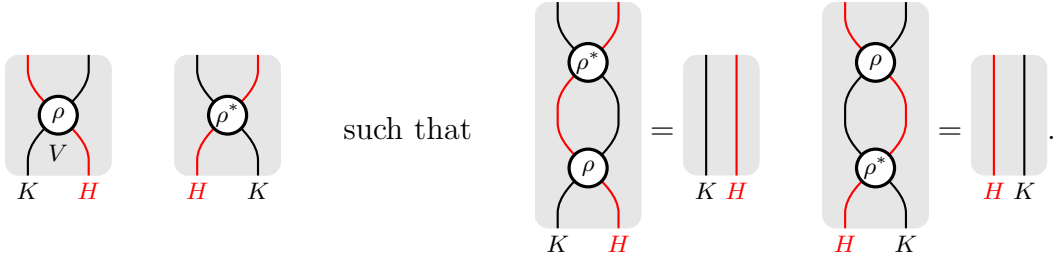
\[
\tikzmath{
\begin{scope}
\clip[rounded corners = 5] (-.7,-.8) rectangle (.7,.8);
\filldraw[\rColor] (-.7,-.8) rectangle (.7,.8);
\end{scope}
\draw[thick] (-.4,-.8) -- (-.4,-.6) .. controls ++(90:.45cm) and ++(270:.45cm) .. (.4,.6) -- (.4,.8);
\draw[\HColor,thick] (.4,-.8) -- (.4,-.6) .. controls ++(90:.45cm) and ++(270:.45cm) .. (-.4,.6) -- (-.4,.8);
\roundRbox{unshaded}{(0,0)}{.25}{\scriptsize{$\rho$}};
\node at (-.4,-1) {\scriptsize{$K$}};
\node[\HColor] at (.4,-1) {\scriptsize{$H$}};
\node at (0,-.5) {\scriptsize{$V$}};
}
\qquad
\tikzmath{
\begin{scope}
\clip[rounded corners = 5] (-.7,-.8) rectangle (.7,.8);
\filldraw[\rColor] (-.7,-.8) rectangle (.7,.8);
\end{scope}
\draw[\HColor,thick] (-.4,-.8) -- (-.4,-.6) .. controls ++(90:.45cm) and ++(270:.45cm) .. (.4,.6) -- (.4,.8);
\draw[thick] (.4,-.8) -- (.4,-.6) .. controls ++(90:.45cm) and ++(270:.45cm) .. (-.4,.6) -- (-.4,.8);
\roundRbox{unshaded}{(0,0)}{.25}{\scriptsize{$\rho^*$}};
\node at (.4,-1) {\scriptsize{$K$}};
\node[\HColor] at (-.4,-1) {\scriptsize{$H$}};
}
\qquad\text{such that}\qquad
\tikzmath{
\begin{scope}
\clip[rounded corners = 5] (-.7,-1.5) rectangle (.7,1.5);
\filldraw[\rColor] (-.7,-1.5) rectangle (.7,1.5);
\end{scope}
\draw[thick] (-.4,-1.5) -- (-.4,-1.3) .. controls ++(90:.45cm) and ++(270:.45cm) .. (.4,-.1) -- (.4,.1) .. controls ++(90:.45cm) and ++(270:.45cm) .. (-.4,1.3) -- (-.4,1.5);
\draw[\HColor,thick] (.4,-1.5) -- (.4,-1.3) .. controls ++(90:.45cm) and ++(270:.45cm) .. (-.4,-.1) -- (-.4,.1) .. controls ++(90:.45cm) and ++(270:.45cm) .. (.4,1.3) -- (.4,1.5);
\roundRbox{unshaded}{(0,-.7)}{.25}{\scriptsize{$\rho$}};
\roundRbox{unshaded}{(0,.7)}{.25}{\scriptsize{$\rho^*$}};
\node at (-.4,-1.7) {\scriptsize{$K$}};
\node[\HColor] at (.4,-1.7) {\scriptsize{$H$}};
}
=
\tikzmath{
\begin{scope}
\clip[rounded corners = 5] (-.5,-.8) rectangle (.5,.8);
\filldraw[\rColor] (-.5,-1.5) rectangle (.5,1.5);
\end{scope}
\draw[thick] (-.2,-.8) -- (-.2,.8);
\draw[\HColor,thick] (.2,-.8) -- (.2,.8);
\node at (-.2,-1) {\scriptsize{$K$}};
\node[\HColor] at (.2,-1) {\scriptsize{$H$}};
}
\qquad
\tikzmath{
\begin{scope}
\clip[rounded corners = 5] (-.7,-1.5) rectangle (.7,1.5);
\filldraw[\rColor] (-.7,-1.5) rectangle (.7,1.5);
\end{scope}
\draw[\HColor,thick] (-.4,-1.5) -- (-.4,-1.3) .. controls ++(90:.45cm) and ++(270:.45cm) .. (.4,-.1) -- (.4,.1) .. controls ++(90:.45cm) and ++(270:.45cm) .. (-.4,1.3) -- (-.4,1.5);
\draw[thick] (.4,-1.5) -- (.4,-1.3) .. controls ++(90:.45cm) and ++(270:.45cm) .. (-.4,-.1) -- (-.4,.1) .. controls ++(90:.45cm) and ++(270:.45cm) .. (.4,1.3) -- (.4,1.5);
\roundRbox{unshaded}{(0,-.7)}{.25}{\scriptsize{$\rho^*$}};
\roundRbox{unshaded}{(0,.7)}{.25}{\scriptsize{$\rho$}};
\node at (.4,-1.7) {\scriptsize{$K$}};
\node[\HColor] at (-.4,-1.7) {\scriptsize{$H$}};
}
=
\tikzmath{
\begin{scope}
\clip[rounded corners = 5] (-.5,-.8) rectangle (.5,.8);
\filldraw[\rColor] (-.5,-1.5) rectangle (.5,1.5);
\end{scope}
\draw[\HColor,thick] (-.2,-.8) -- (-.2,.8);
\draw[thick] (.2,-.8) -- (.2,.8);
\node at (.2,-1) {\scriptsize{$K$}};
\node[\HColor] at (-.2,-1) {\scriptsize{$H$}};
}.
\]

Note that $\rho_{ab}:(K\otimes_V H)_{ab}\to (H\otimes_V K)_{ab}$.
Let $\rho^c_{ab}: K_{ac}\otimes H_{cb}\hookrightarrow (K\otimes_V H)_{ab} \to (H\otimes_V K)_{ab}$, 
and $(\rho^c_{ab})^*: (H\otimes_V K)_{ab}\to (K\otimes_V H)_{ab} \twoheadrightarrow K_{ac}\otimes H_{cb}$,
we have
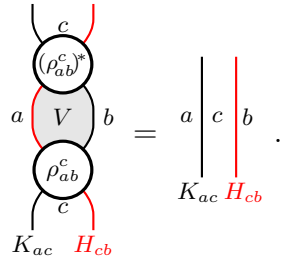
\[
\tikzmath{
\begin{scope}
\clip[rounded corners = 5] (-.7,-1.5) rectangle (.7,1.5);
\filldraw[\rColor] (-.4,.1) -- (-.4,-.1) .. controls ++(270:.3cm) and ++(135:.15cm) .. (0,-.7) .. controls ++(45:.15cm) and ++(270:.3cm) .. (.4,-.1) -- (.4,.1).. controls ++(90:.3cm) and ++(-45:.15cm) .. (0,.7) .. controls ++(-135:.15cm) and ++(90:.3cm) .. (-.4,.1);
\end{scope}
\draw[thick] (-.4,-1.5) -- (-.4,-1.3) .. controls ++(90:.45cm) and ++(270:.45cm) .. (.4,-.1) -- (.4,.1) .. controls ++(90:.45cm) and ++(270:.45cm) .. (-.4,1.3) -- (-.4,1.5);
\draw[\HColor,thick] (.4,-1.5) -- (.4,-1.3) .. controls ++(90:.45cm) and ++(270:.45cm) .. (-.4,-.1) -- (-.4,.1) .. controls ++(90:.45cm) and ++(270:.45cm) .. (.4,1.3) -- (.4,1.5);
\roundRbox{unshaded}{(0,-.7)}{.38}{\scriptsize{$\rho_{ab}^c$}};
\roundRbox{unshaded}{(0,.7)}{.38}{\tiny{$( \rho_{\! ab}^{c}\!)^{\! *}$}};
\node at (-.4,-1.7) {\scriptsize{$K_{ac}$}};
\node[\HColor] at (.4,-1.7) {\scriptsize{$H_{cb}$}};
\node at (0,0) {\scriptsize{$V$}};
\node at (0,-1.2) {\scriptsize{$c$}};
\node at (0,1.2) {\scriptsize{$c$}};
\node at (-.6,0) {\scriptsize{$a$}};
\node at (.6,0) {\scriptsize{$b$}};
}
=
\tikzmath{
\draw[thick] (-.2,-.8) -- (-.2,.8);
\draw[\HColor,thick] (.25,-.8) -- (.25,.8);
\node at (-.25,-1) {\scriptsize{$K_{ac}$}};
\node[\HColor] at (.35,-1) {\scriptsize{$H_{cb}$}};
\node at (0,0) {\scriptsize{$c$}};
\node at (-.4,0) {\scriptsize{$a$}};
\node at (.4,0) {\scriptsize{$b$}};
}.
\]

The main difference between our version of MPOs and the ones typically found in the physics literature, is that the bi-grading on the physical and virtual bond spaces are explicit in our formalism, whereas they are typically implicit in the physical formulation. Taking products over the bi-grading is strictly necessary in order to actually obtain fusion category symmetry, and is done implicitly in the literature. We have chosen to make this explicit for the convenience of giving formal definitions of fusion category symmetry. 

\begin{defn}
Let $\cC$ be a fusion category. Then an MPO-action of $\cC$ on a spin chain with local Hilbert space $H$ consists of a (unitary) monoidal functor $\cC\rightarrow \MPO(H, V)$ for some bi-grading $V$ on $H$.
\end{defn}

We note that the actual physics of the system is described by much more that just the physical local Hilbert space $H$. We also need a \textit{state} for the system. Traditionally this is accomplished by describing a local Hamiltonian, and taking either ground states or equilibrium states for the resulting interaction. Alternate approaches using tensor networks have recently gained traction. In either case, to obtain fusion category symmetry of the low energy effective field theory, we would need to consider the notion of compatibility of states with MPOs (see for example \cite{https://doi.org/10.48550/arxiv.2203.12563}).

In this paper, we take an abstract look at fusion category actions on spin chains. The definition given above using MPOs is a natural extension of ``on-site" symmetries of group actions. We propose that this sort of symmetry is really a special case of AF-actions of fusion categories, and that the latter notion describes a more robust, ``coordinate free" approach to fusion category actions on spin chains that allows for spatial inhomogeneity. Our main result can then be interpreted as a classification of fusion category actions on spin chains in terms of K-theoretic invariants, and can be used as a starting point for a systematic exploration of the phase space of fusion category symmetries of states.

In the remainder of this section, we will show that MPO-actions of fusion categories on spin chains are really the same data as a homogeneous AF-action. The first step is to define the appropriate AF $\rm C^*$-algebra. In the usual operator algebraic approach to spin systems \cite{MR1441540}, if the local Hilbert space is given by $H$, then we would usually describe the $n$-site Hilbert space $H^{\otimes n}$, and thus the algebra of operators localized at these $n$-sites would be $B(H^{\otimes n})=B(H)^{\otimes n}$. Then the quasi-local algebra of observables for the system would be the UHF $\rm C^*$-algebra obtained by taking the inductive limit over all finite sites.

However, we know that this algebra cannot be used to host fusion category symmetries in general. There are basic K-theoretic obstructions (e.g. \cite[\S 5.1]{MR4419534}). At the MPO level, this is reflected by the fact that we have to work with $V$-bigraded Hilbert spaces, rather than ordinary Hilbert spaces, to implement categorical symmetry. In our approach, we will use the quasi-local algebra for the half-infinite lattice. 

In particular, the $n$-site Hilbert space should initially be given by $H^{\otimes^{n}_{V}}:=H\otimes _{V} H\otimes_V \cdots \otimes_{V} H$,
i.e. the $n$-fold tensor power internal to the category $\BigHilb^{V}$. 
The first obvious candidates for how to define the local algebra arises by simply taking $B(H^{\otimes^{n}_{V}})$.
This suffers from the fact that there is no natural inclusion $B(H^{\otimes^{n}_{V}})\hookrightarrow B(H^{\otimes^{n+1}_{V}})$ and so we cannot form an inductive limit algebra. 
Intuitively, what is going on is that the bigrading on $H^{\otimes^{n}_{V}}$ should be though of as left and right boundary conditions. 
In order to be consistent when we extend our observables to larger regions, 
we should \textit{preserve} boundary conditions.
Since we are modeling a half-infinite spin chain (extending to the right), 
we only need to ask our observables to preserve the right boundary condition, 
while we fix a left boundary condition. A left boundary condition consists of a right (singly) $V$-graded finite dimensional Hilbert space $H^{0}=\bigoplus_{a\in V} H^0_a\in\Hilb(V)$.

Then we define
$$A_{n}:= \{ f\in B(H^{0}\otimes_{V} H^{\otimes^{n}_{V}})\ :\ f( (H^{0}\otimes_{V}H^{\otimes^{n}_{V}})_{b}\subseteq (H^{0}\otimes_{V} H^{\otimes^{n}_{V}})_{b}\ \text{for all}\ b\in V\}$$
where $H^0\otimes_V L = \bigoplus_{a,b\in V} H^0_a \otimes L_{ab}$ and $(H^0\otimes_V L)_b = \bigoplus_{a\in V} H^0_a \otimes L_{ab}$ for $L\in \BigHilb^V$. 

Then $A_{n}$ is a finite dimensional $\rm C^*$-algebra, with
$$A_{n}\cong \bigoplus_{b\in V} B((H^0\otimes_V H^{\otimes^{n}_{V}})_{b})$$

We have a natural unital inclusion $A_{n}\hookrightarrow A_{n+1}$ given by $f\mapsto f\otimes_{V} 1_{H}$. Then taking $A_{\infty}:=\varinjlim A_{n}$ gives an AF $\rm C^*$-algebra which we call the quasi-local $\rm C^*$-algebra. 

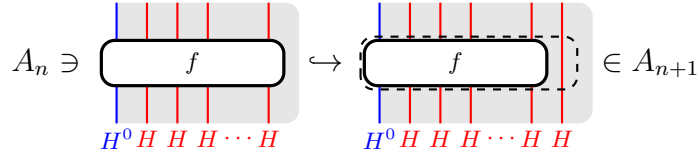
\[
A_n\ni
\begin{tikzpicture}[baseline=-.1cm]
\begin{scope}
\clip[rounded corners = 5] (-.5,-.8) rectangle (2.2,.8);
\filldraw[\rColor] (-.2,-1.5) rectangle (2.2,1.5);
\end{scope}
\draw[\YColor,thick] (-.2,-.8) -- (-.2,.8);
\draw[\HColor,thick] (.2,-.8) -- (.2,.8);
\draw[\HColor,thick] (.6,-.8) -- (.6,.8);
\draw[\HColor,thick] (1,-.8) -- (1,.8);
\draw[\HColor,thick] (1.8,-.8) -- (1.8,.8);
\roundNbox{unshaded}{(.8,0)}{.3}{.9}{.9}{\scriptsize{$f$}};
\node[\YColor] at (-.2,-1) {\scriptsize{$H^0$}};
\node[\HColor] at (.2,-1) {\scriptsize{$H$}};
\node[\HColor] at (.6,-1) {\scriptsize{$H$}};
\node[\HColor] at (1,-1) {\scriptsize{$H$}};
\node[\HColor] at (1.4,-1) {\scriptsize{$\cdots$}};
\node[\HColor] at (1.8,-1) {\scriptsize{$H$}};
\end{tikzpicture}
\hookrightarrow
\begin{tikzpicture}[baseline=-.1cm]
\begin{scope}
\clip[rounded corners = 5] (-.5,-.8) rectangle (2.6,.8);
\filldraw[\rColor] (-.2,-1.5) rectangle (2.6,1.5);
\end{scope}
\draw[\YColor,thick] (-.2,-.8) -- (-.2,.8);
\draw[\HColor,thick] (.2,-.8) -- (.2,.8);
\draw[\HColor,thick] (.6,-.8) -- (.6,.8);
\draw[\HColor,thick] (1,-.8) -- (1,.8);
\draw[\HColor,thick] (1.8,-.8) -- (1.8,.8);
\draw[\HColor,thick] (2.2,-.8) -- (2.2,.8);
\roundNbox{unshaded}{(.8,0)}{.3}{.9}{.9}{\scriptsize{$f$}};
\draw[thick,dashed,rounded corners = 5] (-.45,-.35) rectangle (2.4,.35);
\node[\YColor] at (-.2,-1) {\scriptsize{$H^0$}};
\node[\HColor] at (.2,-1) {\scriptsize{$H$}};
\node[\HColor] at (.6,-1) {\scriptsize{$H$}};
\node[\HColor] at (1,-1) {\scriptsize{$H$}};
\node[\HColor] at (1.4,-1) {\scriptsize{$\cdots$}};
\node[\HColor] at (1.8,-1) {\scriptsize{$H$}};
\node[\HColor] at (2.2,-1) {\scriptsize{$H$}};
\end{tikzpicture}
\in A_{n+1}
\]

We say that an enriched Bratteli diagram is \textit{homogeneous} if $\cM_n\cong\cM_{n+1}\cong\cM$ and $I_n\cong I_{n+1}\cong I$. 
The examples we have in Section \ref{Ex:HilbZ4onZ1/4} are all homogeneous.

\begin{thm}
An MPO-action of a fusion category on a spin chain (up to natural isomorphism) together with a left boundary condition $H^{0}$ is the same data as a homogeneous enriched Bratteli diagram.
\end{thm}

\begin{proof}
Given an MPO-action $\cC\to \MPO(H,V)$. Each $c\in \cC$ has an MPO $(K^{c}, \rho^{c})$, where $K^{c}$ is a $V$-bigraded Hilbert space making $\Hilb(V)$ into a $\cC$-module category $U\lhd c = U\otimes_V K^c$,
and $\rho^{c}$ assemble to make the functor $-\otimes_{V} H$ into a $\cC$-module functor. 
The left boundary condition gives the generating object at level zero.

Conversely, given a homogeneous enriched Bratteli diagram $(\cM,I,x_0)$. 
Let $V:=\Irr(\cM)$ be the finite set of isomorphism classes of simple objects in $\cM$.
Let $H^0$ be the $V$-graded Hilbert space, where the $a$-th component Hilbert space $H^0_a:= \cM(a\to x_0)$, $a\in V$. 
Let $H_{ab}$ be the $V$-bigraded Hilbert space, where the $(a,b)$-th component Hilbert space $H_{ab}:= \cM(a\to I(b))$.
For $c\in \cC$, let $K^c$ be a $V$-bigraded Hilbert space, where the $(a,b)$-th component Hilbert space $K^c_{ab}:=\cM(a\to b\lhd c)$.
In the end, $\rho^c: K^c\otimes_V H = \cM(a\to I(b\lhd c))\xrightarrow{\sim} \cM(a\to I(b)\lhd c) = H\otimes_V K^c$ by using $I(-\lhd c)\cong I(-)\lhd c$.
\end{proof}

Thus in particular, MPO-actions give rise to AF-actions of fusion categories on AF $\rm C^*$-algebras. From the proof above, it is also easy to see that we can interpret non-homogeneous Bratteli diagrams as a kind of MPO-action which is spatially non-homogeneous. Considering these kinds of actions allows us to explore a larger universality class of possible states.  

\begin{rem} The notion of equivalence of actions allows categorical symmetries to be conjugate up to an arbitrary automorphism, which need not reflect spatial structure of the lattice. This is too weak of an equivalence relation from a physical perspective. Indeed by \cite{kishimoto_ozawa_sakai_2003}, automorphisms of $\rm C^*$-algebras act transitively on pure states, so generic automorphisms do not remember even coarse physical information such as the dimensionality of the underlying lattice. Thus it would be interesting to classify AF-actions up to a ``locality preserving'' version of equivalence. 
\end{rem}

\subsection{Classification}
\mbox{}

In this section, we show that the pointed invariant $\left(\widehat{F}, [A]\right)$ is a complete invariant for AF-actions. A key ingredient in Elliott's original proof for AF algebras is that any positive, pointed homomorphism between $K_{0}$ invariants of finite dimensional $\rm C^*$-algebras arises from a unital homomorphism between the algebras, and furthermore any two such homomorphisms are unitarily conjugate \cite[\S 4]{MR397420}. This type of ``existence and uniqueness of morphisms" is a crucial in $\rm C^*$-algebra classification programs. The next lemma provides an analogue of this result in our setting, which is essentially a converse to Corollary \ref{cor:StrEqPtInv} in the finite dimensional setting.

\begin{lem}\label{Lem:Lifting}
Suppose $A$ and $B$ are finite dimensional $\rm C^*$-algebras with actions $\cC\overset{F}{\curvearrowright} A$ and $\cC\overset{G}{\curvearrowright} B$. Let $\theta: \left(\widehat{F},\ [A]\right)\Rightarrow \left(\widehat{G},\  [B]\right)$ be a pointed natural transformation. 
Then there exists a unital $*$-homomorphism $\phi: A\rightarrow B$ together with a $\cC$-equivariant structure $\{h^c\}$ on ${}_\phi B\in\rCorr(A\to B)$ such that $\theta$ is implemented by $- \boxtimes_{A} {}_{\phi} B$. 
Furthermore, $\left({}_{\phi} B,\ \{h^c\}\right)\in \rCorr_{\cC}(F\to G) $ is unique up to unitary isomorphism.
\end{lem}
\begin{proof}
By Theorem \ref{thm:EquivariantCrossedprod}, we have that 
\begin{align*}
&\widehat{F}\cong [\Mod(\cC)(\Psi(-)\to \cM_{A,F})] \quad\qquad \text{ and } && \widehat{G}\cong [\Mod(\cC)(\Psi(-)\to \cM_{B,G})]
\end{align*}
as functors. 

Since $\Psi: [\QSys(\cC)]\rightarrow [\fssMod(\cC)]$ is a 1-contravariant equivalence onto $\fssMod(\cC)$ (see Remark \ref{rem:unitaryostrik}), and $\cM_{A,F}$ and $\cM_{B,G}$ are objects in $\fssMod(\cC)$, then by the Yoneda Lemma, there is a unique morphism $[H]\in [\fssMod(\cC)(\cM_{A,F}\to \cM_{B,G})]=[\Mod(\cC)(\cM_{A,F}\to \cM_{B,G})]$ such that $\theta$ is given by post composition with $H$.
By Proposition \ref{prop:rCorr_CEquivCmod}, $[H]$ uniquely determines an isomorphism class of correspondence $_{A} X_{B}\in \rCorr_{\cC}(A\to B)$ such that $X_{B}\cong H(A)$ (as right Hilbert $B$-modules).
But since $\theta$ is pointed, $H(A)\cong B$ as a right $B$-module, and thus the correspondence defines a unital $*$-homomorphism $\phi: A\rightarrow B$ so that $\theta \cong - \boxtimes_{A} {}_{\phi} B$ (see Definition \ref{defn:hombimod} and proceeding discussion). 
\end{proof}

\begin{lem}\label{lem:NTtolimfinite}
Let $\cD$ be a category enriched in abelian monoids such that there are finitely many isomorphism classes of indecomposable objects, and every object can be written as a direct sum of indecomposables. Let $H_{n}\in \Fun_{+}(\cD\to \PAbG)$ be an inductive system with $\nu_{n}: H_n\Rightarrow H_{n+1}$ injective natural transformations, with $\lambda_n: H_n\Rightarrow H:=\varinjlim H_{n}$.  
Let $E\in \Fun_{+}(\cD \to \PAbG)$ be such that $E(Q)$ is a free abelian group of finite rank for all $Q\in \cD$. Then for any natural transformation $\alpha: E\Rightarrow H$, there exists $m\in\bbN$ and a natural transformation $\alpha^{\prime}:E\Rightarrow H_{m}$ with $\alpha=\lambda_{m}\circ \alpha^{\prime}$.
\end{lem}
\begin{proof}
Let $\InDec(\cD)$ be a set of representatives of isomorphism classes of indecomposable objects in $\cD$.
Then for each $P\in \InDec(\cD)$, there is an $m_{P}$ such that $\alpha_P[E(P)]\subseteq \lambda_{m_p}(H_{m_P}(P))$.
Indeed since $E(P)$ is a finite rank free abelian group, $(\alpha_P)|_{E(P)}$ is completely determined by where it sends generators, each of which occurs at a finite level.
Define $m:=\max_{P\in \InDec(\cD)} \{m_{P}\}$. 

Now let $Q\in\cD$, and pick a decomposition $Q\cong \bigoplus_{i=1}^N P_i,$ where each $P_i\in \InDec(\cD)$ 
and $\pi_i: Q\to P_i$ and $\iota_{i}:P_{i}\rightarrow Q$ in $\cD$ such that 
$\pi_{j}\iota_{i}=\delta_{i=j} \id_{P_i}$
and $\sum_{i}\iota_{i}\pi_{i}=\id_{Q}$ from an direct sum decomposition of $Q$ in $\cD.$ 
Note $E(\iota_{i}): E(P_i)\rightarrow E(Q)$ is injective, and we have $E(Q)\cong \bigoplus_{i} E(P_{i})$ 
where the isomorphism from right to left is given by $(x_{1}, \dots , x_{n})\mapsto \sum^{n}_{i=1}E(\iota_{i})(x_{i})$. In particular, every $y\in E(Q)$ can be written as
$$y=\sum_{i} E(\iota_{i})(x_{i})$$
for uniquely specified $x_{i}\in E(P_{i})$.
Then we have 
$$\alpha_{Q}(y)=\sum_{i} \alpha_{Q}E(\iota_i)(x_{i})=H(\iota_{i})\alpha_{P_i}(x_{i}).$$
But $\alpha_{P_i}(x_i)\in \lambda_{m}(H_{m}(P_{i}))$ and by the definition of inductive limit, for any element $h\in \lambda_{m}(H_{m}(P_{i}))$, 
$$H(\iota_{i})(h)= \lambda_{m}(H_{m}(\iota_{i})(h)).$$ 
Thus $\alpha_{Q}(y)\in H_{m}(Q)$, hence $\alpha_{Q}(E(Q))\subseteq \lambda_{m} H_{m}(Q)$ as desired.

We've show that there exists an $m$ such that for every object $Q\in \cD$, $\alpha\circ E(Q)\subseteq \lambda_{m}(H_{m}(Q))$. Now by the definition of inductive limit functor, we see that for each morphism $f\in \cD(P\to Q)$, $\alpha\circ E(f)=H(f)\circ \alpha=\lambda_{m}(H_{m}(f))\circ \alpha$. Since $\lambda_{m}: H_{m}\rightarrow H$ is an isomorphism onto its image, it has a left inverse $\lambda^{-1}_{m}$ defined on its image. Thus the natural transformation 
$$\alpha^{\prime}:=\lambda^{-1}_{m}\circ \alpha$$
is well defined and satisfies the conclusion of the lemma.
\end{proof}

Now we arrive at the main result of our paper. At its core, the proof is precisely Elliott's intertwining argument. The purpose of the paper up to now has been to develop the language and tools required to apply it on setting (particularly Lemma \ref{lem:StrucEqStrucs} and Lemma \ref{Lem:Lifting}).

\begin{thm*}[Theorem \ref{thmalp:main}]
Suppose $\cC\overset{F}{\curvearrowright} A$ and $\cC\overset{G}{\curvearrowright} B$ are AF-actions. Then $\left(\widehat{F}, [A]\right)\cong \left(\widehat{G}, [B]\right)$ if and only if $\cC\overset{F}{\curvearrowright} A$ is equivalent to $\cC\overset{G}{\curvearrowright} B$. 
\end{thm*}
\begin{proof}
We only need to prove the forward direction, as the converse was established in Corollary \ref{cor:StrEqPtInv}. 

Let $\cC\overset{F}{\curvearrowright} A$ and $\cC\overset{G}{\curvearrowright} B$ be AF-actions and suppose $\alpha:\widehat{F}\Rightarrow \widehat{G}$ is a natural isomorphism of functors, with $\alpha_{1_{\cC}}([A])=[B].$
Moreover, let $\{F_n\}:\cC\to\IL(\{A_n\}\to\{A_n\})$ and $\{G_n\}:\cC\to\IL(\{B_n\}\to\{B_n\})$ be choices of finite-dimensional approximations for $F$ and $G,$ respectively, as in Definition \ref{defn:AFAction}. 
Without loss of generality, we can assume $F=\varinjlim F_{n}$ and $G=\varinjlim G_{n}$.  

If $\nu^{F}_1:F_1\Rightarrow \varinjlim F_n$ denotes the canonical natural transformation, then $\alpha\circ \nu^{F}_1:\widehat{F}_1\Rightarrow \widehat{G}=\varinjlim \widehat{G}_{n}$. 
By Lemma \ref{lem:NTtolimfinite}, there exists $m_{1}$ and a natural transformation $\alpha_{1}: \widehat{F}_{1}\Rightarrow \widehat{G}_{m_1}$ such that 
$$\alpha \circ \nu^F_1 =\nu^{G}_{m_1}\circ \alpha_{1}.$$ 

Applying the same argument with $\beta:=\alpha^{-1}$ gives an $n_{2}>1$ (this strict inequality can be achieved without loss of generality, simply by increasing $n_2$ if necessary and composing with $\nu^{F}$) and a natural transformation $\beta_{1}: \widehat{G}_{m_1}\Rightarrow \widehat{F}_{n_2}$ with 
$$\beta\circ \nu^G_{m_1}=\nu^{F}_{n_2}\circ \beta_{1}.$$

Repeating this procedure and setting $n_{1}=1$, we obtain the following commuting diagram of natural transformations:
\[
 \begin{tikzcd}[column sep=3.5em]
 \widehat{F}_{n_1}
 \arrow[swap]{d}{\alpha_{1}}
 \arrow{r}{\eta^{F}_{n_{2},n_{1}}}
 &  \widehat{F}_{n_2}
 \arrow[swap]{d}{\alpha_{2}}
 \arrow{r}{\eta^{F}_{n_{3},n_{2}}}
 & \cdots
 \arrow{r}{}
 &  \widehat{F}
 \arrow[shift left]{d}{\alpha} \\
  \widehat{G}_{m_1}
 \arrow[swap]{r}{\eta^{G}_{m_2,m_1}}
 \arrow{ru}{\beta_{1}}
 &  \widehat{G}_{m_2}
 \arrow[swap]{r}{\eta^{G}_{m_3, m_2}}
 \arrow{ru}{\beta_{2}}
 & \cdots
 \arrow{r}{}
 &  \widehat{G}
 \arrow[shift left]{u}{\beta}
 \end{tikzcd}.
\]
For notational convenience, we will re-index the subsequences $(n_1, n_{2},\ \dots )$ and $(m_1, m_2,\ \dots)$ by $(1,2,\ \dots )$.

Let us denote by $\{i^A_n,\lambda^A_n\}$ and $\{i^B_n,\lambda^B_n\}$ the inductive limit maps for $A$ and $B$ respectively. By Lemma \ref{Lem:Lifting}, for each $n\in\bbN$ there exist unital $*$-homomorphisms $\phi_{n}:A_{n}\rightarrow B_{n}$ and $\psi_{n}:B_{n}\rightarrow A_{n+1}$, and $\cC$-equivariant structures $\{h_{n}^{c}\}$ and $\{k^{c}_{n}\}$ for $\phi_{n}$ an $\psi_{n}$ respectively, such that the $1$-morphisms $({}_{\phi_{n}} B_{n}, \{h^{n}_{c}\})$ and $({}_{\psi_{n}} A_{n+1}, \{k^{n}_{c}\})$ in $\rCorr_{\cC}$ implement $\alpha_{n}$ and $\beta_{n}$ respectively.

By commutativity of the zig-zag above, $({}_{\psi_{1}\circ\phi_{1}} A_{2}, k^{c}_{1}\circ h^{c}_{1})$ and $({}_{i^{A}_{1}} A_{2}, j^{F,c}_{1})$ implement the same natural transformation $\eta^{F}_{1}:\widehat{F}_{1}\Rightarrow \widehat{F}_{2}$ and thus by Lemma \ref{Lem:Lifting}, they are isomorphic in $\rCorr_{\cC}$. But by Lemma \ref{lem:StrucEqStrucs}, there exists a unitary $v_{1}\in A_{2}$ with  
$$\Ad(v_{1})\circ (\psi_{1}\circ\phi_{1})=i^{A}_{1},$$ 
and for all $c\in \cC$
$$\Ad(v_{1})\circ (k^{c}_{1}\circ h^{c}_{1})=j^{F,c}_{1}$$
Thus we can replace $(\psi_{1}, k^{c}_{1})$ with the isomorphic equivariant correspondence $(\Ad(v_1)\circ \psi_{1}, \Ad(v_1)\circ k^{c}_{1})$.

We continue this procedure inductively to obtain a zig-zag of $\cC$-equivariant homomorphisms
\[
 \begin{tikzcd}[column sep=3.5em]
 A_{1}
 \arrow[swap]{d}{\phi_{1}}
 \arrow{r}{i^{A}_{1}}
 &  A_{2}
 \arrow[swap]{d}{\phi_{2}}
 \arrow{r}{i^{A}_{2}}
 & \cdots
 \arrow{r}{}
 &  A
 \\
 B_1
 \arrow[swap]{r}{i^{B}_{1}}
 \arrow{ru}{\psi_{1}}
 &  B_2
 \arrow[swap]{r}{i^{B}_{2}}
 \arrow{ru}{\psi_{2}}
 & \cdots
 \arrow{r}{}
 &  B 
 \end{tikzcd},
\]
where the $\phi_{n}$ are equipped with the $\cC$-equivariant structures $\{h^{c}_{n}\}$ and the $\psi_{n}$ are equipped with the $\cC$-equivariant structures $\{k^{c}_{n}\}$.
Also, the $i^{A}_{n}$ (resp. $i^{B}_{n}$) are equipped with $j^{F}_{n}$ (resp. $j^{G}_{n}$), and the structures commute precisely; namely:
\[
\begin{minipage}[]{0.3\textwidth}
\begin{align*}
\psi_{n}\circ \phi_{n} &= i^{A}_{n},\\
\phi_{n+1}  \circ \psi_{n} &= i^{B}_{n}, 
\end{align*}
\end{minipage}
\qquad\text{ and }\qquad
\begin{minipage}[]{0.3\textwidth}
\begin{align*}
k^{c}_{n}\circ h^{c}_{n} &= j^{F,c}_{n},\\
h^{c}_{n+1}\circ k^{c}_{n} &= j^{G,c}_{n}. 
\end{align*}
\end{minipage}
\]

Thus $\varinjlim \phi_{n}=:\phi$ and $\varinjlim \psi_{n}=:\psi$ extend to mutually inverse $\rm C^*$-algebra homomorphisms between $A$ and $B$. To equip $\phi$ with an equivariant structure, we will define
$h^{c}:=\varinjlim h^{c}_{n}: F(c)\rightarrow G(c)$. Technically we haven't made sense of the middle limit yet, since the $h^{c}_{n}$ maps are not morphisms in $\IL$. However, they are isometries by Lemma \ref{lem:StrucEqStrucs}.

Note $\lambda^{G(c)}_{n}\circ h^{c}_{n}: F_{n}(c)\rightarrow G(c)$ satisfies Lemma \ref{lem:StrucEqStrucs} (3), and satisfies $\lambda^{G(c)}_{n}\circ h^{c}_{n}=\lambda^{G(c)}_{n+1}\circ h^{c}_{n+1}\circ j^{F,c}_{n}$. Thus it extends to a bounded map $h^{c}: F(c)\rightarrow G(c)$ also satisfying Lemma \ref{lem:StrucEqStrucs} (3). It is clear that conditions $(1), (2), (4),$ and $(5)$, which occur at each finite $n$, will pass to the limit.
\end{proof}

Our result in particular restricts to a complete classification of AF-actions of finite groups on (unital) AF-algebras. From this perspective, our results have some overlap with previous results in the literature \cite{MR1193924, MR769758, MR1078515, MR1414065, https://doi.org/10.48550/arxiv.1304.0813, MR4206896}. A direct comparison with \cite{MR1414065, https://doi.org/10.48550/arxiv.1304.0813}, which deals with the case of $\bbZ/p\bbZ$ actions, is given in Section \ref{subsec:HilbZ/pZaction}. The other results have a different scope: on the one hand, they consider compact group actions rather than just finite groups, but on the other hand restrict the type of action to ``matrix summand preserving" actions (see \cite{MR4206896} for a summary). It is easy to see how the invariants utilized in these results recover pieces of our invariant in the case of finite group actions.

\subsection{Classification of strongly AF-inclusions}
\mbox{}

Let $A\subseteq B$ be an inclusion of unital $\rm C^*$-algebras with trivial center, such that $A^{\prime} \cap B=\bbC$, equipped with a faithful conditional expectation $E_A:B\rightarrow A$.
This allows us to equip $B$ with the structure of an $A$-$A$ correspondence, with right $A$-valued inner product $\langle b_1 |\  b_{2}\rangle_{A}:=E_A(b^{*}_{1}b_{2})$, left and right $A$-actions given by multiplication.
We say $(A\subseteq B, E_A)$ has \textit{finite depth} if the $\rm C^*$-tensor category tensor generated by ${}_{A}B_{A}$ is a fusion category. In particular, this implies ${}_{A}B_{A}\in \Bim(A)$. 

Now we have the following definition for an equivalence of inclusions of $\rm C^*$-algebras

\begin{defn}
An inclusion $(A\subseteq B, E_A)$ is \textit{equivalent} to the inclusion $(C\subseteq D, E_C)$ if there exists a $\rm C^*$-algebra isomorphism $\phi: B\rightarrow D$ such that $\phi(A)=C$ and $E_C\circ\phi=\phi\circ E_A$.  
\end{defn}

Popa's groundbreaking result in subfactor theory says that equivalence classes of finite depth inclusions of hyperfinite $\rm{II}_{1}$ factors are completely determined by their standard invariant \cite{MR1055708}, which in tensor categorical language is characterized by the fusion category of bimodules generated by ${}_{A}B_{A}\in \Bim(A)$ \cite{MR1424954} together with the Q-system ${}_{A}B_{A}$ \cite{MR1966524,MR3948170, MR3308880}. 

This motivates the question of classifying inclusions of AF $\rm C^*$-algebras. 
Obviously the standard invariant (which makes sense in the $\rm C^*$-setting just as well as the $\rm W^*$-setting \cite{MR4419534}) is a good starting place, however clearly more data will be required. 
A general classification, even for finite depth inclusions, seems out of reach at the moment. 
However, if we put additional finite dimensional approximation assumptions on the inclusion, then we can apply our main result Theorem \ref{thmalp:main} to obtain a classification. We introduce the following definitions:   

\begin{defn} If $A$ is an AF-algebra, a unitary tensor subcategory of $\cC\le \Bim(A)$ is called \textit{strongly AF} if there exists an increasing sequences of finite dimensional subalgebras $\{A_{n}\}\subseteq A$ whose union is dense in $A$ and a unitary tensor functor $F:\cC\rightarrow \IL(\{A_{n}\}\to \{A_{n}\})$ such that $\varinjlim \circ F $ is unitarily monoidally naturally isomorphic to $\id_{\cC}$.
\end{defn}

The strongly AF condition essentially means all the bimodules in the tensor category $\cC$ can be approximately by finite dimensional bimodules \textit{simultaneously}. This is a fairly strong condition. Note that it also implies that $\id_{\cC}\cong \varinjlim \circ F: \cC\to \Bim(A)$ is an AF-action.

\begin{defn}
An inclusion $(A\subseteq B,E_A)$ of AF-algebras is called a \textit{strongly AF-inclusion} if the unitary tensor category generated ${}_{A}B_{A}$ is strongly AF.
\end{defn}

\begin{defn}
Let $(A\subseteq B, E_A)$ be a finite depth inclusion of unital $\rm C^*$-algebras with trivial center. Then we define the \textit{extended standard invariant} 
$$\ExStd(A\subseteq B, E_A):=\left(\cC, Q, \widehat{F}, [A]\right),$$
where
\begin{itemize}
    \item 
    $\cC=\langle {}_{A} B_{A}\rangle$ is the fusion subcategory of $\Bim(A)$ tensor generated by ${}_{A} B_{A},$
    \item
    $Q= {}_{A} B_{A}\in \QSys(\cC)$ is the Q-system in $\cC$ defined by $(A\subseteq B, E_A)$ (\cite{MR4419534}).
    \item
    The pair $\left(\widehat{F}, [A]\right)$ is the invariant for the action of $\cC\overset{F}{\curvearrowright} A$ induced from the defining inclusion $F:\cC\hookrightarrow \Bim(A)$.
\end{itemize}
\end{defn}

Note that we have simply taken the usual standard invariant $(\cC, Q)$ and appended our invariant $\left(\widehat{F}, [A]\right)$ for the action $\cC\overset{F}{\curvearrowright} A$. There is a natural notion of equivalence for two such structures. Namely, an isomorphism
$\left(\cC, Q, \widehat{F}, [A]\right)\cong \left(\cD, R, \widehat{G}, [C]\right)$ consists of a triple $(H, \omega, \alpha),$ where
\begin{itemize}
    \item $H:\cC\cong \cD$ is a unitary equivalence of fusion categories.
    \item $\omega\in \cD(H(Q)\to R)$ is an isomorphism of Q-systems. 
    \item $\alpha:\widehat{F} \cong \widehat{G\circ H}$ is a natural isomorphism, such that $\alpha_{1_\cC}([A])=[C]$. 
\end{itemize}

\begin{thm*}[Theorem \ref{thmalp:ExStd&AFaction}]
Let $(A\subseteq B, E_A)$ and $(C\subseteq D, E_C)$ be strongly AF-inclusions of AF $\rm C^*$-algebras. Then $(A\subseteq B, E_A)$ is equivalent to $(C\subseteq D, E_C)$ if and only if $\ExStd(A\subseteq B, E_A)\cong \ExStd(C\subseteq D, E_C)$.
\end{thm*}
\begin{proof}
First suppose $(A\subseteq B, E_A)$ and $(C\subseteq D, E_C)$ are equivalent via in isomorphism $\phi: D\to B$, which restricts to an isomorphism $C\cong A$.
Conjugating by the invertible $C$-$A$ correspondence ${}_{\phi} A$ gives an equivalence $H:\Bim(A)\cong \Bim(C)$, 
which explicitly sends $X\in \Bim(A)$ to $X^{\phi}:={}_\phi A\boxtimes_A X\boxtimes_A {}_{\phi^{-1}}C\in \Bim(C)$ (as defined in the proof of Lemma \ref{lem:AfunderStrong}). 
Now, consider the Q-system ${}_{A}B_{A}\in \Bim(A)$. Then since $\phi:B\rightarrow D$ is an isomorphism of $\rm C^*$-algebras, we see that $\phi: D\to B^{\phi}$ is an isomorphism of Q-systems in $\Bim(C)$.
This data assembles into an equivalence of extended standard invariants.
Indeed, $H:=\widetilde{H}|_{\langle {}_{A}B_{A}\rangle }: \langle {}_{A}B_{A}\rangle \cong \langle {}_{C}D_{C}\rangle $ yields a unitary monoidal equivalence.
The map $\omega:=\phi^{-1}$ gives an isomorphism of Q-systems $H({}_{A}B_{A})$ and ${}_{C}D_{C}$. Finally, since $H$ is implemented by conjugation by ${}_{\phi} A$, this yields an equivalence of the identity actions $F:\langle {}_{A}B_{A}\rangle\rightarrow \Bim(A)$ and $G\circ H:\langle {}_{A}B_{A}\rangle\rightarrow \Bim(C)$, which gives an equivalence of invariants $\left(\widehat{F},[A]\right)\cong \left(\widehat{G\circ H}, [C]\right)$.

Conversely, suppose $(H,\omega,\alpha):\left(\cC,Q,\widehat{F},[A]\right)\cong \left(\cD,R,\widehat{G},[C]\right)$. 
According to the main Theorem \ref{thmalp:main}, the AF-actions $\cC\overset{F}{\curvearrowright} A$ and $\cC\ \overset{G\circ H}{\curvearrowright} C$ are equivalent. 
This gives a $*$-isomorphism $\varphi:A\to C$ and a family of linear maps \{$h^c:F(c)\to G(H(c))\}_{c\in\cC}$ satisfying the conditions of Lemma \ref{lem:StrucEqStrucs}.
Recall from Section \ref{subsec:realization}), $A\cong F(1_\cC)$, $B\cong F(Q)$, $C\cong G(1_\cD)$ and $D\cong G(R)$, with conditional expectations $E_A$ and $E_C$ given by $F(i^*_Q)$ and $G(i^*_R)$.

Now we show the following diagram commutes:
\[
\begin{tikzcd}
F(Q)
\arrow[swap]{d}{F(i^*_Q)}
\arrow{r}{h^Q}
& G(H(Q))
\arrow[swap]{d}{G(H(i^*_Q))}
\arrow{r}{G(\omega)}
& G(R)
\arrow{d}{G(i^*_R)}
 \\
F(1_\cC)
\arrow{r}{h^{1_\cC}}
& G(H(1_\cC))
\arrow[equal]{r}{G(H^1)}
& G(1_\cD)
\end{tikzcd}
\]
According to Lemma \ref{lem:StrucEqStrucs}(5) and (2), $h^Q:F(Q)\to G(H(Q))$ and $h^{1_{\cC}}:F(1_\cC)\to G(H(1_\cC))$ are isomorphisms of $\rm C^*$-algebras such that the left square commutes.
Since $\omega\in \cD(H(Q)\to R)$ is an isomorphism of Q-systems, it implies the following square of Q-system isomorphisms commutes:
\[
\begin{tikzcd}
H(Q)
\arrow[swap]{d}{H(i^*_Q)}
\arrow{r}{\omega}
& R
\arrow{d}{i^*_R}
 \\
H(1_\cC)
\arrow[equal]{r}{H^1}
& 1_\cD
\end{tikzcd}.
\]
After post-composing $G$, we have the right square of $\rm C^*$-algebra isomorphisms commutes.
Thus, $(A\subseteq B, E_A)$ is equivalent to $(C\subseteq D, E_C)$ as constructed.
\end{proof}

\section{Examples}\label{sec:Examples}

In this section, we will compute our invariant in several examples. $[\QSys(\cC)]$ itself has infinitely many objects, but it has only finitely many isomorphism classes of indecomposables. (See discussion above Remark \ref{rem:unitaryostrik}.) 
Furthermore, every object is a direct sum of these indecomposables. Since our functors $\widehat{F}$ are additive, this means the functor is determined by it's values on the subcategory of indecomposables. This means we can compare invariants restricted to categories with finitely many objects. We introduce some notation to make this precise.

\begin{defn}\label{def:chooseE} 
Let $\cE$ denote a set consisting of a representative of each isomorphism class of object in $[\QSys(\cC)]$, such that $1_{\cC}$ is chosen to represent $[1_{\cC}]$. 
We define $[\QSys(\cC)]_{\cE}$ to be the full subcategory of $[\QSys(\cC)],$ whose objects consist of all isomorphism classes of \emph{indecomposable} Q-systems in $\cC$ modulo Morita equivalence. 
\end{defn}

Clearly $[\QSys(\cC)]_{\cE}$ does not depend on the choice of set $\cE$ up to natural equivalence. 
The advantage of $[\QSys(\cC)]_{\cE}$ is that there are only finitely many objects, and between each pair objects, we have a canonical finite $\bbN$-basis consisting of isomorphism classes of irreducible bimodules. 
The structure of composition is then uniquely determined by a fusion rule between these canonical basis elements, which is a finite list of non-negative integers. 
Indeed, if $\cC$ is torsion-free in the sense of \cite{MR3941472}, then $[\QSys(\cC)]_{\cE}$ has a simple description: there is one object, the canonical basis elements correspond to (representatives of isomorphism classes of) simple objects in $\cC$, and the fusion rule is the usual one, see Section \ref{sec:TorsionFree}. 
For explicit examples in the more general case see Example \ref{ex:PointedCatas}, and Appendix \ref{App:Z4}.

\begin{lem}\label{lem:ReductionIndecomp}
Let $\cD$ be an abelian monoid-enriched category with direct sums. Then for any $\cE$ as in Definition \ref{def:chooseE}, the restriction functor gives an equivalence
$$\Fun_{+}([\QSys(\cC)]\to \cD)\cong \Fun_{+}([\QSys(\cC)]_{\cE}\to \cD).$$
\end{lem}
\begin{proof}
Since $[\QSys(\cC)]$ has direct sums and every object is isomorphic to a direct sum of indecomposables, the additive (i.e. direct sum) completion $[\QSys(\cC)]^{\oplus}_{\cE}$ is equivalent to $[\QSys(\cC)]$. In particular, the natural inclusion functor $i: [\QSys(\cC)]_{\cE} \to [\QSys(\cC)] $ extends to an equivalence $i^{\oplus}: [\QSys(\cC)]^{\oplus}_{\cE} \to [\QSys(\cC)].$

Now to show the restriction functor $\Fun_{+}([\QSys(\cC)]\to \cD)\to \Fun_{+}([\QSys(\cC)]_{\cE}\to \cD)$ is an equivalence, we find an inverse. This is given by the composition
$$\Fun_{+}([\QSys(\cC)]_{\cE}\to \cD)\cong \Fun_{+}([\QSys(\cC)]^{\oplus}_{\cE}\to \cD)\cong \Fun_{+}([\QSys(\cC)]\to \cD),$$
where the first natural isomorphism arises form the universal property of additive completion, and the second is from pre-composing with $(i^{\oplus})^{-1}$.
\end{proof}

We now have the following version of our main result, that allows us to compare invariants on the finite subcategory $[\QSys(\cC)]_{\cE}$.

\begin{cor}\label{cor:Fhatidecompobj}
Let $\cE$ be as in Definition \ref{def:chooseE}. 
Suppose $\cC\overset{F}{\curvearrowright} A$ and $\cC\overset{G}{\curvearrowright} B$ are AF-actions. 
Then $\left(\widehat{F}|_{[\QSys(\cC)]_{\cE}}, [A]\right)\cong \left(\widehat{G}|_{[\QSys(\cC)]_{\cE}}, [B]\right)$ if and only if $\cC\overset{F}{\curvearrowright} A$ is equivalent to $\cC\overset{G}{\curvearrowright} B$. 
\end{cor}

\subsection{Torsion-free fusion categories}\label{sec:TorsionFree}
\ \\
Following \cite{MR3941472}, we say a fusion category $\cC$ is \emph{torsion-free} if its only indecomposable $\cC$-module category is $\cC_\cC$.
Examples include $\Fib$ and $\Hilb(\bbZ/p\bbZ, \omega)$ where $p$ is prime and $[\omega]\in H^{3}(\bbZ/p\bbZ, \text{U}(1))$ is non-trivial.
In this case, $[\QSys(\cC)]$ has a single isomorphism class of indecomposable object represented by the trivial Q-system $1_{\cC}$. The endomorphism semiring of this object is just the fusion (semi)ring of $\cC$. Given an action $F:\cC\rightarrow \Bim(A)$, we have $\widehat{F}(1_{\cC})=(K_{0}(A), K^{+}_{0}(A))$, and fusion category objects $c\in\cC$ act on this by $-\boxtimes_{A}F(c)$. In particular, the functor $\widehat{F}$ simply reduces to the (pre)-ordered abelian group $K_{0}(A)$ equipped by a module structure over the fusion semiring of $\cC$. 

Our classification result has some immediate consequences for torsion-free fusion categories. We record the following theorem:
\begin{thm}\label{tmh:UniqeTorAction}
Let $A$ be an AF $\rm C^*$-algebra such that $(K_{0}(A), K^{+}_{0}(A))\subseteq (\bbR,\bbR^+)$ as a (scaled) ordered abelian group. If $\cC$ is a torsion-free fusion category, there is at most one AF-action on $A$.
\end{thm}
\begin{proof}
Suppose $\cC$ acts on $A$ via $F$. By \cite[Prop. 5.2]{MR4419534}, there is a state $\phi$ on $(K_{0}(A), K^{+}_{0}(A))$ such that $\phi([X\boxtimes_{A} F(a)])=\FPdim(a)\phi([X])$. 
Since $(K_{0}(A), K^{+}_{0}(A))\subseteq (\bbR,\bbR^+)$, there is a unique state, which corresponds to the inclusion $(K_{0}(A), K^{+}_{0}(A))\subseteq (\bbR,\bbR^+)$. Thus the fusion ring module structure is determined by $[- \boxtimes_{A} F(a)]=\FPdim(a)[- ]$. Thus there is a unique module structure. By Corollary \ref{cor:Fhatidecompobj}, this is a complete invariant for AF-actions. 
\end{proof}

Let $\cC=\Fib,$ the \emph{unitary Fibonacci category}, whose simple objects are $1,\tau,$ and $\tau\otimes\tau \cong 1\oplus \tau.$ Consider the AF-algebra $A_{\varphi}$ whose Bratteli diagram is given by 
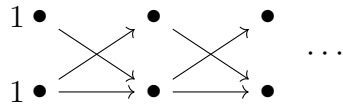
\begin{align*}
    \tikzmath{
\node (W11) at (0,0) {$\bullet$};
\node (W12) at (0,1) {$\bullet$};
\node (W21) at (1.5,0) {$\bullet$};
\node (W22) at (1.5,1) {$\bullet$};
\node (W31) at (3,0) {$\bullet$};
\node (W32) at (3,1) {$\bullet$};
\draw[->] (W11) -- (W21);
\draw[->] (W11) -- (W22);
\draw[->] (W12) -- (W21);
\draw[->] (W21) -- (W31);
\draw[->] (W21) -- (W32);
\draw[->] (W22) -- (W31);
\node at (3.75,.5) {$\cdots$};
\node at (-.3,0) {$1$};
\node at (-.3,1) {$1$};
}
\end{align*}
Recall that $K_0(A_\varphi)=\bbZ+\varphi\bbZ\subset \bbR,$ where $\varphi = \dfrac{1+\sqrt{5}}{2}$ is the golden ratio. 

\begin{cor}\label{cor:FibFib}
There is a unique AF-action of $\Fib$ on $A_{\varphi}$.
\end{cor}

\begin{proof}
By Theorem \ref{tmh:UniqeTorAction}, since $\Fib$ is torsion-free (e.g. \cite{MR3941472}) there is at most one AF-action on $A_{\varphi}$. Thus it suffices to show there exists one.

$[\QSys(\Fib)]$ has one isomorphism class of indecomposable object, which is represented by the trivial $\Fib$-module $\Fib$, and $[\QSys(\Fib)](\Fib\to \Fib)\cong \bbN[L^1]\oplus \bbN[L^\tau]$, where $L^c=c\otimes -$ is a $\Fib$-module functor for any $c\in\Fib$.
Consider the (homogeneous) enriched Bratteli diagram
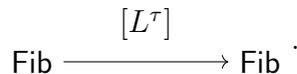
\[
\tikzmath{
\node (V1) at (-1.5,0) {$\Fib$};
\node (V2) at (1.5,0) {$\Fib$};
\draw[->] (V1) -- (V2);
\node at (0,.5) {$[L^\tau]$};
}.
\]
By applying the fusion rules we directly see that the underlying ordinary Bratteli diagram is precisely what we want:
\[
\tikzmath{
\node (W00) at (0,0) {$[\tau]$};
\node (W01) at (0,1) {$[1_{\Fib}]$};
\node (W10) at (2.5,0) {$[\tau]$};
\node (W11) at (2.5,1) {$[1_{\Fib}]$};
\draw[->] (W00) -- (W10);
\draw[->] (W00) -- (W11);
\draw[->] (W01) -- (W10);
\node at (3.25,.5) {$\cdots$};
}
\qquad
\Longrightarrow
\qquad
\tikzmath{
\node (W11) at (0,0) {$\bullet$};
\node (W12) at (0,1) {$\bullet$};
\node (W21) at (1.5,0) {$\bullet$};
\node (W22) at (1.5,1) {$\bullet$};
\draw[->] (W11) -- (W21);
\draw[->] (W11) -- (W22);
\draw[->] (W12) -- (W21);
\node at (2.25,.5) {$\cdots$};
\node at (-.3,0) {$1$};
\node at (-.3,1) {$1$};
}.
\]
\end{proof}

\subsection{$\Hilb(\bbZ/p\bbZ)$-actions}\label{subsec:HilbZ/pZaction}
\mbox{}

In this section, we explain our invariant for $\Hilb(\bbZ/p\bbZ)$ actions, and the close relation to the invariants of \cite{https://doi.org/10.48550/arxiv.1304.0813,MR1414065}.

Let $p$ be a prime. Then $[\QSys(\Hilb(\bbZ/p\bbZ))]$ has two isomorphism classes of indecomposable objects, which can be represented by the trivial Q-system $Q_{1}:=1_{\Hilb(\bbZ/p\bbZ)},$ and the group algebra $Q_{2}:=\bbC[\bbZ/p\bbZ]$. By Proposition \ref{prop:idcpHKbimod}, the indecomposable bimodules between them are given by
\begin{align*}
[\QSys(\Hilb(\bbZ/p\bbZ))](Q_{1}\to Q_{1}) &\cong \bbN[\bbZ/p\bbZ] \ni M_{1-1,g},   \\
[\QSys(\Hilb(\bbZ/p\bbZ))](Q_{2}\to Q_{2}) &\cong \bbN[\widehat{\bbZ/p\bbZ}] \ni M_{2-2,\pi},\\
[\QSys(\Hilb(\bbZ/p\bbZ))](Q_{1}\to Q_{2}) &\cong \bbN[ M_{1-2}], \\
[\QSys(\Hilb(\bbZ/p\bbZ))](Q_{2}\to Q_{1}) &\cong \bbN[ M_{2-1}].
\end{align*}

Here $M_{1-2}: = \bigoplus_{g\in \bbZ/p\bbZ} \bbC_g$ is the only indecomposable $Q_1-Q_2$ bimodule (up to isomorphism) and $M_{2-1}: = \bigoplus_{g\in \bbZ/p\bbZ} \bbC_g$ is the only indecomposable $Q_2-Q_1$ bimodule. The fusion rules are given by
\begin{align*}
    M_{1-1,g}\boxtimes_{Q_1} M_{1-1,h} &= M_{1-1,gh}, \\
    M_{2-2,\pi}\boxtimes_{Q_2} M_{2-2,\delta} &= M_{2-2,\pi\otimes \delta}, \\
    M_{2-1}\boxtimes_{Q_{1}} M_{1-1,g}=M_{2-1} &= M_{2-2,\pi}\boxtimes_{Q_2}M_{2-1}, \\
    M_{2-1}\boxtimes_{Q_1} M_{1-2} &= \bigoplus_{\pi\in \widehat{\bbZ/p\bbZ}} M_{2-2,\pi}, \\
    M_{1-2}\boxtimes_{Q_2} M_{2-1} &= \bigoplus_{g\in \bbZ/p\bbZ} M_{1-1,g}, \\
    M_{1-2}\boxtimes_{Q_2} M_{2-2,\pi} &= M_{1-2}.
\end{align*}

Let $G\overset{F}{\curvearrowright} A$ be an action of a finite group $G$ on the unital $\rm C^*$-algebra $A$. Then as described in Section \ref{subsec:realization}, this gives an action $\Hilb(G)\overset{F}{\curvearrowright} A$. 
The invariant in the case that $G=\bbZ/p\bbZ$ consists of the following information:
\begin{itemize}
    \item 
    The (pre)-ordered abelian group $\left(\widehat{F}(1_{\Hilb(\bbZ/p\bbZ)}), [A]\right)=(K_{0}(A), K^{+}_{0}(A), [A])$ with order unit $[A]$, equipped with the structure of a $\bbZ/p\bbZ$-module.
    \item
     The (pre)-ordered abelian group $\widehat{F}(\bbC[\bbZ/p\bbZ])=(K_{0}(A\rtimes \bbZ/p\bbZ), K^{+}_{0}(A\rtimes \bbZ/p\bbZ))$, equipped with the structure of a $\widehat{\bbZ/p\bbZ}$-module.
     \item
     $\widehat{F}(M_{1-2}):\widehat{F}(1_{\Hilb(\bbZ/p\bbZ)})\rightarrow \widehat{F}(\bbC[\bbZ/p\bbZ])$, which corresponds to the homomorphism $K_{0}(A)\hookrightarrow K_{0}(A\rtimes \bbZ/p\bbZ)$ induced by the natural inclusion.
     \item
     $\widehat{F}(M_{2-1}):\widehat{F}(\bbC[\bbZ/p\bbZ])\rightarrow \widehat{F}(1_{\Hilb(\bbZ/p\bbZ)})$, which corresponds to the homomorphism $K_{0}(A\rtimes \bbZ/p\bbZ)\rightarrow K_{0}(A)$ induced from the standard conditional expectation $E_{A}:A\rtimes \bbZ/p\bbZ\rightarrow A$.
\end{itemize}

From this description, it is clear that our invariant very closely resembles the invariant utilized in the classification results of \cite{MR1414065,https://doi.org/10.48550/arxiv.1304.0813} with two minor differences: firstly, we do not require specifying a distinguished order unit in $K_{0}(A\rtimes \bbZ/p\bbZ)$; and secondly, we require specification of the homomorphism $K_{0}(A\rtimes \bbZ/p\bbZ)\rightarrow K_{0}(A)$ induced by the conditional expectation.

\subsection{Computing the invariant in practice}\label{Ex:HilbZ4onZ1/4}
\mbox{}

In this section, we will explain how to use enriched Bratteli diagrams (Definition \ref{defn:EBD}) to calculate the invariant $\widehat{F}$. First note that an (unenriched) Bratteli diagram can be thought of as a Bratteli diagram enriched over $\Hilb$, i.e. a sequence of finitely semisimple categories
$$(\cN_0,y_0)\xrightarrow{J_0} (\cN_1,y_1)\xrightarrow{J_1} (\cN_2,y_2)\rightarrow\cdots$$
where $\cN_n$ is a finitely semisimple $\rm C^*$-category, $J_n:\cN_n\to \cN_{n+1}$ is a $*$-functor and $y_n\in\cN_n$ is a generating object such that $J_n(y_n)\cong y_{n+1}$.

This is not the usual description of a (connected) Bratteli diagram so we will briefly record the translation:
\begin{itemize}
    \item 
    For each $n\in \bbZ_{\ge 0}$ set $V_{n}$ to be a set of representatives of isomorphism classes of simple objects in $\cN_{n},$ 
    \item
    There are $k$ edges with source $v\in V_{n}$ and target $w\in V_{n+1}$, $k=\dim(\cN_{n+1}(w\to J_n(v))),$
    \item
    For each $n$, define the weight function $l:V_{n}\rightarrow \bbN$ by $l(v)=\dim(\cN_{n}(v\to y_{n}))$.
\end{itemize}
The above recovers the usual definition of a (connected) Bratteli diagram.
A Bratteli diagram without a weight function $l$ is called an \textit{unweighted} Bratteli diagram, which we denote as $(\cN_n,J_n)_n$. 
From this data, one can still compute the ordered $K$-group $K_{0}((\cN_{n}, J_{n})_n)$, but without the weight function we lack a distinguished choice of order unit. 


Now given an enriched Bratteli diagram 
$$(\cM_0,x_0)\xrightarrow{I_0} (\cM_1,x_1)\xrightarrow{I_1} (\cM_2,x_2)\rightarrow\cdots$$ describing an AF-action $F$ we want to compute $\widehat{F}(P)$ for any Q-system $P$. Let ${}_{P}\cC$ denote the corresponding $\cC$-module category of left $P$-modules (see Remark \ref{rem:unitaryostrik}).
Then define the finitely semisimple $\rm C^*$-categories 
$$\cM^{P}_{n}:=\fssMod(\cC)({}_{P} \cC\to \cM_{n}),$$ 
and the functors 
$$I^{P}_{n}:=  I_n \circ - :\cM^{P}_{n}\rightarrow \cM^{P}_{n+1}.$$
This gives an unweighted Bratteli diagram $(\cM^P_n, I^P_n)_n$.
Furthermore, according to Corollary \ref{Cor:FhatLimitFnhat} and Theorem \ref{thm:EquivariantCrossedprod},
we have:
\begin{align}\label{eqn:FHatP}
\begin{split}
    \widehat{F}(P) & = \varinjlim \widehat{F}_n(P) \\
    & = \varinjlim K_0(|P|_{F_n})  \\
    & \cong \varinjlim \cG\circ [\fssMod(\cC)]({}_P \cC\to \cM_n) \\
    & \cong K_{0}\left((\cM^P_n, I^P_n)_n\right).
\end{split}
\end{align}

Now let $X$ be a $P-Q$ bimodule. We need to build a positive homomorphism $\widehat{F}([X]):\widehat{F}(P)\rightarrow \widehat{F}(Q)$. For each $n$, we have the functor $L^{X}_{n}:\cM^{P}_n\rightarrow \cM^{Q}_{n}$, defined for $Y\in \cM^{P}_{n}=\fssMod({}_{P} \cC\to \cM_{n})$ and ${}_{Q} Z\in {}_{Q}\cC$
such that
$$L^{X}_{n}(Y)( {}_{Q}Z):=Y(X\otimes_{Q} Z).$$
We obtain a sequence of squares (commuting up to isomorphism):
\begin{equation}\label{eqn:FHatBim}
\begin{tikzcd}[column sep=3.5em]
 \cM^{P}_{0}
 \arrow[swap]{d}{L^X_0}
 \arrow{r}{I_0^P}
 &  \cM^P_1
 \arrow[swap]{d}{L^X_1}
 \arrow{r}{I_1^P}
 &  \cM^P_2
 \arrow[swap]{d}{L^X_2}
 \arrow{r}{I_2^P}
 & \cdots
 \arrow{r}{}
 &\widehat{F}(P)
 \arrow{d}{\widehat{F}([X])}\\
 \cM^Q_0
 \arrow[swap]{r}{I_0^Q}
 &  \cM^Q_1
 \arrow[swap]{r}{I_1^Q}
 &  \cM^Q_2
 \arrow[swap]{r}{I_2^Q}
 & \cdots
 \arrow{r}{}
 & \widehat{F}(Q)
 \end{tikzcd}.
\end{equation}

The positive homomorphism $\widehat F([X])=\varinjlim [L^X_{n}]$ can then be explicitly computed from universal properties of inductive limits through this diagram. 

\hspace{.5cm}

We will demonstrate this procedure for $\cC=\Hilb(\bbZ/4\bbZ)$. Although our formalism for enriched Bratteli diagrams uses the category $[\Mod(\cC)]$ to decorate data, we prefer for computational reasons to use the category $[\QSys(\cC)]$, which is contravariantly equivalent by Ostrik's theorem (see Remark \ref{rem:unitaryostrik}).
We choose this example since $4$ is the first composite natural.
There are three indecomposable Q-systems in $\Hilb(\bbZ/4\bbZ)$, which we denote $Q_1, Q_2, Q_3$. 
These correspond to the indecomposable Q-systems built from subgroups $\{0\}$, $\{0,2\}\cong\bbZ/2\bbZ$, and $\bbZ/4\bbZ,$ respectively. 
The morphisms and fusion rules can be found in Appendix \ref{App:Z4}.
Let $\cD$ denote the full subcategory of $[\QSys(\Hilb(\bbZ/4\bbZ))]$ whose objects are $Q_{1}, Q_{2}$ and $Q_{3}$ as in Definition \ref{def:chooseE}.

At this stage, it might be helpful for the reader to revisit Example \ref{ex:PointedCatas} and Proposition \ref{prop:idcpHKbimod}, and review Appendix \ref{App:Z4} to familiarize with the notation.
We now provide an example of three equivalent AF-actions $\Hilb(\bbZ/4\bbZ)\overset{}{\curvearrowright} A$ on the same UHF $\rm C^*$-algebra $A\cong \bbM_{4^\infty}$ (with $K_0(A)\cong\bbZ\!\left[\frac{1}{4}\right]$) given by the different (homogeneous) enriched Bratteli diagrams:

\begin{itemize}
\item For each level, the vertex $V_n = Q_1$ and edge from $V_n$ to $V_{n+1}$ is $ \bigoplus\limits_{g\in\bbZ/4\bbZ}M_{1-1,g}$;
\item For each level, the vertex $W_n = Q_2$ and edge from $W_n$ to $W_{n+1}$ is $\bigoplus\limits_{\substack{\ell=\triv,\sign\\ k=1,2} }M_{2-2,k}^{\ell}$; 
\item For each level, the vertex $U_n = Q_3$ and edge from $U_n$ to $U_{n+1}$ is $ \bigoplus\limits_{k=0}^3M_{3-3}^{\chi_k}$;
\end{itemize}

From left to right, we can represent each connection in the diagram graphically as
\[
\tikzmath{
\node (V1) at (-1.5,0) {$Q_1$};
\node (V2) at (1.5,0) {$Q_1$};
\draw[->] (V1) -- (V2);
\node at (0,.5) {$\bigoplus\limits_{g\in\bbZ/4\bbZ}M_{1-1,g}$};
}
\qquad\qquad
\tikzmath{
\node (V1) at (-1.5,0) {$Q_2$};
\node (V2) at (1.5,0) {$Q_2$};
\draw[->] (V1) -- (V2);
\node at (0,.5) {$\bigoplus\limits_{\substack{\ell=\triv,\sign\\ k=1,2} }M_{2-2,k}^{\ell}$};
}
\qquad\qquad
\tikzmath{
\node (V1) at (-1.5,0) {$Q_3$};
\node (V2) at (1.5,0) {$Q_3$};
\draw[->] (V1) -- (V2);
\node at (0,.55) {$\bigoplus\limits_{k=0}^3M_{3-3}^{\chi_k}$};
}
\]

We denote $\Hilb(\bbZ/4\bbZ)\overset{F}{\curvearrowright} A$, $\Hilb(\bbZ/4\bbZ)\overset{G}{\curvearrowright} A$ and $\Hilb(\bbZ/4\bbZ)\overset{H}{\curvearrowright} A$ to be the AF-actions induced by above enriched Bratteli diagrams from left to right, respectively. 
Afterwards, we will introduce a a fourth action denoted by $\Hilb(\bbZ/4\bbZ)\overset{E}{\curvearrowright} A$, which will turn out to be trivial, and will show none of $F,G$ or $H$ are equivalent to $E.$
\medskip

Let us compute the invariant $\widehat{G}$ on all of $\cD$; first on the objects:
\begin{itemize}
\item For $Q=Q_1$: There is only one vertex $Q_2$ in each level, i.e., $\cM_n={}_{Q_2}\cC$, so by Equation (\ref{eqn:FHatP}),
$$\widehat{G}(Q_1) = \varinjlim \cG\circ [\Mod(\Hilb(\bbZ/4\bbZ))]({}_{Q_1}\cC\to {}_{Q_2}\cC)= \varinjlim \cG\circ  \cD(Q_2\to Q_1), $$
where the $n$-th level $\cG\circ \cD(Q_2\to Q_1)_n$ has $\bbZ$-basis $M_{2-1,0}$ and $M_{2-1,1}$. 
According to Equation (\ref{eqn:FHatP}), the connecting map $I_n^{Q_1}$ is implemented by $$\left(\bigoplus\limits_{\substack{\ell=\triv,\sign\\ k=1,2} } M_{2-2,k}^{\ell}\right) \boxtimes_{Q_2} -: \cD(Q_2\to Q_1)_n \to \cD(Q_2\to Q_1)_{n+1}.$$
We now compute the image of this map on the $\bbZ$-basis elements. 
Note that, by Appendix \ref{App:Z4}, $M_{2-2,k}^{\ell}\boxtimes_{Q_2} M_{2-1,j} = M_{2-1,k+j}$, $k,j=0,1$, $\ell=\triv,\sign,$ and the sum $k+j$ is meant modulo 2. 
By looking at the action of the connecting map of the decomposition of the $n$-th and $(n+1)$-th levels using the given $\bbZ$-basis, we can capture the behavior of the connecting maps using a Bratteli diagram, modulo remembering the labeling of each block:
\[
\tikzmath{
\node (W00) at (0,0) {$M_{2-1,0}$};
\node (W01) at (0,1) {$M_{2-1,1}$};
\node (W10) at (2,0) {$M_{2-1,0}$};
\node (W11) at (2,1) {$M_{2-1,1}$};
\draw[->] (W00.10) -- (W10.170);
\draw[->] (W00.350) -- (W10.190);
\draw[->] (W01.10) -- (W11.170);
\draw[->] (W01.350) -- (W11.190);
\draw[postaction={transform canvas={yshift=-.77mm,xshift=1.27mm},draw}]
[->] (W00.50) to (W11.220);
\draw[postaction={transform canvas={yshift=.77mm,xshift=1.27mm},draw}]
[->] (W01.310) to (W10.140);
\node at (3,.5) {$\cdots$};
}
\qquad \Longrightarrow \qquad
\tikzmath{
\node (W11) at (0,0) {$\bullet$};
\node (W12) at (0,1) {$\bullet$};
\node (W21) at (1.5,0) {$\bullet$};
\node (W22) at (1.5,1) {$\bullet$};
\draw[->] (W11.15) -- (W21.165);
\draw[->] (W11.345) -- (W21.195);
\draw[->] (W12.15) -- (W22.165);
\draw[->] (W12.345) -- (W22.195);
\draw[postaction={transform canvas={yshift=.77mm,xshift=-.77mm},draw}]
[->] (W11.35) to (W22.235);
\draw[postaction={transform canvas={yshift=.77mm,xshift=.77mm},draw}]
[->] (W12.305) to (W21.145);
\node at (2.25,.5) {$\cdots$};
}.
\] 

According to Equation (\ref{eqn:FHatP}) for $\widehat{G}(Q_1),$
and the well-known $K_0$ computation of this Bratteli diagram, we can define a group isomorphism 
\begin{align*}
    \phi_1 : \widehat{G}(Q_1) &\to \bbZ\!\left[\textstyle\frac{1}{4}\right] \\
    (x,y)&\mapsto 4^{-n}\cdot\dfrac{x+y}{2}.
\end{align*}
Here, we take $n$ as the first positive integer for which $(x,y)\in \cG\circ\cD(Q_2\to Q_1)_n = \bbZ M_{2-1,0} \oplus \bbZ M_{2-1,1}$ appears in the inductive system.
\item For $Q=Q_2$: 
$\widehat{G}(Q_2) = \varinjlim \cG\circ  \cD(Q_2\to Q_2)$, 
where $\cG\circ\cD(Q_2\to Q_2)_n$ has $\bbZ$-basis $\left\{ M_{2-2,k}^{\ell} \right\}_{k=0,1}^{\ell = \triv, \sign}$. 
Similarly as we did above, by decomposing the connecting map, which is given by tensoring by $\bigoplus\limits_{\substack{\ell=\triv,\sign\\ k=1,2} } M_{2-2,k}^{\ell}$ 
on the left, we obtain the following Bratteli diagram:
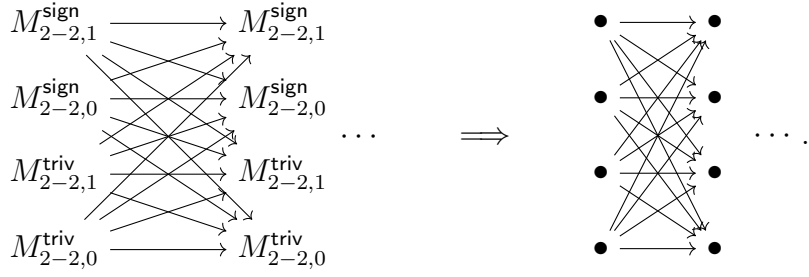
\[
\tikzmath{
\node (W00) at (0,0) {$M_{2-2,0}^\triv$};
\node (W01) at (0,1) {$M_{2-2,1}^\triv$};
\node (W02) at (0,2) {$M_{2-2,0}^\sign$};
\node (W03) at (0,3) {$M_{2-2,1}^\sign$};
\node (W10) at (3,0) {$M_{2-2,0}^\triv$};
\node (W11) at (3,1) {$M_{2-2,1}^\triv$};
\node (W12) at (3,2) {$M_{2-2,0}^\sign$};
\node (W13) at (3,3) {$M_{2-2,1}^\sign$};
\draw[->] (W00) -- (W10);
\draw[->] (W00) -- (W11);
\draw[->] (W00) -- (W12);
\draw[->] (W00) -- (W13);
\draw[->] (W01) -- (W10);
\draw[->] (W01) -- (W11);
\draw[->] (W01) -- (W12);
\draw[->] (W01) -- (W13);
\draw[->] (W02) -- (W10);
\draw[->] (W02) -- (W11);
\draw[->] (W02) -- (W12);
\draw[->] (W02) -- (W13);
\draw[->] (W03) -- (W10);
\draw[->] (W03) -- (W11);
\draw[->] (W03) -- (W12);
\draw[->] (W03) -- (W13);
\node at (4,1.5) {$\cdots$};
}
\qquad\Longrightarrow\qquad
\tikzmath{
\node (U11) at (0,0) {$\bullet$};
\node (U12) at (0,1) {$\bullet$};
\node (U13) at (0,2) {$\bullet$};
\node (U14) at (0,3) {$\bullet$};
\node (U21) at (1.5,0) {$\bullet$};
\node (U22) at (1.5,1) {$\bullet$};
\node (U23) at (1.5,2) {$\bullet$};
\node (U24) at (1.5,3) {$\bullet$};
\draw[->] (U11) -- (U21);
\draw[->] (U11) -- (U22);
\draw[->] (U11) -- (U23);
\draw[->] (U11) -- (U24);
\draw[->] (U12) -- (U21);
\draw[->] (U12) -- (U22);
\draw[->] (U12) -- (U23);
\draw[->] (U12) -- (U24);
\draw[->] (U13) -- (U21);
\draw[->] (U13) -- (U22);
\draw[->] (U13) -- (U23);
\draw[->] (U13) -- (U24);
\draw[->] (U14) -- (U21);
\draw[->] (U14) -- (U22);
\draw[->] (U14) -- (U23);
\draw[->] (U14) -- (U24);
\node at (2.25,1.5) {$\cdots$};
}.
\]
Similarly, we obtain a group isomorphism 
\begin{align*}
    \phi_2 : \widehat{G}(Q_2) &\to \bbZ\!\left[\textstyle\frac{1}{4}\right] \\
    (x,y,z,w)&\mapsto 4^{-n}\cdot\frac{x+y+z+w}{4}
\end{align*}
Here, $(x,y,z,w)\in \cG\circ\cD(Q_2\to Q_2)_n = \bbZ M_{2-2,0}^\triv\oplus\bbZ M_{2-2,1}^\triv\oplus \bbZ M_{2-2,0}^\sign\oplus \bbZ M_{2-2,1}^\sign$ and $n$ corresponds to the level at which these first appear.

\item For $Q=Q_3$:
$\widehat{G}(Q_3) = \varinjlim \cG\circ  \cD(Q_2\to Q_3)$, 
where $\cG\circ \cD(Q_2\to Q_3)_n$ has $\bbZ$-basis $M_{2-3}^\triv$ and $M_{2-3}^\sign$. 
Similarly as we did above, by decomposing the connecting map, which is the same as before, obtain the following Bratteli diagram:
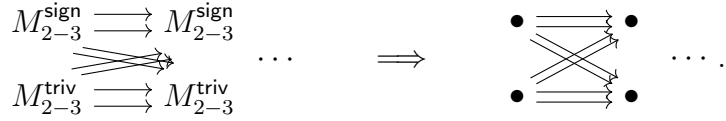
\[
\tikzmath{
\node (W00) at (0,0) {$M_{2-3}^\triv$};
\node (W01) at (0,1) {$M_{2-3}^\sign$};
\node (W10) at (2,0) {$M_{2-3}^\triv$};
\node (W11) at (2,1) {$M_{2-3}^\sign$};
\draw[->] (W00.10) -- (W10.170);
\draw[->] (W00.350) -- (W10.190);
\draw[->] (W01.10) -- (W11.170);
\draw[->] (W01.350) -- (W11.190);
\draw[postaction={transform canvas={yshift=-.77mm,xshift=1.27mm},draw}]
[->] (W00.50) to (W11.220);
\draw[postaction={transform canvas={yshift=.77mm,xshift=1.27mm},draw}]
[->] (W01.310) to (W10.140);
\node at (3,.5) {$\cdots$};
}
\qquad \Longrightarrow \qquad
\tikzmath{
\node (W11) at (0,0) {$\bullet$};
\node (W12) at (0,1) {$\bullet$};
\node (W21) at (1.5,0) {$\bullet$};
\node (W22) at (1.5,1) {$\bullet$};
\draw[->] (W11.15) -- (W21.165);
\draw[->] (W11.345) -- (W21.195);
\draw[->] (W12.15) -- (W22.165);
\draw[->] (W12.345) -- (W22.195);
\draw[postaction={transform canvas={yshift=.77mm,xshift=-.77mm},draw}]
[->] (W11.35) to (W22.235);
\draw[postaction={transform canvas={yshift=.77mm,xshift=.77mm},draw}]
[->] (W12.305) to (W21.145);
\node at (2.25,.5) {$\cdots$};
}.
\]
Once again we obtain a group isomorphism 
\begin{align*}
    \phi_3 : \widehat{G}(Q_3) &\to \bbZ\!\left[\textstyle\frac{1}{4}\right] \\
    (x,y)&\mapsto 4^{-n}\cdot \frac{x+y}{2}.
\end{align*}
Here, we take $n$ as the first positive integer for which $(x,y)\in \cG\circ\cD(Q_2\to Q_3)_n = \bbZ M_{2-3}^\triv \oplus \bbZ M_{2-3}^\sign$ appears in the inductive system.
\end{itemize}

Now we compute $\widehat{G}$ of morphisms in $\cD$.
Again, it suffices to compute $\widehat{G}$ for those indecomposable bimodules. 

\begin{itemize}
\item For $M_{1-2,0}$: 
Note that $\widehat{G}(M_{1-2,0})$ is a group homomorphism from $\widehat{G}(Q_1)$ to $\widehat{G}(Q_2)$, which is given by the inductive limits as shown above.
We are going to unpack the Inductive Sequence (\ref{eqn:FHatBim}) here. Since we are working with a homogeneous enriched Bratteli diagram, it suffices to discuss the one commuting square between the $n$-th level and the $(n+1)$-th level.

For $Y\in \cD(Q_2\to Q_1)_n\hookrightarrow \widehat{G}(Q_1)$, $$L^{M_{1-2,0}}_n (Y) = Y\boxtimes_{Q_1} M_{1-2,0}\in \cD(Q_2\to Q_2)_n \hookrightarrow\widehat{G}(Q_2).$$
When picking up the $\bbZ$-basis $M_{2-1,0}$, $M_{2-1,1}$ for $\cG\circ \cD(Q_2\to Q_1)_n$ and $\bbZ$-basis $M_{2-2,0}^\triv$, $M_{2-2,1}^\triv, M_{2-2,0}^\sign, M_{2-2,1}^\sign$ for $\cG\circ \cD(Q_2\to Q_2)_n$, we have
\begin{align*}
    L_n^{M_{1-2,0}}(M_{2-1,0}) &= M_{2-1,0}\boxtimes_{Q_2} M_{1-2,0} = M_{2-2,0}^\triv \oplus M_{2-2,0}^\sign, \\
    L_n^{M_{1-2,0}}(M_{2-1,1}) &= M_{2-1,0}\boxtimes_{Q_2} M_{1-2,0} = M_{2-2,1}^\triv \oplus M_{2-2,1}^\sign.
\end{align*}
The expressions for $L_{n+1}^{M_{1-2,0}}$ are the same as the above. 
In the diagram below, we depict the commuting diagram from the Inductive Sequence (\ref{eqn:FHatBim}) at the levels $n$ and $n+1,$ where the actions of $L_n^{M_{1-2,0}}$ and $L_{n+1}^{M_{1-2,0}}$are depicted by the thick red arrows:
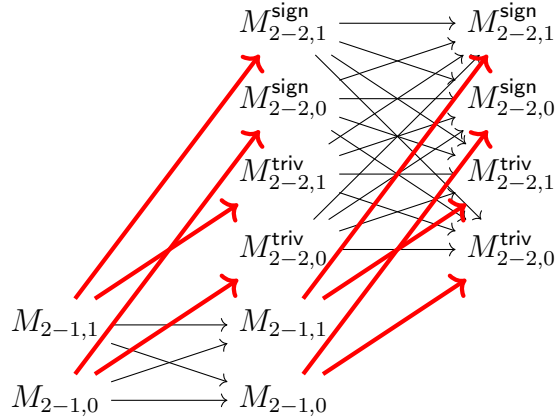
\[
\tikzmath{
\node (W00) at (0,0) {$M_{2-1,0}$};
\node (W01) at (0,1) {$M_{2-1,1}$};
\node (W10) at (3,0) {$M_{2-1,0}$};
\node (W11) at (3,1) {$M_{2-1,1}$};
\node (W12) at (3,2) {$M_{2-2,0}^\triv$};
\node (W13) at (3,3) {$M_{2-2,1}^\triv$};
\node (W14) at (3,4) {$M_{2-2,0}^\sign$};
\node (W15) at (3,5) {$M_{2-2,1}^\sign$};
\node (W22) at (6,2) {$M_{2-2,0}^\triv$};
\node (W23) at (6,3) {$M_{2-2,1}^\triv$};
\node (W24) at (6,4) {$M_{2-2,0}^\sign$};
\node (W25) at (6,5) {$M_{2-2,1}^\sign$};
\draw[->] (W00) -- (W10);
\draw[->] (W00) -- (W11);
\draw[->] (W01) -- (W10);
\draw[->] (W01) -- (W11);
\draw[->] (W12) -- (W22);
\draw[->] (W12) -- (W23);
\draw[->] (W12) -- (W24);
\draw[->] (W12) -- (W25);
\draw[->] (W13) -- (W22);
\draw[->] (W13) -- (W23);
\draw[->] (W13) -- (W24);
\draw[->] (W13) -- (W25);
\draw[->] (W14) -- (W22);
\draw[->] (W14) -- (W23);
\draw[->] (W14) -- (W24);
\draw[->] (W14) -- (W25);
\draw[->] (W15) -- (W22);
\draw[->] (W15) -- (W23);
\draw[->] (W15) -- (W24);
\draw[->] (W15) -- (W25);
\draw[->,\XColor,ultra thick] (W00) -- (W12);
\draw[->,\XColor,ultra thick] (W00) -- (W14);
\draw[->,\XColor,ultra thick] (W01) -- (W13);
\draw[->,\XColor,ultra thick] (W01) -- (W15);
\draw[->,\XColor,ultra thick] (W10) -- (W22);
\draw[->,\XColor,ultra thick] (W10) -- (W24);
\draw[->,\XColor,ultra thick] (W11) -- (W23);
\draw[->,\XColor,ultra thick] (W11) -- (W25);
}
\]
In order to compute $\widehat{G}(M_{1-2,0})$, we use the isomorphism $\phi_1:\widehat{G}(Q_1)\to \bbZ\!\left[\frac{1}{4}\right]$ and $\phi_2:\widehat{G}(Q_2)\to \bbZ\!\left[\frac{1}{4}\right]$ above in order to compute it on the given $\bbZ$-basis. 
Let $f:\bbZ\!\left[\frac{1}{4}\right]\to \bbZ\!\left[\frac{1}{4}\right]$ be a group homomorphism such that the following diagram commute: 
\[
\begin{tikzcd}[column sep=4em]
\widehat{G}(Q_1)
\arrow[swap]{d}{\phi_1}
\arrow{r}{\widehat{G}(M_{1-2,0})}
& \widehat{G}(Q_2)
\arrow{d}{\phi_2}
\\
\bbZ\!\left[\frac{1}{4}\right]
\arrow[]{r}{f}
& \bbZ\!\left[\frac{1}{4}\right]
\end{tikzcd}
\]
For $(x,y)\in \cG\circ \cD(Q_2\to Q_1)_n=\bbZ M_{2-1,0}\oplus \bbZ M_{2-1,1}$, 
by the previous diagram,
we have $\widehat{G}(M_{1-2,0})(x,y) = (x,y,x,y)\in \cG\circ \cD(Q_2\to Q_2)_n=\bbZ M_{2-2,0}^\triv \oplus\bbZ M_{2-2,1}^\triv\oplus\bbZ M_{2-2,0}^\sign\oplus\bbZ M_{2-2,1}^\sign$. 
Recall the isomorphisms $\phi_1$ and $\phi_2$ to $\bbZ\!\left[\frac{1}{4}\right]$,
we can see 
$$\phi_2\circ \widehat{G}(M_{1-2,0})(x,y) = \phi_2 (x,y,x,y) = 4^{-n}\cdot\frac{2x+2y}{4} = \phi_1(x,y) = f\circ \phi_1(x,y).$$

Therefore, $f$ is an identity map and $\widehat{G}(M_{1-2,0})= \phi_2^{-1}\circ f\circ \phi_1 = \phi_2^{-1}\circ \phi_1$. 
We say $\widehat{G}(M_{1-2,0})$ is induced by an identity map.
Similarly, $\widehat{G}(M_{1-2,1})$ is induced by an identity map.
\item For $M_{1-1,g}$ $M_{2-2,0/1}^{\triv/\sign}$ and $M_{3-3}^{\chi_k}$: Similarly, it is easy to see their $\widehat{G}$'s are induced by identity maps.
\item For $M_{2-1,0}$: 
When picking up the $\bbZ$-basis $M_{2-2,0}^\triv$, $M_{2-2,1}^\triv$, $M_{2-2,0}^\sign$, $M_{2-2,1}^\sign$ for $\cG\circ \cD(Q_2\to Q_2)_n$ and $\bbZ$-basis $M_{2-1,0}$, $M_{2-1,1}$ for $\cG\circ \cD(Q_2\to Q_1)_n$, we have
\begin{align*}
    L_n^{M_{2-1,0}}(M_{2-2,i}^\ell) = M_{2-2,i}^\ell\boxtimes_{Q_2} M_{2-1,0} = M_{2-1,i} \qquad i=0,1,\ \ell=\triv/\sign.
\end{align*}
We then obtain the following diagram:
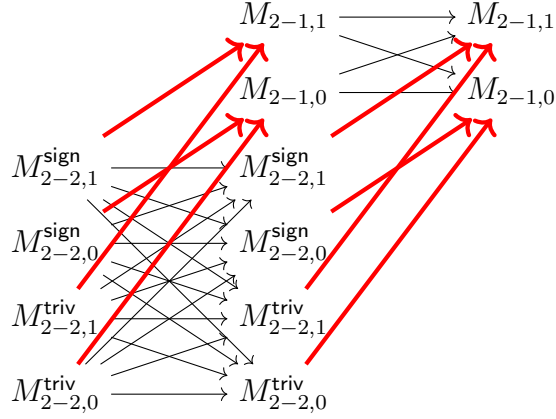
\[
\tikzmath{
\node (W00) at (0,0) {$M_{2-2,0}^\triv$};
\node (W01) at (0,1) {$M_{2-2,1}^\triv$};
\node (W02) at (0,2) {$M_{2-2,0}^\sign$};
\node (W03) at (0,3) {$M_{2-2,1}^\sign$};
\node (W10) at (3,0) {$M_{2-2,0}^\triv$};
\node (W11) at (3,1) {$M_{2-2,1}^\triv$};
\node (W12) at (3,2) {$M_{2-2,0}^\sign$};
\node (W13) at (3,3) {$M_{2-2,1}^\sign$};
\node (W14) at (3,4) {$M_{2-1,0}$};
\node (W15) at (3,5) {$M_{2-1,1}$};
\node (W24) at (6,4) {$M_{2-1,0}$};
\node (W25) at (6,5) {$M_{2-1,1}$};
\draw[->] (W00) -- (W10);
\draw[->] (W00) -- (W11);
\draw[->] (W00) -- (W12);
\draw[->] (W00) -- (W13);
\draw[->] (W01) -- (W10);
\draw[->] (W01) -- (W11);
\draw[->] (W01) -- (W12);
\draw[->] (W01) -- (W13);
\draw[->] (W02) -- (W10);
\draw[->] (W02) -- (W11);
\draw[->] (W02) -- (W12);
\draw[->] (W02) -- (W13);
\draw[->] (W03) -- (W10);
\draw[->] (W03) -- (W11);
\draw[->] (W03) -- (W12);
\draw[->] (W03) -- (W13);
\draw[->] (W14) -- (W24);
\draw[->] (W14) -- (W25);
\draw[->] (W15) -- (W24);
\draw[->] (W15) -- (W25);
\draw[->,\XColor,ultra thick] (W00) -- (W14);
\draw[->,\XColor,ultra thick] (W01) -- (W15);
\draw[->,\XColor,ultra thick] (W02) -- (W14);
\draw[->,\XColor,ultra thick] (W03) -- (W15);
\draw[->,\XColor,ultra thick] (W10) -- (W24);
\draw[->,\XColor,ultra thick] (W11) -- (W25);
\draw[->,\XColor,ultra thick] (W12) -- (W24);
\draw[->,\XColor,ultra thick] (W13) -- (W25);
}
\]

For $(x,y,z,w)\in \cG\circ \cD(Q_2\to Q_2)_n=\bbZ M_{2-2,0}^\triv \oplus\bbZ M_{2-2,1}^\triv\oplus\bbZ M_{2-2,0}^\sign\oplus\bbZ M_{2-2,1}^\sign$,
$\widehat{G}(x,y,z,w) = (x+z,y+w) \in \cG\circ \cD(Q_2\to Q_1)_n = \bbZ M_{2-1,0} \oplus\bbZ M_{2-1,1}$.
Recalling the isomorphisms $\phi_1$ and $\phi_2$ to $\bbZ\!\left[\frac{1}{4}\right]$,
we can see 
$$\phi_1\circ \widehat{G}(M_{2-1,0})(x,y,z,w) = \phi_1 (x+z,y+w) = 4^{-n}\cdot\frac{x+y+z+w}{2} = 2\phi_2(x,y,z,w).$$
Therefore, $\widehat{G}(M_{2-1,0})$ is induced by a map of multiplication by 2.
Similarly, $\widehat{G}(M_{2-1,1})$ is induced by a map of multiplication by 2.

To see the consistency of our computations so far, 
note that $M_{1-2,0}\boxtimes_{Q_2} M_{2-1,0} = M_{1-1,0}\oplus M_{1-1,2}$ and $M_{2-1,0}\boxtimes_{Q_1} M_{1-2,0} = M_{2-2,0}^\triv\oplus M_{2-2,0}^\sign$,
we can see 
\begin{align*}
    \widehat{G}(M_{1-2,0}\boxtimes_{Q_2} M_{2-1,0}) &= \widehat{G}(M_{1-1,0})+\widehat{G}(M_{1-1,2}) = \widehat{G}(M_{1-2,0})\circ \widehat{G}(M_{2-1,0}), \\
    \widehat{G}(M_{2-1,0}\boxtimes_{Q_1} M_{1-2,0}) &= \widehat{G}(M_{2-2,0}^\triv)+ \widehat{G}(M_{2-2,0}^\sign) = \widehat{G}(M_{2-1,0})\circ \widehat{G}(M_{1-2,0})
\end{align*}
as maps of multiplication by 2.
\item For $M_{1-3}:$
When picking up the $\bbZ$-basis $M_{2-1,0}$ and $M_{2-1,1}$ for $\cG\circ \cD(Q_2\to Q_1)_n$ and $\bbZ$-basis $M_{2-3}^\triv$ and $M_{2-3}^\sign$ for $\cG\circ \cD(Q_2\to Q_3)_n$,
we have
\begin{align*}
    L^{M_{1-3}}_n( M_{2-1,0})= M_{2-1,0}\boxtimes_{Q_1} M_{1-3}= M_{2-3}^\triv\oplus M_{2-3}^\sign = L^{M_{1-3}}_n( M_{2-1,1}).
\end{align*}
We then obtain the following diagram:
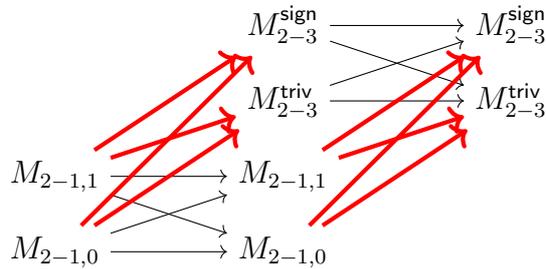
\[
\tikzmath{
\node (W00) at (0,0) {$M_{2-1,0}$};
\node (W01) at (0,1) {$M_{2-1,1}$};
\node (W10) at (3,0) {$M_{2-1,0}$};
\node (W11) at (3,1) {$M_{2-1,1}$};
\node (W12) at (3,2) {$M_{2-3}^\triv$};
\node (W13) at (3,3) {$M_{2-3}^\sign$};
\node (W22) at (6,2) {$M_{2-3}^\triv$};
\node (W23) at (6,3) {$M_{2-3}^\sign$};
\draw[->] (W00) -- (W10);
\draw[->] (W00) -- (W11);
\draw[->] (W01) -- (W10);
\draw[->] (W01) -- (W11);
\draw[->] (W12) -- (W22);
\draw[->] (W12) -- (W23);
\draw[->] (W13) -- (W22);
\draw[->] (W13) -- (W23);
\draw[->,\XColor,ultra thick] (W00) -- (W12);
\draw[->,\XColor,ultra thick] (W00) -- (W13);
\draw[->,\XColor,ultra thick] (W01) -- (W12);
\draw[->,\XColor,ultra thick] (W01) -- (W13);
\draw[->,\XColor,ultra thick] (W10) -- (W22);
\draw[->,\XColor,ultra thick] (W10) -- (W23);
\draw[->,\XColor,ultra thick] (W11) -- (W22);
\draw[->,\XColor,ultra thick] (W11) -- (W23);
}
\]
For $(x,y)\in \cG\circ \cD(Q_2\to Q_1)_n = \bbZ M_{2-1,0}\oplus\bbZ M_{2-1,1}$,
$\widehat{G}(M_{1-3})(x,y) = (2x,2y)\in \cG\circ \cD(Q_2\to Q_3)_n = \bbZ M_{2-3}^{\triv}\oplus \bbZ M_{2-3}^\sign$.
Recalling the isomorphism $\phi_1$ and $\phi_3$ to $\bbZ\!\left[\frac{1}{4}\right]$,
we can see 
$$\phi_3\circ \widehat{G}(M_{1-3})(x,y) = \phi_3 (2x,2y) = 4^{-n}\cdot\frac{2x+2y}{2} = 2\phi_1(x,y).$$
Therefore, $\widehat{G}(M_{1-3})$ is induced by a map of multiplication by 2.
\item For $M_{3-1}$: 
Note that $L^{M_{3-1}}_n (M_{2-3}^\triv) = M_{2-1,0}\oplus M_{2-1,1} = L^{M_{3-1}}_n ( M_{2-3}^\sign)$, 
similarly, 
it is easy to see $\widehat{G}(M_{3-1})$ is also induced by a map of multiplication by 2.

\hspace{.5cm}

Note that $M_{1-3}\boxtimes_{Q_3} M_{3-1} = \bigoplus_{g\in\bbZ/4\bbZ} M_{1-1,g}$ and $M_{3-1}\boxtimes_{Q_1} M_{1-3} = \bigoplus_{k=0}^3 M_{3-3}^{\chi_k}$,
we can see 
\begin{align*}
    \widehat{G}(M_{1-3}\boxtimes_{Q_3} M_{3-1}) &= \sum_{g\in\bbZ/4\bbZ} \widehat{G}(M_{1-1,g}) = \widehat{G}(M_{1-3})\circ \widehat{G}(M_{3-1}) \\
    \widehat{G}(M_{3-1}\boxtimes_{Q_1} M_{1-3}) &= \sum_{k=0}^3 \widehat{G}(M_{3-3}^{\chi_k}) = \widehat{G}(M_{3-1})\circ \widehat{G}(M_{1-3})
\end{align*}
as maps of multiplication by 4 as expected.
\item
For $M_{3-2}^\triv$: 
When picking up the $\bbZ$-basis $M_{2-3}^\triv$, $M_{2-3}^\sign$ for $\cG\circ \cD(Q_2\to Q_3)_n$ and $\bbZ$-basis $M_{2-2,0}^\triv$, $M_{2-2,1}^\triv$, $M_{2-2,0}^\sign$, $M_{2-2,1}^\sign$ for $\cG\circ \cD(Q_2\to Q_2)_n$, we have
\begin{align*}
    L^{M_{3-2}^\triv}_n (M_{2-3}^\triv) &= M_{2-3}^\triv\boxtimes_{Q_3} M_{3-2}^\triv = M_{2-2,0}^\triv\oplus M_{2-2,1}^\triv \\
    L^{M_{3-2}^\triv}_n (M_{2-3}^\sign) &= M_{2-3}^\sign\boxtimes_{Q_3} M_{3-2}^\triv = M_{2-2,0}^\sign\oplus M_{2-2,1}^\sign
\end{align*}
we obtain the following diagram:
\[
\tikzmath{
\node (W00) at (0,0) {$M_{2-3}^\triv$};
\node (W01) at (0,1) {$M_{2-3}^\sign$};
\node (W10) at (3,0) {$M_{2-3}^\triv$};
\node (W11) at (3,1) {$M_{2-3}^\sign$};
\node (W12) at (3,2) {$M_{2-2,0}^\triv$};
\node (W13) at (3,3) {$M_{2-2,1}^\triv$};
\node (W14) at (3,4) {$M_{2-2,0}^\sign$};
\node (W15) at (3,5) {$M_{2-2,1}^\sign$};
\node (W22) at (6,2) {$M_{2-2,0}^\triv$};
\node (W23) at (6,3) {$M_{2-2,1}^\triv$};
\node (W24) at (6,4) {$M_{2-2,0}^\sign$};
\node (W25) at (6,5) {$M_{2-2,1}^\sign$};
\draw[->] (W00) -- (W10);
\draw[->] (W00) -- (W11);
\draw[->] (W01) -- (W10);
\draw[->] (W01) -- (W11);
\draw[->] (W12) -- (W22);
\draw[->] (W12) -- (W23);
\draw[->] (W12) -- (W24);
\draw[->] (W12) -- (W25);
\draw[->] (W13) -- (W22);
\draw[->] (W13) -- (W23);
\draw[->] (W13) -- (W24);
\draw[->] (W13) -- (W25);
\draw[->] (W14) -- (W22);
\draw[->] (W14) -- (W23);
\draw[->] (W14) -- (W24);
\draw[->] (W14) -- (W25);
\draw[->] (W15) -- (W22);
\draw[->] (W15) -- (W23);
\draw[->] (W15) -- (W24);
\draw[->] (W15) -- (W25);
\draw[->,\XColor,ultra thick] (W00) -- (W12);
\draw[->,\XColor,ultra thick] (W00) -- (W13);
\draw[->,\XColor,ultra thick] (W01) -- (W14);
\draw[->,\XColor,ultra thick] (W01) -- (W15);
\draw[->,\XColor,ultra thick] (W10) -- (W22);
\draw[->,\XColor,ultra thick] (W10) -- (W23);
\draw[->,\XColor,ultra thick] (W11) -- (W24);
\draw[->,\XColor,ultra thick] (W11) -- (W25);
}
\]
For $(x,y)\in \cG\circ \cD(Q_2\to Q_3)_n = \bbZ M_{2-3}^\triv\oplus \bbZ M_{2-3}^\sign$,
$\widehat{G}(M_{3-2}^\triv)(x,y) = (x,x,y,y)\in \cG\circ \cD(Q_2\to Q_2)_n=\bbZ M_{2-2,0}^\triv \oplus\bbZ M_{2-2,1}^\triv\oplus\bbZ M_{2-2,0}^\sign\oplus\bbZ M_{2-2,1}^\sign$.
Recall the isomorphism $\phi_1$ and $\phi_2$ to $\bbZ\!\left[\frac{1}{4}\right]$,
we can see 
$$\phi_2\circ \widehat{G}(M_{3-2}^\triv)(x,y) = \phi_2 (x,x,y,y) = 4^{-n}\cdot\frac{2x+2y}{4} = \phi_1(x,y).$$
Therefore, $\widehat{G}(M_{3-2}^\triv)$ is induced by an identity map.
Similarly, $\widehat{G}(M_{3-2}^\sign)$ is induced by an identity map.
\item For $M_{2-3}^{\triv/\sign}$: Similar to the case $M_{2-1,0}$, $\widehat{G}(M_{2-3}^\triv)$ and $\widehat{G}(M_{2-3}^\sign)$ are induced by maps of multiplication by 2.
\end{itemize}
\medskip

Next, we compute the invariant $\widehat{F}$.
\begin{itemize}
\item For $Q=Q_1$: There is only one vertex $Q_1$ in each level,
so $\widehat{F}(Q_1) = \varinjlim \cG\circ  \cD(Q_1\to Q_1)$, 
where $\cG\circ \cD(Q_1\to Q_1)_n$ has $\bbZ$-basis $\{M_{1-1,g}\}_{g\in\bbZ/4\bbZ}$. 

Note that the connecting map is given by $(\bigoplus_{g\in\bbZ/4\bbZ} M_{1-1,g})\boxtimes_{Q_1} - $,
we have the following Bratteli diagram:
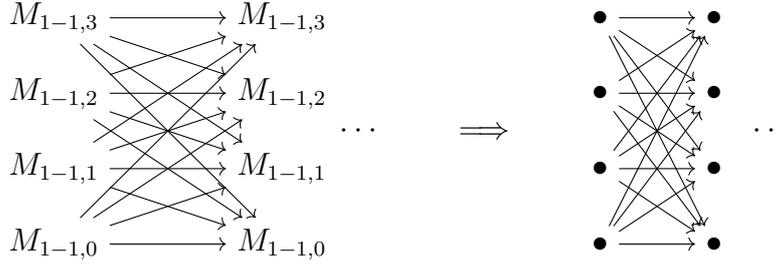
\[
\tikzmath{
\node (W00) at (0,0) {$M_{1-1,0}$};
\node (W01) at (0,1) {$M_{1-1,1}$};
\node (W02) at (0,2) {$M_{1-1,2}$};
\node (W03) at (0,3) {$M_{1-1,3}$};
\node (W10) at (3,0) {$M_{1-1,0}$};
\node (W11) at (3,1) {$M_{1-1,1}$};
\node (W12) at (3,2) {$M_{1-1,2}$};
\node (W13) at (3,3) {$M_{1-1,3}$};
\draw[->] (W00) -- (W10);
\draw[->] (W00) -- (W11);
\draw[->] (W00) -- (W12);
\draw[->] (W00) -- (W13);
\draw[->] (W01) -- (W10);
\draw[->] (W01) -- (W11);
\draw[->] (W01) -- (W12);
\draw[->] (W01) -- (W13);
\draw[->] (W02) -- (W10);
\draw[->] (W02) -- (W11);
\draw[->] (W02) -- (W12);
\draw[->] (W02) -- (W13);
\draw[->] (W03) -- (W10);
\draw[->] (W03) -- (W11);
\draw[->] (W03) -- (W12);
\draw[->] (W03) -- (W13);
\node at (4,1.5) {$\cdots$};
}
\qquad\Longrightarrow\qquad
\tikzmath{
\node (U11) at (0,0) {$\bullet$};
\node (U12) at (0,1) {$\bullet$};
\node (U13) at (0,2) {$\bullet$};
\node (U14) at (0,3) {$\bullet$};
\node (U21) at (1.5,0) {$\bullet$};
\node (U22) at (1.5,1) {$\bullet$};
\node (U23) at (1.5,2) {$\bullet$};
\node (U24) at (1.5,3) {$\bullet$};
\draw[->] (U11) -- (U21);
\draw[->] (U11) -- (U22);
\draw[->] (U11) -- (U23);
\draw[->] (U11) -- (U24);
\draw[->] (U12) -- (U21);
\draw[->] (U12) -- (U22);
\draw[->] (U12) -- (U23);
\draw[->] (U12) -- (U24);
\draw[->] (U13) -- (U21);
\draw[->] (U13) -- (U22);
\draw[->] (U13) -- (U23);
\draw[->] (U13) -- (U24);
\draw[->] (U14) -- (U21);
\draw[->] (U14) -- (U22);
\draw[->] (U14) -- (U23);
\draw[->] (U14) -- (U24);
\node at (2.25,1.5) {$\cdots$};
}
\]
Similarly, we obtain a group isomorphism
\begin{align*}
    \psi_1 : \widehat{F}(Q_1) &\to \bbZ\!\left[\textstyle\frac{1}{4}\right] \\
    (x,y,z,w)&\mapsto 4^{-n}\cdot\frac{x+y+z+w}{4}
\end{align*}
Here, $(x,y,z,w)\in \cG\circ \cD(Q_1\to Q_1)_n = \bbZ M_{1-1,0}\oplus\bbZ M_{1-1,1}\oplus \bbZ M_{1-1,2}\oplus \bbZ M_{1-1,3}$ and $n$ corresponds to the level at which these first appear.

\item For $Q=Q_2$: Similarly, 
$\widehat{F}(Q_2) = \varinjlim \cG\circ  \cD(Q_1\to Q_2)$, 
where $\cG\circ \cD(Q_1\to Q_2)_n$ has $\bbZ$-basis $M_{1-2,0}$ and $M_{1-2,1}$. 
Note that the connecting map is given by $(\bigoplus_{g\in\bbZ/4\bbZ} M_{1-1,g})\boxtimes_{Q_1} - $,
we have the following Bratteli diagram:
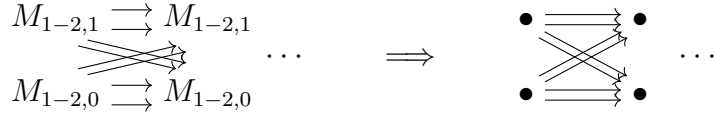
\[
\tikzmath{
\node (W00) at (0,0) {$M_{1-2,0}$};
\node (W01) at (0,1) {$M_{1-2,1}$};
\node (W10) at (2,0) {$M_{1-2,0}$};
\node (W11) at (2,1) {$M_{1-2,1}$};
\draw[->] (W00.10) -- (W10.170);
\draw[->] (W00.350) -- (W10.190);
\draw[->] (W01.10) -- (W11.170);
\draw[->] (W01.350) -- (W11.190);
\draw[postaction={transform canvas={yshift=-.77mm,xshift=1.27mm},draw}]
[->] (W00.50) to (W11.220);
\draw[postaction={transform canvas={yshift=.77mm,xshift=1.27mm},draw}]
[->] (W01.310) to (W10.140);
\node at (3,.5) {$\cdots$};
}
\qquad \Longrightarrow \qquad
\tikzmath{
\node (W11) at (0,0) {$\bullet$};
\node (W12) at (0,1) {$\bullet$};
\node (W21) at (1.5,0) {$\bullet$};
\node (W22) at (1.5,1) {$\bullet$};
\draw[->] (W11.15) -- (W21.165);
\draw[->] (W11.345) -- (W21.195);
\draw[->] (W12.15) -- (W22.165);
\draw[->] (W12.345) -- (W22.195);
\draw[postaction={transform canvas={yshift=.77mm,xshift=-.77mm},draw}]
[->] (W11.35) to (W22.235);
\draw[postaction={transform canvas={yshift=.77mm,xshift=.77mm},draw}]
[->] (W12.305) to (W21.145);
\node at (2.25,.5) {$\cdots$};
}
\]
Similarly, we obtain a group isomorphism 
\begin{align*}
    \psi_2 : \widehat{F}(Q_2) &\to \bbZ\!\left[\textstyle\frac{1}{4}\right] \\
    (x,y)&\mapsto 4^{-n}\cdot \frac{x+y}{2}.
\end{align*}
Here, we take $n$ as the first positive integer for which $(x,y)\in \cG\circ\cD(Q_1\to Q_2)_n = \bbZ M_{1-2,0} \oplus \bbZ M_{1-2,1}$ appears in the inductive system.
\item For $Q=Q_3$: Similarly,
$\widehat{F}(Q_3) = \varinjlim \cG\circ  \cD(Q_1\to Q_3)$, 
where $\cG\circ \cD(Q_1\to Q_3)_n$ has $\bbZ$-basis $M_{1-3}$. 
We have the following Bratteli diagram:
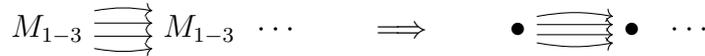
\[
\tikzmath{
\node (V1) at (0,0) {$M_{1-3}$};
\node (V2) at (2,0) {$M_{1-3}$};
\draw[->] (V1.08) -- (V2.172);
\draw[->] (V1.352) -- (V2.188);
\draw[->] (V1.22) to [bend left=10] (V2.158);
\draw[->] (V1.338) to [bend right=10] (V2.202);
\node at (3,0) {$\cdots$};
}
\qquad \Longrightarrow \qquad
\tikzmath{
\node (V1) at (0,0) {$\bullet$};
\node (V2) at (1.5,0) {$\bullet$};
\draw[->] (V1.15) -- (V2.165);
\draw[->] (V1.345) -- (V2.195);
\draw[->] (V1.35) to [bend left=10] (V2.145);
\draw[->] (V1.325) to [bend right=10] (V2.215);
\node at (2.25,0) {$\cdots$};
}
\]
we obtain a group isomorphism 
\begin{align*}
    \psi_3 : \widehat{F}(Q_3) &\to \bbZ\!\left[\textstyle\frac{1}{4}\right] \\
    x &\mapsto 4^{-n}\cdot x.
\end{align*}
Here, we take $n$ as the first positive integer for which $x\in \cG\circ\cD(Q_1\to Q_3)_n = \bbZ M_{1-3}$ appears in the inductive system.
\end{itemize}

\hspace{.5cm}

Now we compute $\widehat{F}$ of morphisms in $\cD$.
\begin{itemize}
\item For $M_{1-2,0}$: 
When picking up the $\bbZ$-basis $\{M_{1-1,g}\}_{g\in \bbZ/4\bbZ}$ for $\cG\circ \cD(Q_1\to Q_1)_n$ and $\bbZ$-basis $M_{1-2,0}$, $M_{1-2,1}$ for $\cG\circ \cD(Q_1\to Q_2)_n$, we have
\begin{align*}
    L^{M_{1-2,0}}(M_{1-1,0/2}) & = M_{1-1,0/2}\boxtimes_{Q_1} M_{1-2,0}= M_{1-2,0} \\
    L^{M_{1-2,0}}(M_{1-1,1/3}) & = M_{1-1,1/3}\boxtimes_{Q_1} M_{1-2,0}= M_{1-2,1}
\end{align*}
We obtain the following diagram:
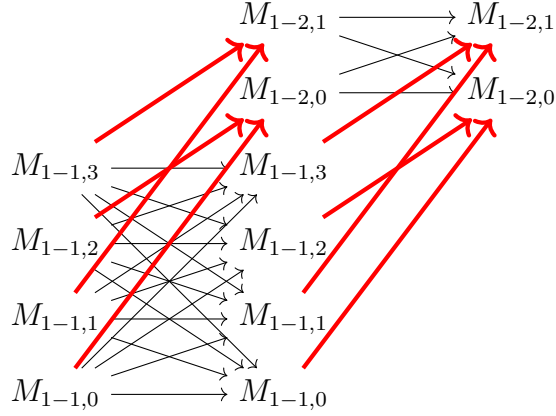
\[
\tikzmath{
\node (V00) at (0,0) {$M_{1-1,0}$};
\node (V01) at (0,1) {$M_{1-1,1}$};
\node (V02) at (0,2) {$M_{1-1,2}$};
\node (V03) at (0,3) {$M_{1-1,3}$};
\node (V10) at (3,0) {$M_{1-1,0}$};
\node (V11) at (3,1) {$M_{1-1,1}$};
\node (V12) at (3,2) {$M_{1-1,2}$};
\node (V13) at (3,3) {$M_{1-1,3}$};
\node (V14) at (3,4) {$M_{1-2,0}$};
\node (V15) at (3,5) {$M_{1-2,1}$};
\node (V24) at (6,4) {$M_{1-2,0}$};
\node (V25) at (6,5) {$M_{1-2,1}$};
\draw[->] (V00) -- (V10);
\draw[->] (V00) -- (V11);
\draw[->] (V00) -- (V12);
\draw[->] (V00) -- (V13);
\draw[->] (V01) -- (V10);
\draw[->] (V01) -- (V11);
\draw[->] (V01) -- (V12);
\draw[->] (V01) -- (V13);
\draw[->] (V02) -- (V10);
\draw[->] (V02) -- (V11);
\draw[->] (V02) -- (V12);
\draw[->] (V02) -- (V13);
\draw[->] (V03) -- (V10);
\draw[->] (V03) -- (V11);
\draw[->] (V03) -- (V12);
\draw[->] (V03) -- (V13);
\draw[->] (V14) -- (V24);
\draw[->] (V14) -- (V25);
\draw[->] (V15) -- (V24);
\draw[->] (V15) -- (V25);
\draw[->,\XColor,ultra thick] (V00) -- (V14);
\draw[->,\XColor,ultra thick] (V01) -- (V15);
\draw[->,\XColor,ultra thick] (V02) -- (V14);
\draw[->,\XColor,ultra thick] (V03) -- (V15);
\draw[->,\XColor,ultra thick] (V10) -- (V24);
\draw[->,\XColor,ultra thick] (V11) -- (V25);
\draw[->,\XColor,ultra thick] (V12) -- (V24);
\draw[->,\XColor,ultra thick] (V13) -- (V25);
}
\]
For $(x,y,z,w)\in \cG\circ \cD(Q_1\to Q_1)_n = \bbZ M_{1-1,0}\oplus \bbZ M_{1-1,1}\oplus \bbZ M_{1-1,2}\oplus \bbZ M_{1-1,3}$, we have 
$\widehat{F}(M_{1-2,0})(x,y,z,w) = (x+z,y+w)\in \cG\circ \cD(Q_1\to Q_2)_n = \bbZ M_{1-2,0}\oplus \bbZ M_{1-2,1}$.
Recalling the isomorphism $\psi_1$ and $\psi_2$ to $\bbZ\!\left[\frac{1}{4}\right]$,
we can see 
$$\psi_2\circ \widehat{F}(M_{1-2,0})(x,y,z,w) = \psi_2 (x+z,y+w) = 4^{-n}\cdot\frac{x+y+z+w}{2} = 2\psi_1(x,y,z,w).$$
Therefore, $\widehat{F}(M_{1-2,0})$ is induced by a map of multiplication by 2.
Similarly, $\widehat{F}(M_{1-2,1})$ is induced by a map of multiplication by 2.
\item For $M_{2-1,0}$: 
When picking up the $\bbZ$-basis $M_{1-2,0}$, $M_{1-2,1}$ for $\cG\circ \cD(Q_1\to Q_2)_n$ and $\bbZ$-basis $\{M_{1-1,g}\}_{g\in \bbZ/4\bbZ}$ for $\cG\circ \cD(Q_1\to Q_1)_n$, we have
\begin{align*}
    L^{M_{2-1,0}}_n (M_{1-2,0}) & = M_{1-2,0}\boxtimes_{Q_2} M_{2-1,0} = M_{1-1,0}\oplus M_{1-1,2} \\
    L^{M_{2-1,0}}_n (M_{1-2,1}) & = M_{1-2,1}\boxtimes_{Q_2} M_{2-1,0}= M_{1-1,1}\oplus M_{1-1,3}
\end{align*}
we obtain the following diagram:
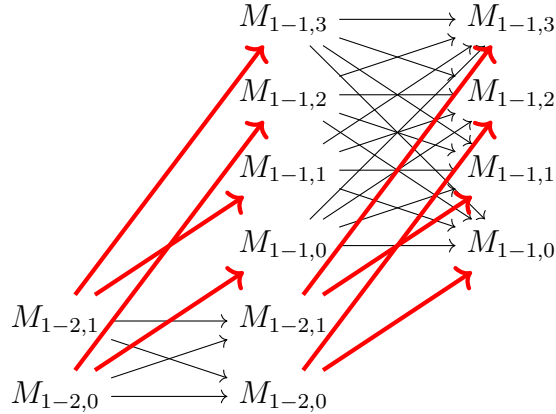
\[
\tikzmath{
\node (V00) at (0,0) {$M_{1-2,0}$};
\node (V01) at (0,1) {$M_{1-2,1}$};
\node (V10) at (3,0) {$M_{1-2,0}$};
\node (V11) at (3,1) {$M_{1-2,1}$};
\node (V12) at (3,2) {$M_{1-1,0}$};
\node (V13) at (3,3) {$M_{1-1,1}$};
\node (V14) at (3,4) {$M_{1-1,2}$};
\node (V15) at (3,5) {$M_{1-1,3}$};
\node (V22) at (6,2) {$M_{1-1,0}$};
\node (V23) at (6,3) {$M_{1-1,1}$};
\node (V24) at (6,4) {$M_{1-1,2}$};
\node (V25) at (6,5) {$M_{1-1,3}$};
\draw[->] (V00) -- (V10);
\draw[->] (V00) -- (V11);
\draw[->] (V01) -- (V10);
\draw[->] (V01) -- (V11);
\draw[->] (V12) -- (V22);
\draw[->] (V12) -- (V23);
\draw[->] (V12) -- (V24);
\draw[->] (V12) -- (V25);
\draw[->] (V13) -- (V22);
\draw[->] (V13) -- (V23);
\draw[->] (V13) -- (V24);
\draw[->] (V13) -- (V25);
\draw[->] (V14) -- (V22);
\draw[->] (V14) -- (V23);
\draw[->] (V14) -- (V24);
\draw[->] (V14) -- (V25);
\draw[->] (V15) -- (V22);
\draw[->] (V15) -- (V23);
\draw[->] (V15) -- (V24);
\draw[->] (V15) -- (V25);
\draw[->,\XColor,ultra thick] (V00) -- (V12);
\draw[->,\XColor,ultra thick] (V00) -- (V14);
\draw[->,\XColor,ultra thick] (V01) -- (V13);
\draw[->,\XColor,ultra thick] (V01) -- (V15);
\draw[->,\XColor,ultra thick] (V10) -- (V22);
\draw[->,\XColor,ultra thick] (V10) -- (V24);
\draw[->,\XColor,ultra thick] (V11) -- (V23);
\draw[->,\XColor,ultra thick] (V11) -- (V25);
}
\]
For $(x,y)\in \cG\circ \cD(Q_1\to Q_2)_n = \bbZ M_{1-2,0}\oplus \bbZ M_{1-2,1}$, we have
$\widehat{F}(M_{2-1,0})(x,y) = (x,y,x,y)\in \cG\circ \cD(Q_1\to Q_1)_n = \bbZ M_{1-1,0}\oplus \bbZ M_{1-1,1}\oplus \bbZ M_{1-1,2}\oplus \bbZ M_{1-1,3}$. 
Recalling the isomorphism $\psi_1$ and $\psi_2$ to $\bbZ\!\left[\frac{1}{4}\right]$,
we can see 
$$\psi_1\circ \widehat{F}(M_{2-1,0})(x,y) = \psi_1 (x,y,x,y) = 4^{-n}\cdot\frac{2x+2y}{4} = \psi_2(x,y).$$
Therefore, $\widehat{F}(M_{2-1,0})$ is induced by an identity map.
Similarly, $\widehat{F}(M_{2-1,1})$ is induced by an identity map.
\item For $M_{1-1,g}$, $M_{2-2,0/1}^{\triv/\sign}$ and $M_{3-3}^{\chi_k}$: Similarly, it is easy to see their $\widehat{F}$'s are induced by identity maps.
\item For $M_{1-3}$: 
When picking up the $\bbZ$-basis $\{M_{1-1,g}\}_{g\in \bbZ/4\bbZ}$ for $\cG\circ \cD(Q_1\to Q_1)_n$ and $\bbZ$-basis $M_{1-2,0}$, $M_{1-2,1}$ for $\cG\circ \cD(Q_1\to Q_2)_n$, we have
\begin{align*}
    L^{M_{1-3}}_n(M_{1-1,g}) = M_{1-1,g}\boxtimes_{Q_1} M_{1-3}=M_{1-3},\qquad g\in\bbZ/4\bbZ,
\end{align*}
we obtain the following diagram:
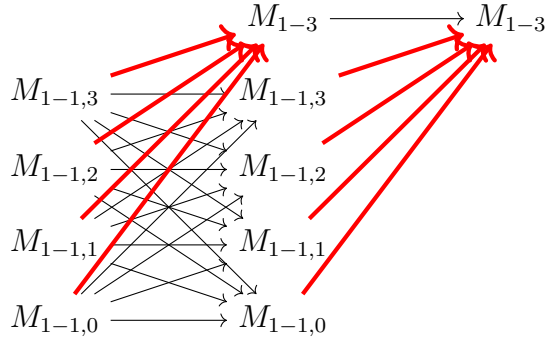
\[
\tikzmath{
\node (V00) at (0,0) {$M_{1-1,0}$};
\node (V01) at (0,1) {$M_{1-1,1}$};
\node (V02) at (0,2) {$M_{1-1,2}$};
\node (V03) at (0,3) {$M_{1-1,3}$};
\node (V10) at (3,0) {$M_{1-1,0}$};
\node (V11) at (3,1) {$M_{1-1,1}$};
\node (V12) at (3,2) {$M_{1-1,2}$};
\node (V13) at (3,3) {$M_{1-1,3}$};
\node (V14) at (3,4) {$M_{1-3}$};
\node (V24) at (6,4) {$M_{1-3}$};
\draw[->] (V00) -- (V10);
\draw[->] (V00) -- (V11);
\draw[->] (V00) -- (V12);
\draw[->] (V00) -- (V13);
\draw[->] (V01) -- (V10);
\draw[->] (V01) -- (V11);
\draw[->] (V01) -- (V12);
\draw[->] (V01) -- (V13);
\draw[->] (V02) -- (V10);
\draw[->] (V02) -- (V11);
\draw[->] (V02) -- (V12);
\draw[->] (V02) -- (V13);
\draw[->] (V03) -- (V10);
\draw[->] (V03) -- (V11);
\draw[->] (V03) -- (V12);
\draw[->] (V03) -- (V13);
\draw[->] (V14) -- (V24);
\draw[->,\XColor,ultra thick] (V00) -- (V14);
\draw[->,\XColor,ultra thick] (V01) -- (V14);
\draw[->,\XColor,ultra thick] (V02) -- (V14);
\draw[->,\XColor,ultra thick] (V03) -- (V14);
\draw[->,\XColor,ultra thick] (V10) -- (V24);
\draw[->,\XColor,ultra thick] (V11) -- (V24);
\draw[->,\XColor,ultra thick] (V12) -- (V24);
\draw[->,\XColor,ultra thick] (V13) -- (V24);
}
\]
For $(x,y,z,w)\in \cG\circ \cD(Q_1\to Q_1)_n \bbZ M_{1-1,0}\oplus \bbZ M_{1-1,1}\oplus \bbZ M_{1-1,2}\oplus \bbZ M_{1-1,3}$,
$\widehat{F}(M_{1-3})(x,y,z,w) = x+y+z+w\in \cG\circ \cD(Q_1\to Q_3)_n = \bbZ M_{1-3}$.
Recalling the isomorphism $\psi_1$ and $\psi_3$ to $\bbZ\!\left[\frac{1}{4}\right]$,
we can see 
$$\psi_3\circ \widehat{F}(M_{1-3})(x,y,z,w) = \psi_3 (x+y+z+w) = 4^{-n}\cdot(x+y+z+w) = 4\psi_1(x,y,z,w).$$
Therefore, $\widehat{F}(M_{1-3})$ is induced by a map of multiplication by 4.
\item For $M_{3-1}$: 
When picking up the $\bbZ$-basis $M_{1-3}$ for $\cG\circ \cD(Q_1\to Q_3)_n$ and $\bbZ$-basis $\{M_{1-1,g}\}_{g\in\bbZ/4\bbZ}$ for $\cG\circ \cD(Q_1\to Q_1)_n$, we have
$$L^{M_{3-1}}_n(M_{1-3}) = M_{1-3}\boxtimes_{Q_3} M_{3-1}=\bigoplus_{g\in\bbZ/4\bbZ} M_{1-1,g},$$
we obtain the following diagram
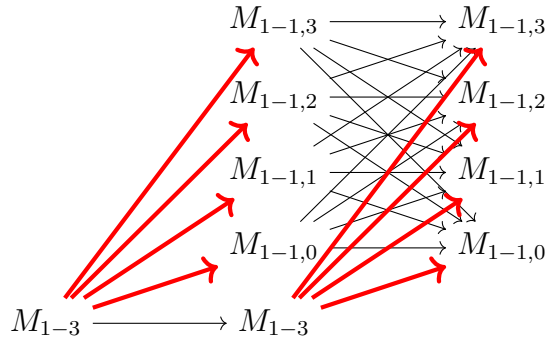
\[
\tikzmath{
\node (V00) at (0,0) {$M_{1-3}$};
\node (V10) at (3,0) {$M_{1-3}$};
\node (V12) at (3,1) {$M_{1-1,0}$};
\node (V13) at (3,2) {$M_{1-1,1}$};
\node (V14) at (3,3) {$M_{1-1,2}$};
\node (V15) at (3,4) {$M_{1-1,3}$};
\node (V22) at (6,1) {$M_{1-1,0}$};
\node (V23) at (6,2) {$M_{1-1,1}$};
\node (V24) at (6,3) {$M_{1-1,2}$};
\node (V25) at (6,4) {$M_{1-1,3}$};
\draw[->] (V00) -- (V10);
\draw[->] (V12) -- (V22);
\draw[->] (V12) -- (V23);
\draw[->] (V12) -- (V24);
\draw[->] (V12) -- (V25);
\draw[->] (V13) -- (V22);
\draw[->] (V13) -- (V23);
\draw[->] (V13) -- (V24);
\draw[->] (V13) -- (V25);
\draw[->] (V14) -- (V22);
\draw[->] (V14) -- (V23);
\draw[->] (V14) -- (V24);
\draw[->] (V14) -- (V25);
\draw[->] (V15) -- (V22);
\draw[->] (V15) -- (V23);
\draw[->] (V15) -- (V24);
\draw[->] (V15) -- (V25);
\draw[->,\XColor,ultra thick] (V00) -- (V12);
\draw[->,\XColor,ultra thick] (V00) -- (V13);
\draw[->,\XColor,ultra thick] (V00) -- (V14);
\draw[->,\XColor,ultra thick] (V00) -- (V15);
\draw[->,\XColor,ultra thick] (V10) -- (V22);
\draw[->,\XColor,ultra thick] (V10) -- (V23);
\draw[->,\XColor,ultra thick] (V10) -- (V24);
\draw[->,\XColor,ultra thick] (V10) -- (V25);
}
\]
For $x\in \cG\circ \cD(Q_1\to Q_3)_n=\bbZ M_{1-3}$,
$\widehat{F}(M_{3-1})(x) = (x,x,x,x)\in \cG\circ \cD(Q_1\to Q_1)_n = \bbZ M_{1-1,0}\oplus \bbZ M_{1-1,1}\oplus \bbZ M_{1-1,2}\oplus \bbZ M_{1-1,3}$ corresponds to $M_{1-1,g}$.
Recall the isomorphism $\psi_1$ and $\psi_3$ to $\bbZ\!\left[\frac{1}{4}\right]$,
we can see 
$$\psi_1\circ \widehat{F}(M_{3-1})(x) = \psi_1 (x,x,x,x) = 4^{-n}\cdot\frac{x+x+x+x}{4} = \psi_3(x).$$
Therefore, $\widehat{F}(M_{3-1})$ is induced by an identity map.
\item For $M_{2-3}^\triv$: 
When picking up the $\bbZ$-basis $M_{1-3}$ for $\cG\circ \cD(Q_1\to Q_3)_n$ and $\bbZ$-basis $M_{1-2,0}$ and $M_{1-2,1}$ for $\cG\circ \cD(Q_1\to Q_2)_n$, we have
$$L^{M_{2-3}^\triv}_n(M_{1-2,i}) = M_{1-2,i}\boxtimes_{Q_2} M_{2-3}^\triv = M_{1-3}, \qquad i=0,1$$
We then obtain the following diagram:
\[
\tikzmath{
\node (W00) at (0,0) {$M_{1-2,0}$};
\node (W01) at (0,1) {$M_{1-2,1}$};
\node (W10) at (3,0) {$M_{1-2,0}$};
\node (W11) at (3,1) {$M_{1-2,1}$};
\node (W12) at (3,2) {$M_{1-3}$};
\node (W22) at (6,2) {$M_{1-3}$};
\draw[->] (W00) -- (W10);
\draw[->] (W00) -- (W11);
\draw[->] (W01) -- (W10);
\draw[->] (W01) -- (W11);
\draw[->] (W12) -- (W22);
\draw[->,\XColor,ultra thick] (W00) -- (W12);
\draw[->,\XColor,ultra thick] (W01) -- (W12);
\draw[->,\XColor,ultra thick] (W10) -- (W22);
\draw[->,\XColor,ultra thick] (W11) -- (W22);
}
\]
For $(x,y)\in \cG\circ \cD(Q_1\to Q_2)_n =\bbZ M_{1-2,0}\oplus \bbZ M_{1-2,1}$,
$\widehat{F}(M_{2-3}^\triv)(x,y) = x+y\in \cG\circ \cD(Q_1\to Q_3)_n = \bbZ M_{1-3}$. 
Recalling the isomorphism $\psi_2$ and $\psi_3$ to $\bbZ\!\left[\frac{1}{4}\right]$,
we can see 
$$\psi_3\circ \widehat{F}(M_{2-3}^\triv)(x,y) = \psi_3 (x+y) = 4^{-n}\cdot(x+y) = 2\psi_2(x,y).$$
Therefore, $\widehat{F}(M_{2-3}^\triv)$ is induced by a map of multiplication by 2.
Similarly, $\widehat{F}(M_{2-3}^\sign)$ is induced by a map of multiplication by 2.
\item For $M_{3-2}^\triv$:
When picking up the $\bbZ$-basis $M_{1-3}$ for $\cG\circ \cD(Q_1\to Q_3)_n$ and $\bbZ$-basis $M_{1-2,0}$ and $M_{1-2,1}$ for $\cG\circ \cD(Q_1\to Q_2)_n$, we have
$$L^{M_{3-2}^\triv}_n(M_{1-3}) = M_{1-3}\boxtimes_{Q_3} M_{3-2}^\triv = M_{1-2,0}\oplus M_{1-2,1},$$
we then obtain the following diagram:
\[
\tikzmath{
\node (V00) at (0,0) {$M_{1-3}$};
\node (V10) at (3,0) {$M_{1-3}$};
\node (V12) at (3,1) {$M_{1-2,0}$};
\node (V13) at (3,2) {$M_{1-2,1}$};
\node (V22) at (6,1) {$M_{1-2,0}$};
\node (V23) at (6,2) {$M_{1-2,1}$};
\draw[->] (V00) -- (V10);
\draw[->] (V12) -- (V22);
\draw[->] (V12) -- (V23);
\draw[->] (V13) -- (V22);
\draw[->] (V13) -- (V23);
\draw[->,\XColor,ultra thick] (V00) -- (V12);
\draw[->,\XColor,ultra thick] (V00) -- (V13);
\draw[->,\XColor,ultra thick] (V10) -- (V22);
\draw[->,\XColor,ultra thick] (V10) -- (V23);
}
\]
For $x\in \cG\circ \cD(Q_1\to Q_3) = \bbZ M_{1-3}$,
$\widehat{F}(M_{3-2}^\triv)(x) = (x,x)\in \cG\circ \cD(Q_1\to Q_2) = \bbZ M_{1-2,0}\oplus \bbZ M_{1-2,1}$.
Recall the isomorphism $\psi_2$ and $\psi_3$ to $\bbZ\!\left[\frac{1}{4}\right]$,
we can see 
$$\psi_2\circ \widehat{F}(M_{3-2}^\triv)(x) = \psi_2 (x,x) = 4^{-n}\cdot\frac{x+x}{2} = \psi_3(x).$$
Therefore, $\widehat{F}(M_{3-2}^\triv)$ is induced by an identity map.
Similarly, $\widehat{F}(M_{3-2}^\sign)$ is induced by an identity map.
\end{itemize}
\medskip

For convenience, we record the values of $\widehat{F}$ and $\widehat{G}$ on morphisms in the following table, where each entry denotes the corresponding multiplication map $\bbZ\!\left[\frac{1}{4}\right]\to\bbZ\!\left[\frac{1}{4}\right]$:
\[
\bgroup
\def\arraystretch{1.5}%
\begin{tabular}{c c c c c c c c c c}
\hline
& $M_{1-1,g}$ & $M_{1-2,i}$ & $M_{1-3}$ & $M_{2-1,i}$ & $M_{2-2,i}^{\triv/\sign}$ & $M_{2-3}^{\triv/\sign}$ & $M_{3-1}$ & $M_{3-2}^{\triv/\sign}$ & $M_{3-3}^{\chi_k}$\\
\hline 
$\widehat{F}$ & $1$ & $2$ & $4$ & $1$ & $1$ & $2$ & $1$ & $1$ & $1$ \\
\hline
$\widehat{G}$ & $1$ & $1$ & $2$ & $2$ & $1$ & $2$ & $2$ & $1$ & $1$\\
\hline
\end{tabular}
\egroup
\]

We now construct an explicit natural isomorphism $\alpha:\widehat{F}\Rightarrow \widehat{G}$.
Let $\alpha_1: \widehat{F}(Q_1)\to \widehat{G}(Q_1)$ be induced by a map of multiplication by 2, and $\alpha_2: \widehat{F}(Q_2)\to \widehat{G}(Q_2)$ and $\alpha_3: \widehat{F}(Q_3)\to \widehat{G}(Q_3)$ be induced by identity maps.
Naturality follows from direct computations. 
One can easily verify that 
\[
\begin{tikzcd}[column sep=4em]
\widehat{F}(P)
\arrow[swap]{d}{\alpha_P}
\arrow{r}{\widehat{F}(M_{P-Q})}
& \widehat{F}(Q)
\arrow{d}{\alpha_Q}
 \\
\widehat{G}(P)
\arrow{r}{\widehat{G}(M_{P-Q})}
& \widehat{G}(Q)
\end{tikzcd}
\]
is commutative for all objects $P,Q\in[\QSys(\Hilb(\bbZ/4\bbZ))]$ and morphisms $M_{P-Q}$.
Therefore, $\Hilb(\bbZ/4\bbZ)\overset{F}{\curvearrowright} A$ and $\Hilb(\bbZ/4\bbZ)\overset{G}{\curvearrowright} A$ are equivalent AF-actions by Theorem \ref{thmalp:main}.

Verifying that 
$\Hilb(\bbZ/4\bbZ)\overset{H}{\curvearrowright} A$ and $\Hilb(\bbZ/4\bbZ)\overset{G}{\curvearrowright} A$ are equivalent AF-actions can be done by similar methods. 

\medskip

Finally, we provide an example of a $\Hilb(\bbZ/4\bbZ)$ AF-action $\Hilb(\bbZ/4\bbZ)\overset{E}{\curvearrowright} A$ on $A\cong\bbM_{4^\infty}$ which is not equivalent to any of $F, G$ or $H$ from the previous examples. To construct $E$ we consider the following homogeneous enriched Bratteli diagram: 
\[
\tikzmath{
\node (V1) at (-1.5,0) {$Q_3$};
\node (V2) at (1.5,0) {$Q_3$};
\draw[->] (V1) -- (V2);
\node at (0,.5) {${M_{3-3}^\triv}^{\oplus 4}$};
}.
\]
We shall now compute $\widehat{E}$.

\begin{itemize}
\item For $Q=Q_1$: 
$\widehat{E}(Q_1) = \varinjlim \cG\circ  \cD(Q_3\to Q_1)$, 
where $\cG\circ \cD(Q_3\to Q_1)_n$ has $\bbZ$-basis $M_{3-1}$. 
Note that $M_{3-3}^\triv\boxtimes_{Q_3} M_{3-1} = M_{3-1}$, 
so we have the following Bratteli diagram:
\[
\tikzmath{
\node (V1) at (0,0) {$M_{3-1}$};
\node (V2) at (2,0) {$M_{3-1}$};
\draw[->] (V1.08) -- (V2.172);
\draw[->] (V1.352) -- (V2.188);
\draw[->] (V1.22) to [bend left=10] (V2.158);
\draw[->] (V1.338) to [bend right=10] (V2.202);
\node at (3,0) {$\cdots$};
}
\qquad \Longrightarrow \qquad
\tikzmath{
\node (V1) at (0,0) {$\bullet$};
\node (V2) at (1.5,0) {$\bullet$};
\draw[->] (V1.15) -- (V2.165);
\draw[->] (V1.345) -- (V2.195);
\draw[->] (V1.35) to [bend left=10] (V2.145);
\draw[->] (V1.325) to [bend right=10] (V2.215);
\node at (2.25,0) {$\cdots$};
}.
\]
Computing $K_0$ for this Bratteli diagram yields $\widehat{E}(Q_1)\cong \bbZ\!\left[\frac{1}{4}\right]$.
\item For $Q=Q_2$: 
Similarly, $\widehat{E}(Q_2) = \varinjlim \cG\circ  \cD(Q_3\to Q_2)$, 
where $\cG\circ \cD(Q_3\to Q_2)_n$ has $\bbZ$-basis $M_{3-2}^\triv$ and $M_{3-2}^\sign$. 
Note that $M_{3-3}^\triv\boxtimes_{Q_3} M_{3-2}^{\triv} = M_{3-2}^\triv$ and $M_{3-3}^\triv\boxtimes_{Q_3} M_{3-2}^{\sign} = M_{3-2}^\sign$.
We have the following Bratteli diagram:
\[
\tikzmath{
\node (V10) at (0,0) {$M_{3-2}^\triv$};
\node (V11) at (0,1) {$M_{3-2}^\sign$};
\node (V20) at (3,0) {$M_{3-2}^\triv$};
\node (V21) at (3,1) {$M_{3-2}^\sign$};
\draw[->] (V10.08) -- (V20.172);
\draw[->] (V10.352) -- (V20.188);
\draw[->] (V10.22) to [bend left=10] (V20.158);
\draw[->] (V10.338) to [bend right=10] (V20.202);
\draw[->] (V11.08) -- (V21.172);
\draw[->] (V11.352) -- (V21.188);
\draw[->] (V11.22) to [bend left=10] (V21.158);
\draw[->] (V11.338) to [bend right=10] (V21.202);
\node at (4,.5) {$\cdots$};
}
\qquad \Longrightarrow \qquad
\tikzmath{
\node (V10) at (0,0) {$\bullet$};
\node (V11) at (0,1) {$\bullet$};
\node (V20) at (1.5,0) {$\bullet$};
\node (V21) at (1.5,1) {$\bullet$};
\draw[->] (V10.15) -- (V20.165);
\draw[->] (V10.345) -- (V20.195);
\draw[->] (V10.35) to [bend left=10] (V20.145);
\draw[->] (V10.325) to [bend right=10] (V20.215);
\draw[->] (V11.15) -- (V21.165);
\draw[->] (V11.345) -- (V21.195);
\draw[->] (V11.35) to [bend left=10] (V21.145);
\draw[->] (V11.325) to [bend right=10] (V21.215);
\node at (2.25,.5) {$\cdots$};
}
\]
Computing $K_0$ for this Bratteli diagram yields $\widehat{E}(Q_2)\cong \bbZ\!\left[\frac{1}{4}\right]^{\oplus 2}$.
\item For $Q=Q_3$: Similarly,
$\widehat{E}(Q_3) = \varinjlim \cG\circ  \cD(Q_3\to Q_3)$, 
where $\cG\circ \cD(Q_3\to Q_3)_n$ has $\bbZ$-basis $M_{3-3}^{\chi_k}$, $k=0,1,2,3$. 
Note that $M_{3-3}^\triv\boxtimes_{Q_3} M_{3-3}^{\chi_k} = M_{3-3}^{\chi_k}$.
We have the following Bratteli diagram:
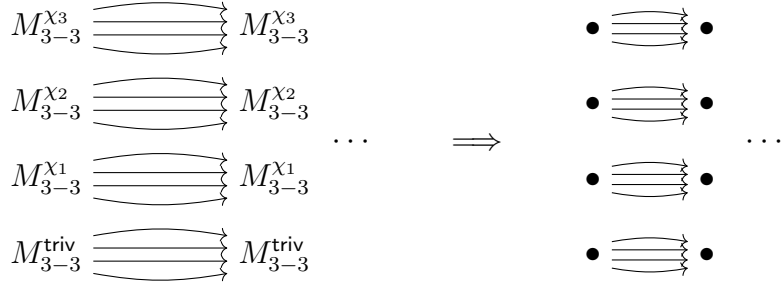
\[
\tikzmath{
\node (V10) at (0,0) {$M_{3-3}^\triv$};
\node (V11) at (0,1) {$M_{3-3}^{\chi_1}$};
\node (V12) at (0,2) {$M_{3-3}^{\chi_2}$};
\node (V13) at (0,3) {$M_{3-3}^{\chi_3}$};
\node (V20) at (3,0) {$M_{3-3}^\triv$};
\node (V21) at (3,1) {$M_{3-3}^{\chi_1}$};
\node (V22) at (3,2) {$M_{3-3}^{\chi_2}$};
\node (V23) at (3,3) {$M_{3-3}^{\chi_3}$};
\draw[->] (V10.08) -- (V20.172);
\draw[->] (V10.352) -- (V20.188);
\draw[->] (V10.22) to [bend left=10] (V20.158);
\draw[->] (V10.338) to [bend right=10] (V20.202);
\draw[->] (V11.08) -- (V21.172);
\draw[->] (V11.352) -- (V21.188);
\draw[->] (V11.22) to [bend left=10] (V21.158);
\draw[->] (V11.338) to [bend right=10] (V21.202);
\draw[->] (V12.08) -- (V22.172);
\draw[->] (V12.352) -- (V22.188);
\draw[->] (V12.22) to [bend left=10] (V22.158);
\draw[->] (V12.338) to [bend right=10] (V22.202);
\draw[->] (V13.08) -- (V23.172);
\draw[->] (V13.352) -- (V23.188);
\draw[->] (V13.22) to [bend left=10] (V23.158);
\draw[->] (V13.338) to [bend right=10] (V23.202);
\node at (4,1.5) {$\cdots$};
}
\qquad \Longrightarrow \qquad
\tikzmath{
\node (V10) at (0,0) {$\bullet$};
\node (V11) at (0,1) {$\bullet$};
\node (V12) at (0,2) {$\bullet$};
\node (V13) at (0,3) {$\bullet$};
\node (V20) at (1.5,0) {$\bullet$};
\node (V21) at (1.5,1) {$\bullet$};
\node (V22) at (1.5,2) {$\bullet$};
\node (V23) at (1.5,3) {$\bullet$};
\draw[->] (V10.15) -- (V20.165);
\draw[->] (V10.345) -- (V20.195);
\draw[->] (V10.35) to [bend left=10] (V20.145);
\draw[->] (V10.325) to [bend right=10] (V20.215);
\draw[->] (V11.15) -- (V21.165);
\draw[->] (V11.345) -- (V21.195);
\draw[->] (V11.35) to [bend left=10] (V21.145);
\draw[->] (V11.325) to [bend right=10] (V21.215);
\draw[->] (V12.15) -- (V22.165);
\draw[->] (V12.345) -- (V22.195);
\draw[->] (V12.35) to [bend left=10] (V22.145);
\draw[->] (V12.325) to [bend right=10] (V22.215);
\draw[->] (V13.15) -- (V23.165);
\draw[->] (V13.345) -- (V23.195);
\draw[->] (V13.35) to [bend left=10] (V23.145);
\draw[->] (V13.325) to [bend right=10] (V23.215);
\node at (2.25,1.5) {$\cdots$};
}
\]
Computing $K_0$ for this Bratteli diagram yields $\widehat{E}(Q_3)\cong \bbZ\!\left[\frac{1}{4}\right]^{\oplus 4}$.
\end{itemize}

Observe that $\widehat{E}(Q_2)\cong \bbZ\!\left[\frac{1}{4}\right]^{\oplus 2}\ncong\bbZ\!\left[\frac{1}{4}\right]\cong \widehat{F}(Q_2)$ and $\widehat{E}(Q_3)\cong \bbZ\!\left[\frac{1}{4}\right]^{\oplus 4}\ncong\bbZ\!\left[\frac{1}{4}\right]\cong \widehat{F}(Q_3)$,
so AF-actions $\Hilb(\bbZ/4\bbZ)\overset{E}{\curvearrowright} A$ and $\Hilb(\bbZ/4\bbZ)\overset{F}{\curvearrowright} A$ are not equivalent. Also note that $E$ is the trivial action, and thus the $F,G,H$ are non-trivial.

\appendix
\section{$[\QSys(\Hilb(\bbZ/4\bbZ))]$}
\label{App:Examples}
\label{App:Z4}
We will describe all the indecomposable objects and indecomposable morphisms in
$[\QSys(\Hilb(\bbZ/4\bbZ))]$ together with its composition law. 
There are only three objects in $[\QSys(\Hilb(\bbZ/4\bbZ))]$, which we denote $Q_1, Q_2, Q_3$. 
These correspond to the Q-systems built from subgroups $\{0\}$, $\{0,2\}\cong\bbZ/2\bbZ$, and $\bbZ/4\bbZ$, respectively. 

Relying on Proposition \ref{prop:idcpHKbimod}, We shall now describe the isomorphism classes of indecomposable bimodules. 
Given $Q_i,Q_j\in[\QSys(\Hilb(G))],$ a $Q_i$-$Q_j$ bimodule is denoted by $M_{i-j,k}^{\ell},$ where $k$ ranges over a set of representatives of the double cosets determined by $Q_i$ and $Q_j$, and $\ell$ determines a representation of the corresponding stabilizer subgroup. (Typically, $\ell = \triv, \sign, \hdots.)$ We will omit the super-index $\ell$ where the only possible representation is trivial, or when it is clear from context. 
We then have that: 
\begin{itemize}
\item $1-1$ bimodules: $M_{1-1,g}=V_g$, where $V_g=\bbC$ for $g\in \bbZ/4\bbZ$.  ($\bbC = \bbC^\triv$.)
\item $1-2$ bimodules: $M_{1-2,0} = V_0\oplus V_2$ and $M_{1-2,1} = V_1\oplus V_3$, where $V_g = \bbC$ for $g\in\bbZ/4\bbZ$.
\item $1-3$ bimodules: $M_{1-3} = \bigoplus_{g\in \bbZ/4\bbZ} V_g$, where $V_g = \bbC$ for $g\in\bbZ/4\bbZ$.
\item $2-1$ bimodules: $M_{2-1,0} = V_0\oplus V_2$ and $M_{2-1,1} = V_1\oplus V_3$, where $V_g = \bbC$ for $g\in\bbZ/4\bbZ$.
\item $2-2$ bimodules: Note that the stabilizer of 0 is $\Stab_{2\times 2^\op}(0)= \{(0,0),(2,2)\}\cong \bbZ/2\bbZ$, so there are four indecomposable bimodules.
$M_{2-2,0}^\triv = V_0^\triv\oplus V_2^\triv$, $M_{2-2,0}^\sign = V_0^\sign\oplus V_2^\sign$ and $M_{2-2,1}^\triv = V_1^\triv\oplus V_3^\triv$, $M_{2-2,1}^\sign = V_1^\sign \oplus V_3^\sign$, where $V_g^\triv = \bbC^\triv$ and $V_g^\sign = \bbC^\sign$ for $g\in\bbZ/4\bbZ$.
\item $2-3$ bimodules: $M_{2-3}^\triv = \bigoplus_{g\in\bbZ/4\bbZ}V_g^\triv$ and $M_{2-3}^\sign = \bigoplus_{g\in\bbZ/4\bbZ}V_g^\sign$, where $V_g^\triv = \bbC^\triv$ and $V_g^\sign = \bbC^\sign$ for $g\in\bbZ/4\bbZ$.
\item $3-1$ bimodules: $M_{3-1} = \bigoplus_{g\in \bbZ/4\bbZ} V_g$, where $V_g = \bbC$ for $g\in\bbZ/4\bbZ$.
\item $3-2$ bimodules: $M_{3-2}^\triv = \bigoplus_{g\in\bbZ/4\bbZ}V_g^\triv$ and $M_{3-2}^\sign = \bigoplus_{g\in\bbZ/4\bbZ}V_g^\sign$, where $V_g^\triv = \bbC^\triv$ and $V_g^\sign = \bbC^\sign$ for $g\in\bbZ/4\bbZ$.
\item $3-3$ bimodules: Note that the stabilizer of 0 is $\Stab_{3\times 3^\op}(0)= \{(g,-g)|g\in\bbZ/4\bbZ\}\cong \bbZ/4\bbZ$, and there are four irreps of $\bbZ/4\bbZ$.
$M_{3-3}^\triv = \bigoplus_{g\in \bbZ/4\bbZ} V_g^\triv$, $M_{3-3}^{\chi_k} = \bigoplus_{g\in \bbZ/4\bbZ} V_g^{\chi_k}$, $k=1,2,3$, where $V_g^\triv = \bbC^\triv$ and $V_g^{\chi_k} = \bbC^{\chi_k}$ for $g\in \bbZ/4\bbZ$. (Here, $\triv=\chi_0$)
\end{itemize}

In the following, we describe all the possible composite bimodules in $[\QSys(\Hilb(\bbZ/4\bbZ))]:$

First, $M_{?-1}\boxtimes_{Q_1} M_{1-?}:$
\[
\scalebox{.8}{
\begin{tabular}{c c c c c c c c}
\hline
& $M_{1-1,0}$ & $M_{1-1,1}$ & $M_{1-1,2}$ & $M_{1-1,3}$ & $M_{1-2,0}$ & $M_{1-2,1}$ & $M_{1-3}$ \\
\hline 
$M_{1-1,0}$ & $M_{1-1,0}$ & $M_{1-1,1}$ & $M_{1-1,2}$ & $M_{1-1,3}$ & $M_{1-2,0}$ & $M_{1-2,1}$ & $M_{1-3}$ \\
\hline
$M_{1-1,1}$ & $M_{1-1,1}$ & $M_{1-1,2}$ & $M_{1-1,3}$ & $M_{1-1,0}$ & $M_{1-2,1}$ & $M_{1-2,0}$ & $M_{1-3}$ \\
\hline
$M_{1-1,2}$ & $M_{1-1,2}$ & $M_{1-1,3}$ & $M_{1-1,0}$ & $M_{1-1,1}$ & $M_{1-2,0}$ & $M_{1-2,1}$ & $M_{1-3}$ \\
\hline
$M_{1-1,3}$ & $M_{1-1,3}$ & $M_{1-1,0}$ & $M_{1-1,1}$ & $M_{1-1,2}$ & $M_{1-2,1}$ & $M_{1-2,0}$ & $M_{1-3}$ \\
\hline
$M_{2-1,0}$ & $M_{2-1,0}$ & $M_{2-1,1}$ & $M_{2-1,0}$ & $M_{2-1,1}$ & $
M_{2-2,0}^\triv\oplus M_{2-2,0}^\sign$ & $M_{2-2,1}^\triv\oplus M_{2-2,1}^\sign$ & $M_{2-3}^\triv\oplus M_{2-3}^\sign$\\
\hline
$M_{2-1,1}$ & $M_{2-1,1}$ & $M_{2-1,0}$ & $M_{2-1,1}$ & $M_{2-1,0}$ & $M_{2-2,1}^\triv\oplus M_{2-2,1}^\sign$ & $M_{2-2,0}^\triv\oplus M_{2-2,0}^\sign$ & $M_{2-3}^\triv\oplus M_{2-3}^\sign$\\
\hline
$M_{3-1}$ & $M_{3-1}$ & $M_{3-1}$ & $M_{3-1}$ & $M_{3-1}$ & $M_{3-2}^\triv\oplus M_{3-2}^\sign$ & $M_{3-2}^\triv\oplus M_{3-2}^\sign$ & $\bigoplus_{k=0}^3 M_{3-3}^{\chi_k}$\\
\hline
\end{tabular}
}
\]
\noindent Second, $M_{?-2}\boxtimes_{Q_2} M_{2-?}:$
\[
\hspace*{-1cm}
\scalebox{.8}{
\begin{tabular}{c c c c c c c c c}
\hline
& $M_{2-1,0}$ & $M_{2-1,1}$ & $M_{2-2,0}^\triv$ & $M_{2-2,0}^\sign$ & $M_{2-2,1}^\triv$ & $M_{2-2,1}^\sign$ & $M_{2-3}^\triv$ & $M_{2-3}^\sign$ \\
\hline 
$M_{1-2,0}$ & $M_{1-1,0}\oplus M_{1-1,2}$ & $M_{1-1,1}\oplus M_{1-1,3}$ & $M_{1-2,0}$ & $M_{1-2,0}$ & $M_{1-2,1}$ & $M_{1-2,1}$ & $M_{1-3}$ & $M_{1-3}$ \\
\hline
$M_{1-2,1}$ & $M_{1-1,1}\oplus M_{1-1,3}$ & $M_{1-1,0}\oplus M_{1-1,2}$ & $M_{1-2,1}$ & $M_{1-2,1}$ & $M_{1-2,0}$ & $M_{1-2,0}$ & $M_{1-3}$ & $M_{1-3}$ \\
\hline
$M_{2-2,0}^\triv$ & $M_{2-1,0}$ & $M_{2-1,1}$ & $M_{2-2,0}^\triv$ & $M_{2-2,0}^\sign$ & $M_{2-2,1}^\triv$ & $M_{2-2,1}^\sign$ & $M_{2-3}^\triv$ & $M_{2-3}^\sign$ \\
\hline
$M_{2-2,0}^\sign$ & $M_{2-1,0}$ & $M_{2-1,1}$ & $M_{2-2,0}^\sign$ & $M_{2-2,0}^\triv$ & $M_{2-2,1}^\sign$ & $M_{2-2,1}^\triv$ & $M_{2-3}^\sign$ & $M_{2-3}^\triv$ \\
\hline
$M_{2-2,1}^\triv$ & $M_{2-1,1}$ & $M_{2-1,0}$ & $M_{2-2,1}^\triv$ & $M_{2-2,1}^\sign$ & $M_{2-2,0}^\triv$ & $M_{2-2,0}^\sign$ & $M_{2-3}^\triv$ & $M_{2-3}^\sign$ \\
\hline
$M_{2-2,1}^\sign$ & $M_{2-1,1}$ & $M_{2-1,0}$ & $M_{2-2,1}^\sign$ & $M_{2-2,1}^\triv$ & $M_{2-2,0}^\sign$ & $M_{2-2,0}^\triv$ & $M_{2-3}^\sign$ & $M_{2-3}^\triv$ \\
\hline
$M_{3-2}^\triv$ & $M_{3-1}$ & $M_{3-1}$ & $M_{3-2}^\triv$ & $M_{3-2}^\sign$ & $M_{3-2}^\triv$ & $M_{3-2}^\sign$ & $M_{3-3}^\triv\oplus M_{3-3}^{\chi_2}$ & $M_{3-3}^{\chi_1}\oplus M_{3-3}^{\chi_3}$\\
\hline
$M_{3-2}^\sign$ & $M_{3-1}$ & $M_{3-1}$ & $M_{3-2}^\sign$ & $M_{3-2}^\triv$ & $M_{3-2}^\sign$ & $M_{3-2}^\triv$ & $M_{3-3}^{\chi_1}\oplus M_{3-3}^{\chi_3}$ & $M_{3-3}^\triv\oplus M_{3-3}^{\chi_2}$\\
\hline
\end{tabular}
}
\]
\noindent Third, $M_{?-3}\boxtimes_{Q_3} M_{3-?}:$
\[
\scalebox{.8}{
\begin{tabular}{c c c c c c c c}
\hline
& $M_{3-1}$ & $M_{3-2}^\triv$ & $M_{3-2}^\sign$ & $M_{3-3}^\triv$ & $M_{3-3}^{\chi_1}$ & $M_{3-3}^{\chi_2}$ & $M_{3-3}^{\chi_3}$ \\
\hline 
$M_{1-3}$ & $\bigoplus_{g\in\bbZ/4\bbZ} M_{1-1,g}$ & $M_{1-2,0}\oplus M_{1-2,1}$ & $M_{1-2,0}\oplus M_{1-2,1}$ & $M_{1-3}$ & $M_{1-3}$ & $M_{1-3}$ & $M_{1-3}$ \\
\hline
$M_{2-3}^\triv$ & $M_{2-1,0}\oplus M_{2-1,1}$ & $M_{2-2,0}^\triv\oplus M_{2-2,1}^\triv$ & $M_{2-2,0}^\sign\oplus M_{2-2,1}^\sign$ & $M_{2-3}^\triv$ & $M_{2-3}^\sign$ & $M_{2-3}^\triv$ & $M_{2-3}^\sign$ \\
\hline
$M_{2-3}^\sign$ & $M_{2-1,0}\oplus M_{2-1,1}$ & $M_{2-2,0}^\sign\oplus M_{2-2,1}^\sign$ & $M_{2-2,0}^\triv\oplus M_{2-2,1}^\triv$ & $M_{2-3}^\sign$ & $M_{2-3}^\triv$ & $M_{2-3}^\sign$ & $M_{2-3}^\triv$ \\
\hline
$M_{3-3}^\triv$ & $M_{3-1}$ & $M_{3-2}^\triv$ & $M_{3-2}^\sign$ & $M_{3-3}^\triv$ & $M_{3-3}^{\chi_1}$ & $M_{3-3}^{\chi_2}$ & $M_{3-3}^{\chi_3}$ \\
\hline
$M_{3-3}^{\chi_1}$ & $M_{3-1}$ & $M_{3-2}^\sign$ & $M_{3-2}^\triv$ & $M_{3-3}^{\chi_1}$ & $M_{3-3}^{\chi_2}$ & $M_{3-3}^{\chi_3}$ & $M_{3-3}^\triv$ \\
\hline
$M_{3-3}^{\chi_2}$ & $M_{3-1}$ & $M_{3-2}^\triv$ & $M_{3-2}^\sign$ & $M_{3-3}^{\chi_2}$ & $M_{3-3}^{\chi_3}$ & $M_{3-3}^\triv$ & $M_{3-3}^{\chi_1}$ \\
\hline
$M_{3-3}^{\chi_3}$ & $M_{3-1}$ & $M_{3-2}^\sign$ & $M_{3-2}^\triv$ & $M_{3-3}^{\chi_3}$ & $M_{3-3}^\triv$ & $M_{3-3}^{\chi_1}$ & $M_{3-3}^{\chi_2}$ \\
\hline
\end{tabular}
}
\]

\bibliographystyle{alpha}
{\footnotesize{
\bibliography{bibliography}}}

\hspace{.5cm}

\Contact

\end{document}